\documentclass[12pt,a4paper,margin=1in]{article}

\usepackage{amsmath}
\usepackage{amssymb}
\usepackage{amsthm}
\usepackage{bm}
\usepackage{mathrsfs}

\usepackage{graphicx}
\usepackage{xcolor}
\usepackage{booktabs}
\usepackage{caption}
\usepackage{subcaption}
\usepackage{float}
\usepackage{placeins}

\usepackage{algorithm}
\usepackage{algpseudocode}
\usepackage{listings}

\usepackage{tikz}
\usepackage{pgfplots}
\usetikzlibrary{
  patterns,
  positioning,
  decorations.markings,
  matrix,
  shapes,
  arrows.meta
}

\usepackage[a4paper,margin=1in]{geometry}
\usepackage{authblk}
\usepackage{etoolbox}

\usepackage{cite}
\usepackage{hyperref}

\newcommand{\jump}[1]{[\![#1]\!]}

\DeclareMathAlphabet{\mathpzc}{OT1}{pzc}{m}{it}

\newtheorem{Theorem}{Theorem}

\newtheorem{remark}{Remark}

\title{hp-Adaptive Finite Element Modeling of Static Cracks in Density-Dependent Media}

\author[1]{Vicente Valle Martinez}
\author[1]{Saugata Ghosh}
\author[2,*]{S. M. Mallikarjunaiah}
\author[1]{Dambaru Bhatta}
\affil[1]{School of Mathematical \& Statistical Sciences,
University of Texas Rio Grande Valley, TX 78539, USA \textit{Emails:} \texttt{vicente.valle01@utrgv.edu}, \texttt{saugata.ghosh01@utrgv.edu}, \texttt{dambaru.bhatta@utrgv.edu}}
\affil[2]{Department of Mathematics \& Statistics,
Texas A\&M University--Corpus Christi, TX 78412, USA  
\textit{Email:} \texttt{M.Muddamallappa@tamucc.edu} (M. Mallikarjunaiah)\textsuperscript{*} \texttt{(Corresponding Author) }
}
\date{}

\begin{document}

\maketitle

\begin{abstract}
The elastostatic crack-tip fields within a density-dependent continuum are numerically investigated under a variety of macroscopic loading configurations. The governing nonlinear boundary value problem is formulated based on an \textit{algebraically nonlinear} hyperelastic constitutive relationship, in which the material moduli are intrinsically coupled to the volumetric strain. The existence and uniqueness of the weak solution for the nonlinear displacement-based variational formulation are rigorously established at the continuous level by the \textit{Minty-Browder theorem}. Subsequently, the variational formulation is linearized at the partial differential equation level to establish a continuous damped Newton-Raphson scheme. To capture the severe field gradients and internal singularities inherent to this model accurately, the resulting sequence of linear operator equations is discretized and solved via a conforming finite element method within an $hp$-adaptive framework. A comprehensive comparative analysis is conducted across five distinct boundary loading profiles—encompassing uniform, piecewise-slope, and concentrated point loads—applied to a stationary cracked geometry. Within each configuration, the density-dependent constitutive parameter  is systematically varied so that its regulatory impact on the fracture process zone can be assessed. Particular attention is directed toward the spatial evolution of the \textit{normal stress}, \textit{normal strain}, and \textit{strain energy density} in the immediate vicinity of the crack tip. It is demonstrated by the numerical results that local kinematic deformability is strictly bounded by the strain-limiting nature of the continuum, whereby static equilibrium is sustained by the material through severe, localized stress amplifications. Furthermore, it is revealed that the mechanical influence of the modeling parameter is fundamentally dictated by the complex topological interplay between the internal crack geometry and the spatial distribution of the applied macroscopic loads.
\end{abstract}


\section{Introduction}

The presence of a macroscopic crack within an elastic continuum profoundly alters its mechanical response, precipitating highly localized stress and strain fields in the vicinity of the crack tip. While classical linear elastic fracture mechanics provides a foundational framework for analyzing these localized fields, its reliance on a strictly linear stress--strain relation, coupled with the infinitesimal-strain kinematic assumption, inherently yields non-physical stress and strain singularities \cite{Rajagopal2010}. This mathematical paradox necessitates the adoption of broader, more sophisticated constitutive theories capable of capturing nonlinear stress--strain behavior while strictly preserving the small-displacement gradient assumption \cite{Rajagopal2011Limiting}. 

A rigorous resolution to this fundamental issue is found within the implicit constitutive framework formalized in \cite{Rajagopal2003,Rajagopal2007,Rajagopal2010}. By characterizing the mechanical response via an implicit relationship between stress and strain---rather than demanding an explicit functional dependence of one upon the other---this paradigm seamlessly accommodates profound material nonlinearity even within the regime of infinitesimal deformation \cite{RajagopalSrinivasa2007,RajagopalSaccomandi2016,Rajagopal2014,Rajagopal2018Linearization}. This approach extends to nondissipative solids defined in the current configuration, relaxing the classical requirement that stress must derive solely from a strain energy potential \cite{RajagopalSrinivasa2009}. Within this broader class, implicit relations for porous elastic media exhibiting density-dependent moduli have been established \cite{Rajagopal2021}. Such models provide a thermodynamically consistent mechanism by which volumetric deformation directly dictates the effective mechanical stiffness \cite{ItouKovtunenkoRajagopal2021,RajagopalSaccomandi2022,AlagappanArumugamRajagopal2023}. 

The numerical treatment of cracked bodies governed by these density-dependent relations has recently garnered significant attention. Computational studies have resolved crack-tip fields in porous elastic solids with preferential stiffness, confirming that material nonlinearity fundamentally alters the local stress and strain response \cite{YoonMallikarjunaiahBhatta2024}. Moreover, sophisticated $hp$-adaptive finite element frameworks have been successfully deployed for stationary cracks subject to pointwise density variations \cite{Mallikarjunaiah2026}. 

Parallel advancements have been achieved within the closely related strain-limiting subclass of implicit elasticity \cite{RajagopalWalton2011,KulvaitMalekRajagopal2013,GouMallikarjunaiahRajagopalWalton2015}. Rigorous analytical results, including existence theorems for anti-plane problems \cite{BulicekMalekRajagopalWalton2015} and local strong ellipticity conditions \cite{MaiWalton2015}, have been complemented by convergence analyses for regularized finite element approximations \cite{BonitoGiraultSuli2020}. Computationally, this framework has been applied to plane-strain fracture in anisotropic solids \cite{MallikarjunaiahWalton2015}, mixed boundary value problems involving notches and interface cracks \cite{YoonMallikarjunaiah2022}, and fully coupled thermo-mechanical edge-crack configurations \cite{YoonVasudevaMallikarjunaiah2022}. Furthermore, phase-field damage formulations have been integrated with implicit elasticity to model evolving anti-plane \cite{YoonLeeMallikarjunaiah2021} and tensile fractures \cite{LeeYoonMallikarjunaiah2022}. Recent attention has also shifted toward direction-dependent strain-limiting materials subjected to nonuniform \cite{GhoshBhattaMallikarjunaiah2026} and compressive loading \cite{BhattaGhoshMallikarjunaiah2026}. Collectively, these investigations demonstrate that crack-tip behavior is highly sensitive to both the chosen implicit constitutive law and the applied load. 

Despite these substantial theoretical and computational developments, the complex interaction between macroscopic loading configurations and localized nonlinear responses in density-dependent elastic models remains largely unexplored. Previous studies have primarily restricted their scope to establishing the constitutive theory, examining preferential stiffness, and resolving fields under idealized canonical fracture modes \cite{Rajagopal2021,YoonMallikarjunaiahBhatta2024,Mallikarjunaiah2026}. However, distinct boundary conditions---such as tension versus compression, uniform versus piecewise-slope loading, or concentrated point loads---induce fundamentally different macroscopic volumetric strain fields. Because the nonlinear constitutive response in the density-dependent model is directly governed by the trace of the infinitesimal strain tensor, the spatial profile of the external displacement-load dictates the local activation of the material nonlinearity. Consequently, observations derived from idealized uniform loadings, or from constitutively distinct strain-limiting models \cite{GhoshBhattaMallikarjunaiah2026,BhattaGhoshMallikarjunaiah2026}, cannot be trivially extrapolated to predict density-dependent crack-tip behavior under arbitrary boundary conditions. 

To bridge this gap, the present study systematically investigates how diverse spatial displacement-loading profiles dictate the localized material response around a stationary crack within a unified density-dependent continuum framework. By fixing the geometry, material description, and numerical discretization strategy, we isolate the interplay between the boundary loading profile and the material nonlinearity. The nonlinear boundary value problem is formulated in an infinitesimal-strain setting, discretized via a conforming $hp$-adaptive finite element method, and solved utilizing an exact Newton scheme. We execute a comprehensive comparative analysis across five distinct loading configurations: uniform tension, uniform compression, piecewise-slope tension, piecewise-slope compression, and a concentrated tensile load. By parametrically varying the constitutive parameter controlling the nonlinear density dependence across these regimes, we demonstrate that the effect of the material nonlinearity is not an intrinsic property of the crack geometry alone, but is profoundly coupled to the macroscopic loading topology. 

The remainder of this paper is organized as follows. Section~\ref{sec:model} introduces the mathematical setting and formally develops the density-dependent constitutive relation. Section~\ref{sec:bvp} establishes the governing boundary value problem and the associated mechanical boundary conditions. 
Section~\ref{sec:solution} details the $hp$-adaptive finite element discretization and the Newton linearization procedure; further, Section~\ref{sec:weak} derives the weak formulation and analyzes its foundational mathematical properties.  Section~\ref{sec:numerical_results} presents the computational experiments across the five loading configurations, scrutinizing the influence of the nonlinear parameter on the near-tip stress and strain fields. Finally, Section~\ref{sec:conclusion} summarizes the principal findings.

\section{Mathematical model}
\label{sec:model}

We consider a stationary cracked elastic body occupying a bounded domain
$\Omega\subset\mathbb{R}^{d}$, with $d\in\{2,3\}$. The numerical
experiments employ the intrinsically two-dimensional form of the model
($d=2$), rather than a plane-stress or plane-strain reduction of a
three-dimensional constitutive law. Let
$\mathbf{u}:\Omega\to\mathbb{R}^{d}$ denote the displacement field,
$\boldsymbol{\varepsilon}(\mathbf{u})$ the infinitesimal strain tensor,
$\boldsymbol{\sigma}$ the Cauchy stress tensor, and $\mathbf{I}$ the identity
tensor. Under the small-displacement gradient assumption,
\begin{equation}
\boldsymbol{\varepsilon}(\mathbf{u})
=\frac{1}{2}\left(\nabla\mathbf{u}+(\nabla\mathbf{u})^{T}\right),
\qquad
\operatorname{tr}\boldsymbol{\varepsilon}(\mathbf{u})
=\nabla\cdot\mathbf{u}.
\label{eq:infinitesimal_strain}
\end{equation}
Here $\nabla\cdot\mathbf{u}$ is the volumetric strain. If $\rho_R$ and
$\rho$ denote the reference and current mass densities, respectively,
conservation of mass gives
\begin{equation}
 \rho_R = \rho\det(\mathbf{I}+\nabla\mathbf{u}),
\end{equation}
and hence, in the small-strain regime
\cite{Rajagopal2021,Mallikarjunaiah2026}, the conservation of mass reduces to
\begin{equation}
\rho/\rho_R=1-\nabla\cdot\mathbf{u}
+O(\|\nabla\mathbf{u}\|^{2}).
\end{equation}
 Thus, density variation is
represented through volumetric strain; density is not treated as an
additional field variable.

The density-dependent constitutive framework considered in
\cite{Rajagopal2021,Mallikarjunaiah2026} includes the relation
\begin{equation}
\left(1+a_3\operatorname{tr}\boldsymbol{\varepsilon}\right)
\boldsymbol{\varepsilon}
=E_1\left(1+a_1\operatorname{tr}\boldsymbol{\varepsilon}\right)
\boldsymbol{\sigma}+E_2\left(1+a_2\operatorname{tr}\boldsymbol{\varepsilon}\right)
\left(\operatorname{tr}\boldsymbol{\sigma}\right)\mathbf{I},
\label{eq:general_density_relation}
\end{equation}
where $E_1$ and $E_2$ are material parameters and $a_1,a_2,a_3$ are
constitutive coefficients. To ensure the mathematical well-posedness of the resulting boundary value problem, we constrain the general implicit constitutive framework by setting $a_1 = a_3 = 0$ and defining $a_2 = \beta$. This specific parameter selection serves a dual purpose. Mathematically, it eliminates the nonlinear coupling in the deviatoric response, isolating the density-dependent nonlinearity strictly within the volumetric strain component. This decoupling is essential for preserving the strong ellipticity and strict monotonicity of the underlying differential operator \cite{Mallikarjunaiah2026,ItouKovtunenkoRajagopal2021,itou2019well,itou2022investigation,itou2023generalization}. Physically, this specialization aligns with established continuum and geomechanical principles, which postulate that purely volumetric deformations---such as local dilation or hydrostatic compression---significantly alter the bulk modulus of a porous solid while its macroscopic shear modulus remains constant. Applying these parameter restrictions yields the signed-volumetric-strain specialization, given by
\begin{equation}\label{eq:magnitude_constitutive}
E_1\boldsymbol{\sigma} = \boldsymbol{\varepsilon} - E_2(1+\beta\operatorname{tr}\boldsymbol{\varepsilon})(\operatorname{tr}\boldsymbol{\sigma})\mathbf{I}.
\end{equation}
In this formulation, the constitutive parameter $\beta$ explicitly modulates the intensity of the nonlinear volumetric response without perturbing the linear shear behavior.
For this particular study, we adopt the above stress-strain relationship. This is a constitutive choice distinct from the signed-trace
specialization above: the nonlinear coefficient depends on the
magnitude of volumetric strain, while the resulting stress retains its
sign dependence. The parameter $\beta$ governs the nonlinear
contribution, and $\beta=0$ defines the corresponding linear reference
model.


Taking the trace of \eqref{eq:magnitude_constitutive} yields
\begin{equation}
\operatorname{tr}\boldsymbol{\sigma}
=\frac{\nabla\cdot\mathbf{u}}
{E_1+dE_2\left(1+\beta|\nabla\cdot\mathbf{u}|\right)}.
\label{eq:stress_trace}
\end{equation}
Substituting \eqref{eq:stress_trace} into \eqref{eq:magnitude_constitutive} gives a new, explicit constitutive relationship:
\begin{equation}
\boldsymbol{\sigma} 
=\frac{1}{E_1}\boldsymbol{\varepsilon}(\mathbf{u})
+H\bigl(\nabla\cdot\mathbf{u}\bigr)\mathbf{I},
\label{eq:stress_strain_H}
\end{equation}
where $H$ is the constitutive function with a scalar value (defined for admissible scalar arguments)
\begin{equation}
H(\xi)
=-\frac{E_2}{E_1}
\frac{\xi\left(1+\beta|\xi|\right)}
{E_1+dE_2\left(1+\beta|\xi|\right)}.
\label{eq:H_definition}
\end{equation}
This representation isolates the nonlinear volumetric contribution and
will be used in the displacement-based weak formulation.

\begin{remark}[Use of the absolute value of volumetric strain and constitutive properties]
\label{rem:abs_volumetric_strain}
Replacing the signed volumetric strain $\operatorname{tr}\boldsymbol{\varepsilon}$ by its magnitude $|\operatorname{tr}\boldsymbol{\varepsilon}| = |\nabla\cdot\mathbf{u}|$ in the nonlinear coefficient of \eqref{eq:magnitude_constitutive} is motivated by both physical and analytical considerations. Throughout this study, we assume $E_1 > 0$ and $E_1 + dE_2 > 0$, and restrict the admissible volumetric strains $\xi = \nabla\cdot\mathbf{u}$ to values for which the constitutive denominator remains strictly positive. For convenience, we set
\begin{equation}
A_0 := E_1 + dE_2, \qquad
B_0 := dE_2\beta, \qquad
D(\xi) := A_0 + B_0|\xi| > 0.
\label{eq:constitutive_denominator}
\end{equation}
On this admissible relation, the scalar constitutive function \eqref{eq:H_definition} can be rewritten in the equivalent form
\begin{equation}
H(\xi) = \frac{1}{d}\left(\frac{\xi}{D(\xi)} - \frac{\xi}{E_1}\right).
\label{eq:H_alternative}
\end{equation}
This magnitude-dependent formulation establishes the following key properties:
\begin{enumerate}
    \item \textbf{Admissibility under both dilation and compression:} In a cracked body subjected to diverse loading scenarios---including tension, compression, and concentrated loads---the volumetric strain $\xi = \nabla\cdot\mathbf{u}$ naturally assumes both positive (local dilation) and negative (local compression) values. Writing the denominator as $D(\xi) = A_0 + B_0|\xi|$ ensures that the nonlinear volumetric modulation depends strictly on the \emph{magnitude} of local volume change. Within the admissible branch $D(\xi) > 0$, this guarantees that the effective bulk modulus remains strictly positive under both local dilation and local compression, avoiding asymmetric loss of ellipticity or unphysical negative bulk moduli.

    \item \textbf{$C^1$ regularity and consistent linearization:} Although $|\xi|$ is not differentiable at $\xi = 0$, the composite function $H(\xi)$ in \eqref{eq:H_alternative} is continuously differentiable wherever $D(\xi) > 0$. Indeed, for $\xi \neq 0$, differentiating $\xi / D(\xi)$ and using $\xi\,\operatorname{sgn}(\xi) = |\xi|$ cancels the non-smooth term in the numerator:
    \[
    \frac{\mathrm{d}}{\mathrm{d}\xi}\left(\frac{\xi}{A_0 + B_0|\xi|}\right)
    = \frac{(A_0 + B_0|\xi|) - \xi B_0 \operatorname{sgn}(\xi)}{(A_0 + B_0|\xi|)^2}
    = \frac{A_0}{D(\xi)^2}.
    \]
    Passing to the limit as $\xi \to 0$ matches the direct difference quotient at the origin, yielding the continuous derivative
    \begin{equation}
    H'(\xi)
    = \frac{1}{d}\left(\frac{A_0}{D(\xi)^2} - \frac{1}{E_1}\right),
    \qquad
    H'(0) = \frac{1}{d}\left(\frac{1}{A_0} - \frac{1}{E_1}\right) = -\frac{E_2}{E_1 A_0}.
    \label{eq:H_derivative_explicit}
    \end{equation}
    Consequently, the absolute-value modification introduces no non-smoothness into the weak formulation, preserving the $C^1$ regularity required for both the Minty--Browder well-posedness analysis (Section~\ref{sec:existence}) and the consistent Newton linearization (Section~\ref{sec:solution}), into which \eqref{eq:H_derivative_explicit} enters directly.

    \item \textbf{Tension--compression symmetry and linear elastic limit:} Because $D(-\xi) = D(\xi)$, the scalar function satisfies $H(-\xi) = -H(\xi)$, which implies the odd symmetry $\boldsymbol{\sigma}(-\boldsymbol{\varepsilon}) = -\boldsymbol{\sigma}(\boldsymbol{\varepsilon})$. Thus, sign-reversed loadings (such as uniform tension versus uniform compression) provide a useful consistency comparison for the present stationary-crack traction-free model, in which crack-face contact is not included. Moreover, setting $\beta = 0$ yields $B_0 = 0$ and $D(\xi) \equiv A_0$, so \eqref{eq:stress_strain_H} reduces directly to the classical linear isotropic relation $\boldsymbol{\sigma} = 2\mu\boldsymbol{\varepsilon} + \lambda(\operatorname{tr}\boldsymbol{\varepsilon})\mathbf{I}$, with Lam\'e parameters $\mu = 1/(2E_1)$ and $\lambda = -E_2/[E_1(E_1 + dE_2)]$.
\end{enumerate}
\end{remark}

\section{Boundary value problem}
\label{sec:bvp}

Next, we formulate the boundary value problem associated with the
constitutive model introduced in Section~\ref{sec:model}. Let
$\Omega\subset\mathbb{R}^{d}$, $d\in\{2,3\}$, denote the cracked
elastic body. In the absence of inertia, the mechanical equilibrium
equation is
\begin{equation}
-\nabla\cdot\boldsymbol{\sigma}
=
\mathbf{f}
\qquad
\text{in }\Omega,
\label{eq:equilibrium}
\end{equation}
where $\mathbf{f}:\Omega\rightarrow\mathbb{R}^{d}$ denotes the
prescribed body force.

The crack domain boundary is decomposed into the Dirichlet
boundary $\Gamma_{D}$, the Neumann boundary $\Gamma_{N}$ and the
 crack boundary without traction $\Gamma_{c}$:
\begin{equation}
\partial\Omega
=
\Gamma_{D}
\cup
\Gamma_{N}
\cup
\Gamma_{c},
\label{eq:boundary_decomposition}
\end{equation}
where these boundary portions are assumed to be mutually disjoint up
to sets of boundary measure zero. The corresponding mechanical
boundary conditions are
\begin{equation}
\mathbf{u}
=
\mathbf{u}_{D}
\quad\text{on }\Gamma_{D},
\qquad
\boldsymbol{\sigma}\mathbf{n}
=
\mathbf{g}
\quad\text{on }\Gamma_{N},
\qquad
\boldsymbol{\sigma}\mathbf{n}
=
\mathbf{0}
\quad\text{on }\Gamma_{c},
\label{eq:boundary_conditions}
\end{equation}
where $\mathbf{n}$ denotes the outward unit normal, $\mathbf{u}_{D}$ is the prescribed displacement, and $\mathbf{g}$
is the prescribed traction. On the two crack faces, the normal is understood separately on each face, and the traction-free condition is imposed accordingly. Here $\Omega$ is the material domain with the crack removed; the two crack faces have independent traces. The stated boundary conditions describe a stationary traction-free crack and do not impose crack-face contact or friction. In particular, the compression calculations are interpreted within this boundary-value model.

Using the constitutive representation
\eqref{eq:stress_strain_H},
\begin{equation}
\boldsymbol{\sigma}
\left(
\boldsymbol{\varepsilon}(\mathbf{u})
\right)
=
\frac{1}{E_{1}}
\boldsymbol{\varepsilon}(\mathbf{u})
+
H(\nabla\cdot\mathbf{u})\mathbf{I},
\end{equation}
together with \[ \operatorname{tr} \boldsymbol{\varepsilon}(\mathbf{u}) =\nabla\cdot\mathbf{u},\]
the equilibrium equation \eqref{eq:equilibrium} can be written entirely in terms of the displacement field as
\begin{equation}
-\nabla\cdot
\left[
\frac{1}{E_{1}}
\boldsymbol{\varepsilon}(\mathbf{u})
+
H(\nabla\cdot\mathbf{u})\mathbf{I}
\right]
=
\mathbf{f}
\qquad
\text{in }\Omega.
\label{eq:displacement_strong_form}
\end{equation}
Thus, the strong form of the nonlinear boundary value problem is to find a displacement field $\mathbf{u}$ that
satisfies \eqref{eq:displacement_strong_form} together with the boundary conditions \eqref{eq:boundary_conditions}. The nonlinearity enters through the scalar constitutive function $H$, defined in \eqref{eq:H_definition}.

The cracked geometry is kept fixed throughout the numerical experiments, whereas the prescribed boundary data are varied to generate the different loading configurations considered in this work. In particular, the uniform, piecewise-slope, and concentrated loading cases are introduced through the appropriate choices of $\mathbf{u}_{D}$ and/or $\mathbf{g}$. The corresponding loading schematics and specific boundary data are presented together with the
numerical experiments in Section~\ref{sec:numerical_results}.


\section{Numerical solution procedure}
\label{sec:solution}

\subsection{Weak formulation}
\label{sec:weak}

Let $\mathbf{u}_{D}$ denote the prescribed displacement on the Dirichlet boundary $\Gamma_D$. 
We introduce the trial and test spaces as
\begin{equation}
\begin{aligned}
V_{\mathbf{u}}&=\bigl\{\mathbf{w}\in[H^1(\Omega)]^d: \mathbf{w}=\mathbf{u}_D\ \text{on }\Gamma_D\bigr\},\\
V_\mathbf{0}&=\bigl\{\mathbf{v}\in[H^1(\Omega)]^d: \mathbf{v}=\mathbf{0}\ \text{on }\Gamma_D\bigr\},
\end{aligned}
\label{eq:trial_space}
\end{equation}
where the boundary conditions are understood in the trace sense. We assume that the prescribed displacement $\mathbf{u}_D$ is the boundary restriction of a function in $[H^1(\Omega)]^d$ and that the Dirichlet conditions eliminate rigid-body motions. The two crack faces possess independent displacement traces. Furthermore, admissible displacement fields are restricted to satisfy the physical admissibility condition:
\begin{equation}
\operatorname*{ess\,inf}_{\mathbf{x}\in\Omega} D\bigl(\nabla\cdot\mathbf{u}(\mathbf{x})\bigr)>0,
\label{eq:weak_admissibility}
\end{equation}
with the denominator $D$ defined in \eqref{eq:constitutive_denominator}.

Multiplying the equilibrium equation \eqref{eq:equilibrium} by a test function $\mathbf{v}\in V_\mathbf{0}$, integrating by parts, and imposing the boundary conditions \eqref{eq:boundary_conditions} yields
\begin{equation}
\int_\Omega \boldsymbol{\sigma}\bigl(\boldsymbol{\varepsilon}(\mathbf{u})\bigr) : \boldsymbol{\varepsilon}(\mathbf{v})\,d\mathbf{x}
=\int_\Omega\mathbf{f}\cdot\mathbf{v}\,d\mathbf{x} + \int_{\Gamma_N}\mathbf{g}\cdot\mathbf{v}\,dS.
\label{eq:weak_stress_form}
\end{equation}
Using the constitutive representation \eqref{eq:stress_strain_H}, the weak formulation is stated as: find an admissible displacement field $\mathbf{u}\in V_{\mathbf{u}}$ such that
\begin{equation}
\mathcal{A}(\mathbf{u};\mathbf{v})=\ell(\mathbf{v}) \qquad\forall\,\mathbf{v}\in V_\mathbf{0},
\label{eq:abstract_weak_form}
\end{equation}
where the nonlinear semi-linear form $\mathcal{A}(\mathbf{u};\mathbf{v})$ is defined by
\begin{equation}
\mathcal{A}(\mathbf{u};\mathbf{v})
=\int_\Omega\left[ \frac{1}{E_1}\boldsymbol{\varepsilon}(\mathbf{u}) : \boldsymbol{\varepsilon}(\mathbf{v}) + H(\nabla\cdot\mathbf{u})(\nabla\cdot\mathbf{v}) \right]d\mathbf{x},
\label{eq:nonlinear_form}
\end{equation}
and the load functional $\ell(\mathbf{v})$ is given by
\begin{equation}
\ell(\mathbf{v}) = \int_\Omega\mathbf{f}\cdot\mathbf{v}\,d\mathbf{x} + \int_{\Gamma_N}\mathbf{g}\cdot\mathbf{v}\,dS.
\label{eq:load_functional}
\end{equation}
The prescribed loads $\mathbf{f}$ and $\mathbf{g}$ are assumed to define a bounded linear functional on $V_\mathbf{0}$.

\subsection{Constitutive assumptions and operator framework}
\label{sec:operator_framework}

To establish the existence and uniqueness of a weak solution to the nonlinear problem \eqref{eq:abstract_weak_form}, we cast the problem into the abstract framework of nonlinear monotone operator theory using the \textit{Minty-Browder Theorem} \cite{showalter2013monotone,brezis2011functional}. Let the trial space be represented as an affine space $V_{\mathbf{u}} = \mathbf{u}_0 + V_\mathbf{0}$, where $\mathbf{u}_0 \in [H^1(\Omega)]^d$ is a fixed lifting of the Dirichlet boundary data $\mathbf{u}_D$ on $\Gamma_D$. We define a nonlinear operator $\mathcal{T}: V_\mathbf{0} \to V_\mathbf{0}^*$ via the duality pairing:
\begin{equation}
\langle \mathcal{T}(\tilde{\mathbf{u}}), \mathbf{v} \rangle = \mathcal{A}(\mathbf{u}_0 + \tilde{\mathbf{u}}; \mathbf{v}) \quad \forall \mathbf{v} \in V_\mathbf{0},
\end{equation}
where $\tilde{\mathbf{u}} = \mathbf{u} - \mathbf{u}_0 \in V_\mathbf{0}$. Similarly, we define the bounded linear functional $F \in V_\mathbf{0}^*$ by $\langle F, \mathbf{v} \rangle = \ell(\mathbf{v})$. Thus, the weak formulation is equivalent to finding $\tilde{\mathbf{u}} \in V_\mathbf{0}$ such that $\mathcal{T}(\tilde{\mathbf{u}}) = F$ in $V_\mathbf{0}^*$.

To guarantee well-posedness and physical admissibility, we impose standard regularity and ellipticity assumptions on the scalar constitutive function $H(\xi)$, where $\xi = \nabla\cdot\mathbf{u}$:
\begin{enumerate}
\item[A1.] \textbf{Regularity:} $H \in C^1(\mathbb{R})$, ensuring that the first derivative $H'(\xi)$ is continuous over the admissible range of volumetric strains.
\item[A2.] \textbf{Strong ellipticity (monotonicity) condition:} There exists a constant $\alpha > 0$ such that for all $\xi \in \mathbb{R}$ satisfying the admissibility constraint~\eqref{eq:weak_admissibility}:
\begin{equation}
\frac{1}{d E_1} + H'(\xi) \ge \alpha > 0,
\label{eq:ellipticity_condition}
\end{equation}
where $d \in \{2, 3\}$ is the spatial dimension. This condition ensures that the tangent bulk modulus remains strictly positive.
\item[A3.] \textbf{Growth restriction:} The derivative $H'(\xi)$ is bounded, meaning there exists a constant $C_H > 0$ such that $\vert H'(\xi)\vert \le C_H$ for all admissible $\xi$. Consequently, $H(\xi)$ exhibits at most linear growth.
\end{enumerate}

\subsection{Existence and uniqueness theorem}
\label{sec:existence}

\begin{Theorem}
\label{thm:existence_uniqueness}
Let $\Omega \subset \mathbb{R}^d$ ($d \in \{2, 3\}$) be a bounded Lipschitz domain with $\vert\Gamma_D\vert > 0$. Assume that the elastic modulus $E_1 > 0$, the external data $\mathbf{f} \in [L^2(\Omega)]^d$ and $\mathbf{g} \in [L^2(\Gamma_N)]^d$ are bounded, and that the constitutive function $H(\xi)$ satisfies the conditions (A1)-(A3). Then, the weak problem admits a unique solution $\mathbf{u} \in V_{\mathbf{u}}$.
\end{Theorem}

\begin{proof}
To apply the Minty-Browder Theorem for pseudomonotone, coercive operators, we must systematically verify three fundamental properties of the operator $\mathcal{T}$ on the real separable Hilbert space $V_\mathbf{0}$: boundedness (and continuity), strict monotonicity, and coercivity.

\textbf{Step 1: Boundedness and continuity.}
By Korn's inequality (valid since $\vert\Gamma_D\vert > 0$), the norm $\|\mathbf{v}\|_{V_\mathbf{0}} = \|\boldsymbol{\varepsilon}(\mathbf{v})\|_{L^2(\Omega)}$ is equivalent to the standard $[H^1(\Omega)]^d$ norm. Using the Cauchy-Schwarz inequality, the operator application can be bounded as:
\begin{equation}
|\langle \mathcal{T}(\tilde{\mathbf{u}}), \mathbf{v} \rangle| \le \frac{1}{E_1} \|\boldsymbol{\varepsilon}(\mathbf{u}_0 + \tilde{\mathbf{u}})\|_{L^2} \|\boldsymbol{\varepsilon}(\mathbf{v})\|_{L^2} + \|H(\nabla\cdot(\mathbf{u}_0 + \tilde{\mathbf{u}}))\|_{L^2} \|\nabla\cdot\mathbf{v}\|_{L^2}.
\end{equation}
Because $H'(\xi)$ is bounded, the function $H(\xi)$ satisfies a linear growth condition 
\[
\vert H(\xi)\vert \le \vert H(0)\vert + C_H\vert\xi\vert.
\]
Consequently, for any $\mathbf{u} \in [H^1(\Omega)]^d$, the mapped term $H(\nabla\cdot\mathbf{u})$ remains bounded in $L^2(\Omega)$. This confirms that $\mathcal{T}$ maps bounded sets in $V_\mathbf{0}$ to bounded sets in $V_\mathbf{0}^*$. Furthermore, because $H \in C^1(\mathbb{R})$ and satisfies the growth condition, the associated \textit{Nemytskii (composition) operator} is continuous from $L^2(\Omega)$ into $L^2(\Omega)$. Therefore, the mapping $\mathcal{T}: V_\mathbf{0} \to V_\mathbf{0}^*$ is continuous, which strictly satisfies the hemicontinuity requirement of the Minty-Browder Theorem.

\textbf{Step 2: Strict monotonicity.}
Let $\mathbf{u}_1, \mathbf{u}_2 \in V_{\mathbf{u}}$ be two distinct, admissible displacement fields. We define their respective strain fields as $\mathbf{E}_1 = \boldsymbol{\varepsilon}(\mathbf{u}_1)$ and $\mathbf{E}_2 = \boldsymbol{\varepsilon}(\mathbf{u}_2)$, and their difference as $\mathbf{E} = \mathbf{E}_1 - \mathbf{E}_2 = \boldsymbol{\varepsilon}(\mathbf{u}_1 - \mathbf{u}_2)$. Examining the operator difference yields:
\begin{equation}
\langle \mathcal{T}(\mathbf{u}_1 - \mathbf{u}_0) - \mathcal{T}(\mathbf{u}_2 - \mathbf{u}_0), (\mathbf{u}_1 - \mathbf{u}_0) - (\mathbf{u}_2 - \mathbf{u}_0) \rangle = \int_\Omega \left[ \frac{1}{E_1} |\mathbf{E}|^2 + \bigl(H(\operatorname{tr}\mathbf{E}_1) - H(\operatorname{tr}\mathbf{E}_2)\bigr)(\operatorname{tr}\mathbf{E}) \right] d\mathbf{x},
\end{equation}
where we have utilized $\nabla \cdot \mathbf{v} = \operatorname{tr}(\boldsymbol{\varepsilon}(\mathbf{v}))$. By the Mean Value Theorem applied to $H$, there exists an intermediate scalar state $\xi_\theta$ such that 
\[
H(\operatorname{tr}\mathbf{E}_1) - H(\operatorname{tr}\mathbf{E}_2) = H'(\xi_\theta)(\operatorname{tr}\mathbf{E}). 
\]
The nonlinear difference term therefore simplifies to $H'(\xi_\theta) (\operatorname{tr}\mathbf{E})^2$.

To bound this from below, we apply the orthogonal decomposition of the strain tensor into its deviatoric and volumetric parts:
\begin{equation}
 \mathbf{E} = \mathbf{E}_D + \frac{1}{d}(\operatorname{tr}\mathbf{E})\mathbf{I}, 
 \end{equation}
 noting that 
 \begin{equation}
 |\mathbf{E}|^2 = |\mathbf{E}_D|^2 + \frac{1}{d}(\operatorname{tr}\mathbf{E})^2. 
 \end{equation}
 Substituting this decomposition into the integrand gives:
\begin{equation}
\frac{1}{E_1} |\mathbf{E}|^2 + H'(\xi_\theta) (\operatorname{tr}\mathbf{E})^2 = \frac{1}{E_1} |\mathbf{E}_D|^2 + \left( \frac{1}{d E_1} + H'(\xi_\theta) \right) (\operatorname{tr}\mathbf{E})^2.
\end{equation}
Invoking the strict ellipticity condition \eqref{eq:ellipticity_condition}, $\frac{1}{d E_1} + H'(\xi) \ge \alpha > 0$, we establish a definitive lower bound:
\begin{align}
\frac{1}{E_1} |\mathbf{E}_D|^2 + \left( \frac{1}{d E_1} + H'(\xi_\theta) \right) (\operatorname{tr}\mathbf{E})^2 &\ge \frac{1}{E_1} |\mathbf{E}_D|^2 + \alpha (\operatorname{tr}\mathbf{E})^2 \nonumber \\
&\ge \min\left(\frac{1}{E_1}, \alpha d\right) \left[ |\mathbf{E}_D|^2 + \frac{1}{d}(\operatorname{tr}\mathbf{E})^2 \right]  \notag \\
&= \tilde{\alpha} |\mathbf{E}|^2,
\end{align}
where $\tilde{\alpha} = \min(1/E_1, \alpha d) > 0$. Integrating this result over the domain $\Omega$ and applying Korn's inequality yields:
\begin{equation}
\langle \mathcal{T}(\tilde{\mathbf{u}}_1) - \mathcal{T}(\tilde{\mathbf{u}}_2), \tilde{\mathbf{u}}_1 - \tilde{\mathbf{u}}_2 \rangle \ge \tilde{\alpha} \|\boldsymbol{\varepsilon}(\tilde{\mathbf{u}}_1 - \tilde{\mathbf{u}}_2)\|_{L^2(\Omega)}^2 \ge \tilde{\alpha} C_K \|\tilde{\mathbf{u}}_1 - \tilde{\mathbf{u}}_2\|_{V_\mathbf{0}}^2 > 0,
\end{equation}
for any $\tilde{\mathbf{u}}_1 \neq \tilde{\mathbf{u}}_2$, where $C_K > 0$ is the optimal Korn's constant. Thus, the operator $\mathcal{T}$ is strictly monotone.

\textbf{Step 3: Coercivity.}
To demonstrate coercivity, we expand the action of the operator on $\tilde{\mathbf{u}}$ around the zero element. Using the established monotonicity bound, we have:
\begin{equation}
\langle \mathcal{T}(\tilde{\mathbf{u}}) - \mathcal{T}(\mathbf{0}), \tilde{\mathbf{u}} - \mathbf{0} \rangle \ge \tilde{\alpha} C_K \|\tilde{\mathbf{u}}\|_{V_\mathbf{0}}^2.
\end{equation}
Rearranging and applying the Cauchy-Schwarz inequality to the linear residual $\mathcal{T}(\mathbf{0})$ gives:
\begin{equation}
\langle \mathcal{T}(\tilde{\mathbf{u}}), \tilde{\mathbf{u}} \rangle \ge \tilde{\alpha} C_K \|\tilde{\mathbf{u}}\|_{V_\mathbf{0}}^2 + \langle \mathcal{T}(\mathbf{0}), \tilde{\mathbf{u}} \rangle \ge \tilde{\alpha} C_K \|\tilde{\mathbf{u}}\|_{V_\mathbf{0}}^2 - \|\mathcal{T}(\mathbf{0})\|_{V_\mathbf{0}^*} \|\tilde{\mathbf{u}}\|_{V_\mathbf{0}}.
\end{equation}
Dividing by the norm of the displacement field yields the coercivity limit:
\begin{equation}
\frac{\langle \mathcal{T}(\tilde{\mathbf{u}}), \tilde{\mathbf{u}} \rangle}{\|\tilde{\mathbf{u}}\|_{V_\mathbf{0}}} \ge \tilde{\alpha} C_K \|\tilde{\mathbf{u}}\|_{V_\mathbf{0}} - \|\mathcal{T}(\mathbf{0})\|_{V_\mathbf{0}^*} \quad \implies \quad \lim_{\|\tilde{\mathbf{u}}\|_{V_\mathbf{0}} \to \infty} \frac{\langle \mathcal{T}(\tilde{\mathbf{u}}), \tilde{\mathbf{u}} \rangle}{\|\tilde{\mathbf{u}}\|_{V_\mathbf{0}}} = \infty.
\end{equation}

\textbf{Step 4: Conclusion.}
Because the operator $\mathcal{T}: V_\mathbf{0} \to V_\mathbf{0}^*$ is bounded, continuous (and thus hemicontinuous), strictly monotone, and coercive on the real separable Hilbert space $V_\mathbf{0}$, the Minty-Browder Theorem guarantees the existence of a unique solution $\tilde{\mathbf{u}} \in V_\mathbf{0}$ to the nonlinear operator equation $\mathcal{T}(\tilde{\mathbf{u}}) = F$. Consequently, a unique, physically admissible weak solution $\mathbf{u} = \mathbf{u}_0 + \tilde{\mathbf{u}} \in V_{\mathbf{u}}$ is guaranteed to exist.
\end{proof}

\subsection{Finite element discretization}

For the two-dimensional computations, let $\mathcal{T}_h$ be a
quadrilateral mesh of $\Omega$, with elementwise polynomial degrees $p_K$.
The conforming $hp$-finite element space is
\begin{equation}
V_h=\left\{\mathbf{w}_h\in[H^1(\Omega)]^2:
\mathbf{w}_h|_K\circ F_K\in[Q_{p_K}(\widehat K)]^2,
\quad K\in\mathcal{T}_h\right\},
\label{eq:hp_space}
\end{equation}
where $F_K:\widehat K\to K$ is the element mapping and
$Q_{p_K}(\widehat K)$ is the tensor-product polynomial space of degree
at most $p_K$ in each reference coordinate. The Hanging-node and
polynomial-degree transition constraints enforce continuity across element
interfaces, without identifying the two crack faces.

Let $V_{h,0}=V_h\cap V_\mathbf{0}$ and
$V_{h,\mathbf{u}}=\mathbf{u}_{D,h}+V_{h,0}$, where
$\mathbf{u}_{D,h}\in V_h$ satisfies the prescribed discrete
displacement conditions on $\Gamma_D$.
The Galerkin approximation seeks an admissible
$\mathbf{u}_h\in V_{h,\mathbf{u}}$ satisfying
\begin{equation}
\mathcal{A}(\mathbf{u}_h;\mathbf{v}_h)=\ell_h(\mathbf{v}_h)
\qquad\forall\,\mathbf{v}_h\in V_{h,0},
\label{eq:discrete_weak_form}
\end{equation}
where $\ell_h$ denotes the discrete load functional.

\subsection{Newton linearization}

Define the nonlinear residual
\begin{equation}
\mathcal{R}_h(\mathbf{u}_h;\mathbf{v}_h)
:=\mathcal{A}(\mathbf{u}_h;\mathbf{v}_h)-\ell_h(\mathbf{v}_h).
\label{eq:residual}
\end{equation}
The directional derivative with respect to $\mathbf{u}_h$ is
\begin{equation}
\begin{aligned}
D\mathcal{R}_h(\mathbf{u}_h;\mathbf{v}_h)
[\delta\mathbf{u}_h]
=\int_\Omega\Bigl[
&\frac{1}{E_1}
\boldsymbol{\varepsilon}(\delta\mathbf{u}_h):
\boldsymbol{\varepsilon}(\mathbf{v}_h)\\
&+H'(\nabla\cdot\mathbf{u}_h)
(\nabla\cdot\delta\mathbf{u}_h)
(\nabla\cdot\mathbf{v}_h)
\Bigr]d\mathbf{x},
\end{aligned}
\label{eq:residual_derivative}
\end{equation}
where $H'$ is given in \eqref{eq:H_derivative_explicit}.
At iteration $k$, the Newton correction
$\delta\mathbf{u}_h^k\in V_{h,0}$ is determined by
\begin{equation}
D\mathcal{R}_h(\mathbf{u}_h^k;\mathbf{v}_h)
[\delta\mathbf{u}_h^k]
=-\mathcal{R}_h(\mathbf{u}_h^k;\mathbf{v}_h)
\qquad\forall\,\mathbf{v}_h\in V_{h,0}.
\label{eq:newton_abstract}
\end{equation}

Let $\{\boldsymbol{\phi}_i\}_{i=1}^{N_{\mathrm{dof}}}$ be a basis of
$V_{h,0}$. The resulting algebraic system is
\begin{equation}
\mathbf{J}^{k}\Delta\mathbf{U}^{k}=-\mathbf{r}^{k},
\qquad
\begin{cases}
r_i^{k}=\mathcal{R}_h(\mathbf{u}_h^k;\boldsymbol{\phi}_i),\\
J_{ij}^{k}=D\mathcal{R}_h(\mathbf{u}_h^k;
\boldsymbol{\phi}_i)[\boldsymbol{\phi}_j],
\end{cases}
\label{eq:assembled_newton_system}
\end{equation}
with
$\delta\mathbf{u}_h^k=\sum_j\Delta U_j^k\boldsymbol{\phi}_j$.
The displacement is updated according to
\begin{equation}
\mathbf{u}_h^{k+1}=\mathbf{u}_h^k
+\alpha^k\delta\mathbf{u}_h^k,
\qquad 0<\alpha^k\le 1,
\label{eq:newton_update}
\end{equation}
where $\alpha^k=1$ corresponds to a full Newton step and a reduced
step may be used to improve convergence. Accepted iterates must satisfy
the constitutive admissibility condition \eqref{eq:weak_admissibility}.
The nonlinear iterations end when
$\|\mathbf{r}^{k}\|_2<\varepsilon_{\mathrm N}$ or the prescribed
maximum number of iterations is reached.

\subsection{\texorpdfstring{$hp$}{hp}-adaptive solution strategy}
\label{subsec:hp_adaptive}

To accurately resolve the highly localized deformations induced by the stationary crack and prescribed boundary loading, the nonlinear boundary value problem is solved within an $hp$-adaptive finite element framework. Unlike uniform mesh refinement, the $hp$-adaptive method dynamically and independently optimizes both the local element size ($h$) and the polynomial approximation order ($p$). This dual strategy is particularly critical for fracture mechanics: geometric $h$-refinement introduces smaller elements to capture the steep gradients and singularities near the crack tip, while spectral $p$-enrichment elevates the polynomial order to achieve exponential convergence in the smooth, far-field regions.  The adaptive algorithm is governed by a rigorous four-stage cycle: 
\[
\textsf{SOLVE} \rightarrow \textsf{ESTIMATE} \rightarrow \textsf{MARK} \rightarrow \textsf{REFINE}. 
\]
To drive the spatial adaptation, local discretization errors are quantified using a residual-based a posteriori Kelly error estimator \cite{KellyGagoZienkiewiczBabuska1983}. For each element $K$ in the active mesh $\mathcal{T}_h$, the error indicator $\eta_K$ evaluates the violation of flux continuity across inter-element boundaries, defined mathematically as:
\begin{equation}
\eta_K = \left( \frac{h_K}{24} \int_{\partial K} \jump{ \nabla \mathbf{u}_h \cdot \mathbf{n} }^2 \, ds \right)^{1/2},
\label{eq:kelly_estimator}
\end{equation}
where $h_K$ denotes the characteristic cell diameter, $\mathbf{n}$ is the unit outward normal vector to the element face $\partial K$, and $\jump{\cdot}$ represents the jump operator across the element interface. In the context of the present constitutive model, regions of severe localization---such as the crack-tip vicinity---generate large flux discontinuities. The indicator $\eta_K$ naturally identifies these regions, ensuring that computational resources are concentrated exclusively where the current finite element approximation $\mathbf{u}_h$ lacks sufficient resolution.

Once the elements are marked for refinement based on the magnitude of $\eta_K$, the algorithm must arbitrate between $h$-subdivision and $p$-enrichment. This decision is dictated by a local regularity estimator that analyzes the decay rate of the Legendre-Fourier expansion coefficients of the local displacement solution \cite{Mallikarjunaiah2026}. If the expansion coefficients exhibit a rapid algebraic decay, the solution is deemed locally smooth, and the element is marked for $p$-enrichment to leverage spectral convergence (up to an optional prescribed upper bound $p_{\max}$). Conversely, if the coefficient decay is slow---indicative of the reduced analytical regularity associated with the crack-tip singularity---the element is flagged for $h$-subdivision to geometrically isolate the non-smooth behavior.

Following the marking phase, the finite element space and its associated algebraic constraints are reconstructed. Hanging-node constraints and polynomial-degree transition matrices are strictly enforced to preserve global conformity between neighboring elements of disparate sizes and polynomial orders, while the two opposing crack faces retain kinematically independent displacement traces. The displacement field from the previous mesh is projected onto the newly adapted space to initialize the Newton-Raphson iterations. During the nonlinear solution phase, accepted Newton iterates are required to satisfy the physical admissibility condition defined in \eqref{eq:weak_admissibility}. Upon satisfying the convergence criterion on the current mesh, the displacement field is utilized to evaluate the stress and strain components via \eqref{eq:stress_strain_H}, and the adaptation cycle continues until the prescribed computational tolerances are achieved. The entire $hp$-adaptive procedure is implemented leveraging the conforming multi-constraint facilities of the \texttt{deal.II} finite element library \cite{ArndtBangerthEtAl2021}. 

The complete computational procedure, including initialization,
Newton iteration, postprocessing, and $hp$-adaptation, is summarized
in Algorithm~\ref{alg:ddm_hp}. Here, $m$ denotes the adaptive-cycle
index and $k$ the Newton iteration index on a fixed mesh.

\begin{algorithm}[H]
\small
\caption{$hp$-adaptive Newton scheme}
\label{alg:ddm_hp}
\begin{algorithmic}

\Statex \textbf{Input:}
Initial mesh, material parameters, and boundary data.
\Statex Newton tolerance $\varepsilon_{\mathrm{N}}$,
iteration limit $k_{\max}$, degree limit $p_{\max}$,
and mesh cycles $N$.

\Statex \textbf{Step 1: Initialization}
\State Read and globally refine the mesh.
\State Construct the displacement space and constraints.
\State Compute an auxiliary linear solution.
\State On each adapted mesh, transfer the current solution and reapply the displacement constraints before the Newton solve.

\For{$m=0,\ldots,N-1$}

    \State \textbf{Step 2: Newton iteration}
    \State Initialize $\mathbf{U}^{0}$ from the current solution.
    \For{$k=0,\ldots,k_{\max}-1$}
        \State Assemble residual $\mathbf{r}^{k}$
        and tangent $\mathbf{J}^{k}$.
        \State \textbf{break} if
        $\|\mathbf{r}^{k}\|_2<\varepsilon_{\mathrm{N}}$.
        \State Solve
        $\mathbf{J}^{k}\Delta\mathbf{U}^{k}=-\mathbf{r}^{k}$.
        \State Select the step length $\alpha^{k}$.
        \State Update
        $\mathbf{U}^{k+1}=\mathbf{U}^{k}
        +\alpha^{k}\Delta\mathbf{U}^{k}$.
        \State Enforce displacement constraints.
    \EndFor
    \State Recompute the residual of the final accepted iterate.
    \If{the convergence criterion is not satisfied}
        \State Report nonlinear-solver failure and stop the current run.
    \EndIf

    \State \textbf{Step 3: Postprocessing}
    \State Compute stress and strain.
    \State Export fields and sampling-line data.

    \If{$m<N-1$}
        \State \textbf{Step 4: $hp$-adaptation}
        \State Compute Kelly indicators and Fourier smoothness.
        \State Mark $h$- and $p$-changes.
        \State Adapt the mesh and polynomial degrees.
        \State Rebuild the space and transfer the solution.
        \State Enforce displacement constraints.
    \EndIf

\EndFor
\end{algorithmic}
\end{algorithm}

\section{Numerical results and discussion} 
\label{sec:numerical_results}

This section systematically investigates the mechanical response of the cracked domain under the five distinct loading configurations introduced previously. The primary objective is to elucidate the role of the nonlinearity parameter $\beta$ in governing the near-tip stress and strain fields. By comparing the material behavior under uniform tensile, compressive, spatially varying, and concentrated loading conditions, we demonstrate how local density variations alter the energetic and kinematic landscape of the fracture process zone.

Let $\mathbf{x} = (x, y)$ denote the spatial coordinates in the reference configuration, with the corresponding displacement field approximated by the finite element solution $\mathbf{u} = (u_1, u_2)$. The kinematic quantities of primary interest for assessing local failure are the vertical normal strain, $\varepsilon_{22} = \partial u_2/\partial y$, and the associated vertical normal stress, $\sigma_{22}$. Under the assumption of infinitesimal deformation in two dimensions, the volumetric strain driving the constitutive nonlinearity is given by $\mathrm{tr}(\boldsymbol{\varepsilon}) = \varepsilon_{11} + \varepsilon_{22}$. 

To accurately capture the severe localized gradients, the normal stress $\sigma_{22}$ and normal strain $\varepsilon_{22}$ are sampled along a predefined line directly approaching the crack tip. For the four distributed loading cases, these sampling profiles are supplemented with full-field contour panels that illustrate the global deformation topology and the extent of crack-tip localization. The response to the concentrated load is reported exclusively through high-resolution sampling-line plots and tabulated extremal values.

\subsection{Computational setup}

To isolate the effects of the volumetric strain-dependent nonlinearity, the applied boundary conditions are held strictly constant within each configuration while the material parameter $\beta$ is systematically varied. This parametric variation reveals the transition from classical linear elastic behavior to the highly compliant, strain-limiting response characteristic of the proposed constitutive model.

For all distributed loading scenarios, the prescribed displacement values are given along the top boundary. Conversely, for the concentrated tensile loading configuration, the applied central force is $P_{\mathrm{center}}=0.005$. It should be noted that these quantities specify fundamentally different classes of boundary loading; they are calibrated to probe the nonlinear regime of the material without inducing numerical instability, rather than to enforce statically equivalent resultant forces across configurations. The full-field spatial distributions, highlighting the profound localized variations in the mechanical response, are visualized at a high nonlinearity threshold of $\beta=5$.

\begin{figure}[H]
\centering
\begin{subfigure}[t]{0.3\textwidth}
    \centering
    \includegraphics[width=\linewidth,keepaspectratio]
    {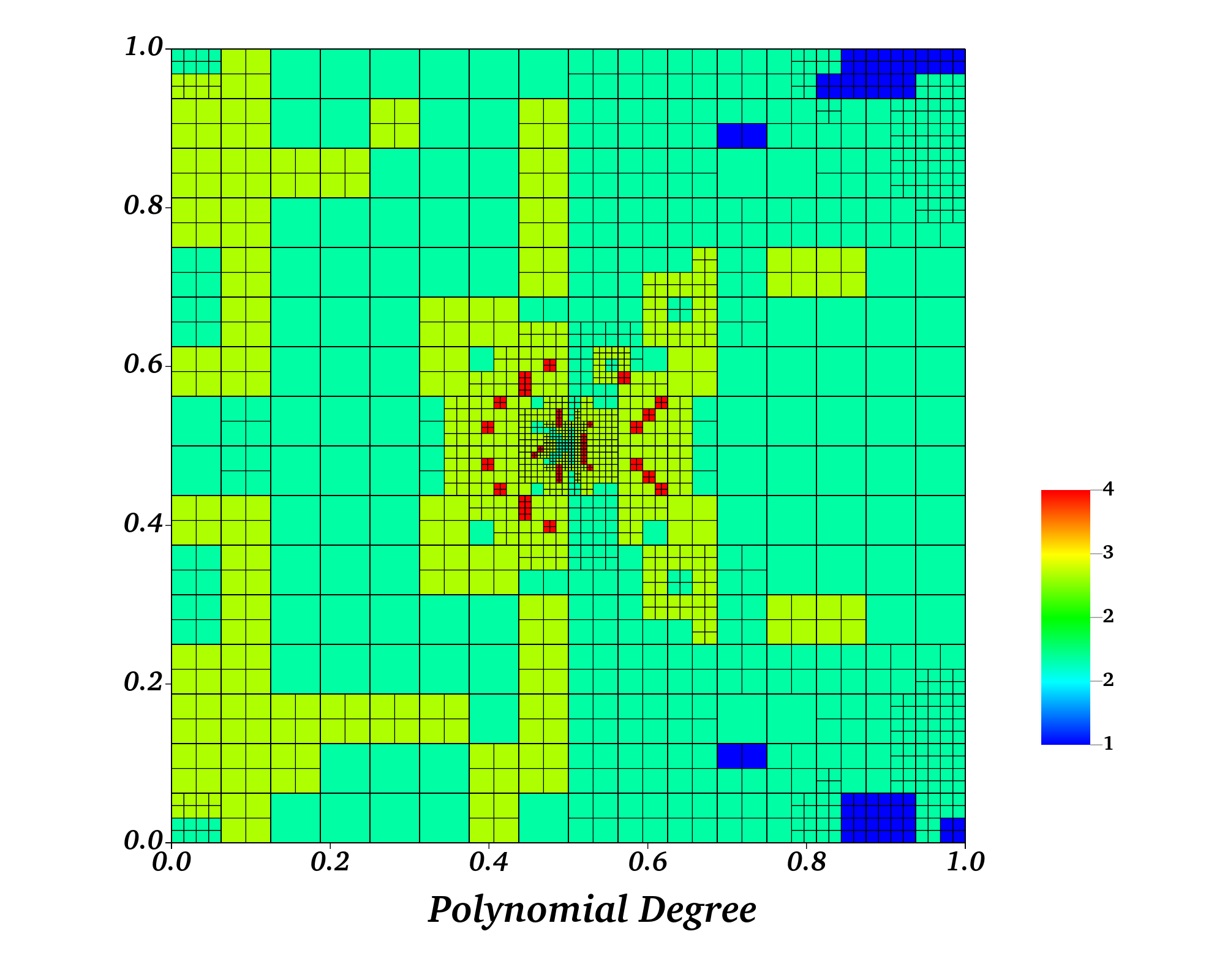}
    \caption{Uniform tension.}
    \label{fig:hp_uniform_tension}
\end{subfigure}
\begin{subfigure}[t]{0.3\textwidth}
    \centering
    \includegraphics[width=\linewidth,keepaspectratio]
    {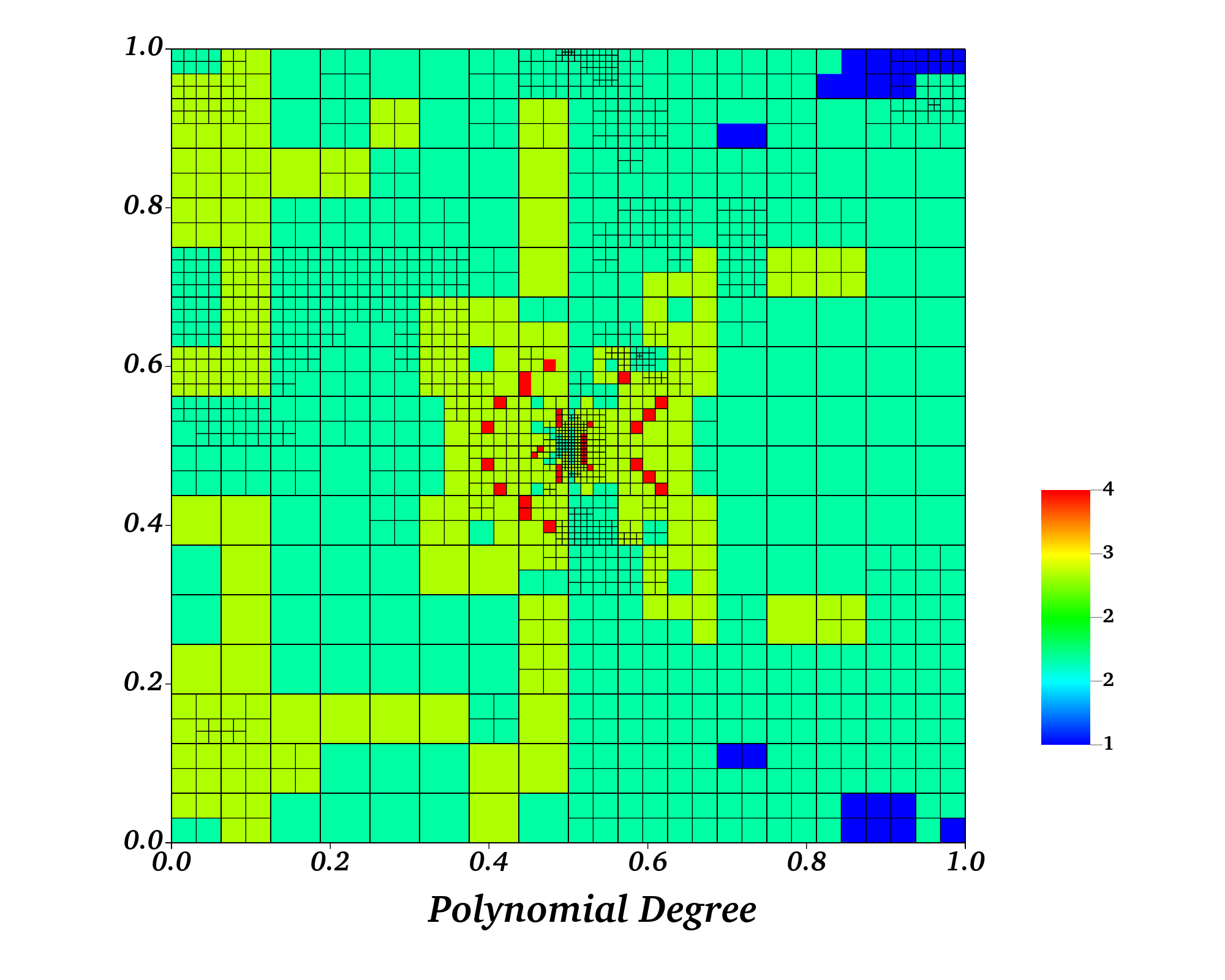}
    \caption{Piecewise-slope tension.}
    \label{fig:hp_piecewise_tension}
\end{subfigure}
\medskip
\begin{subfigure}[t]{0.3\textwidth}
    \centering
    \includegraphics[width=\linewidth,keepaspectratio]
    {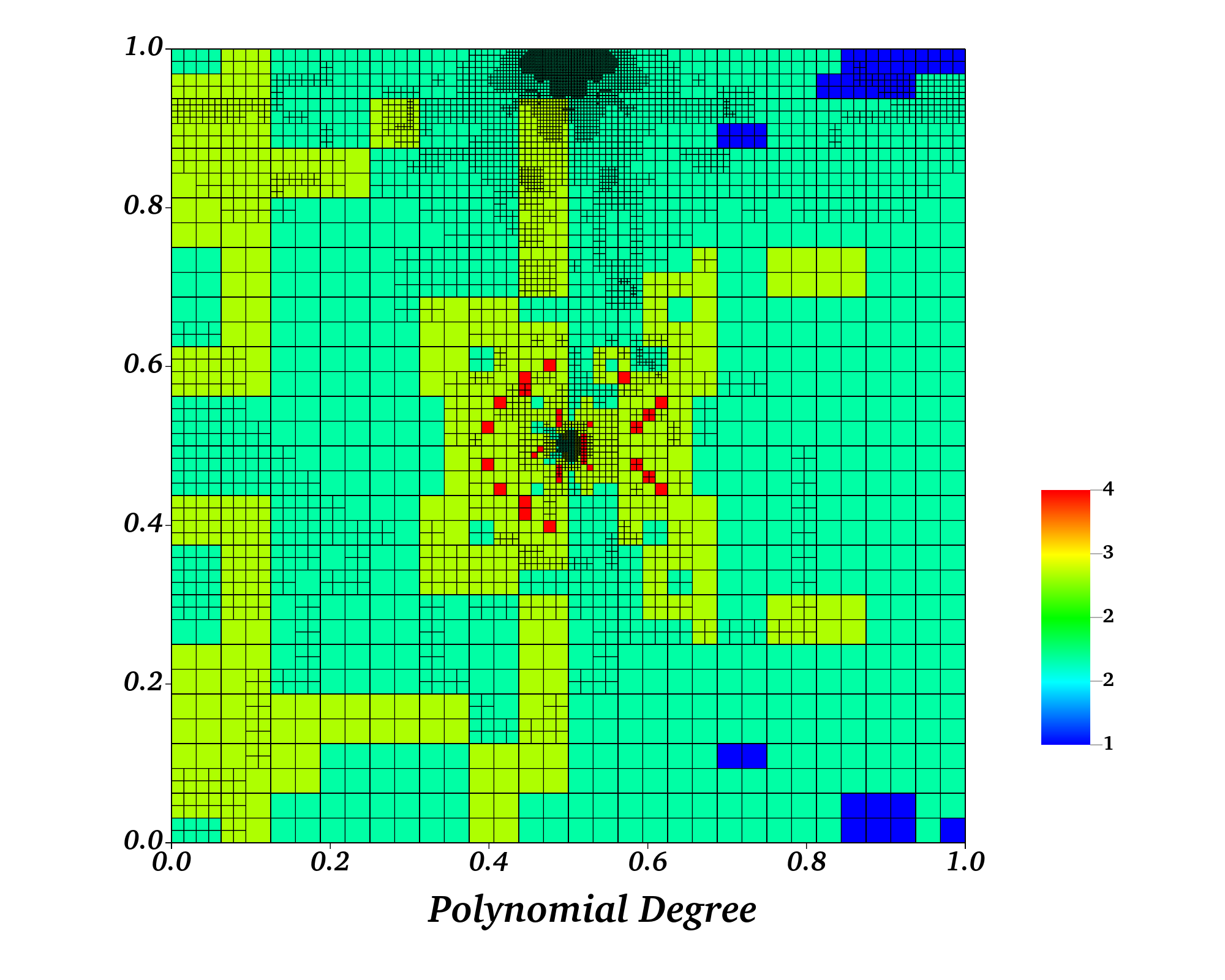}
    \caption{Central point load. }
    \label{fig:hp_point_load}
\end{subfigure}
\caption{Adapted computational meshes and elementwise polynomial degree distributions for three representative loading configurations at $\beta=5$ (adaptive cycle 5). The dynamic refinement isolates the crack-tip singularity while enriching the smooth bulk domain.}
\label{fig:hp_loading_comparison}
\end{figure}

Figure~\ref{fig:hp_loading_comparison} illustrates the terminal state of the $hp$-adaptive discretization (at adaptive cycle 5) for three representative loading scenarios evaluated at $\beta=5$. The geometric boundaries indicate the intensity of the local $h$-refinement, while the color map denotes the polynomial degree ($p$) assigned to each element. This dual-refinement strategy autonomously isolates the singular continuum behavior at the crack tip through dense spatial subdivision while exploiting high-order polynomial bases to achieve spectral convergence in the smooth far-field regions.

\subsubsection{Postprocessing of crack-tip fields.}
Upon convergence of the Newton--Raphson iteration on a given adapted mesh, the discrete displacement field $\mathbf{u}_h \in V_{h,\mathbf{u}}$ is postprocessed to recover the local kinematic, kinetic, and energetic quantities of interest. First, the discrete infinitesimal strain tensor is computed directly from the symmetric gradient of the finite element displacement field in accordance with \eqref{eq:infinitesimal_strain}:
\begin{equation}
\boldsymbol{\varepsilon}_h
:= \boldsymbol{\varepsilon}(\mathbf{u}_h)
= \frac{1}{2}\left(\nabla\mathbf{u}_h + (\nabla\mathbf{u}_h)^{T}\right),
\qquad
\operatorname{tr}\boldsymbol{\varepsilon}_h = \nabla\cdot\mathbf{u}_h.
\label{eq:post_strain}
\end{equation}
Next, the corresponding Cauchy stress tensor $\boldsymbol{\sigma}_h$ is evaluated pointwise by substituting $\boldsymbol{\varepsilon}_h$ and its trace into the explicit density-dependent constitutive relation \eqref{eq:stress_strain_H}:
\begin{equation}
\boldsymbol{\sigma}_h
:= \boldsymbol{\sigma}(\boldsymbol{\varepsilon}_h)
= \frac{1}{E_1}\boldsymbol{\varepsilon}_h
+ H(\nabla\cdot\mathbf{u}_h)\mathbf{I},
\label{eq:post_stress}
\end{equation}
where $H(\cdot)$ is defined in \eqref{eq:H_definition}. Finally, to quantify the local energetic concentration within the fracture process zone, we evaluate the local strain energy density (SED) measure from the tensor contraction of the computed stress and strain fields:
\begin{equation}
W_h := \boldsymbol{\sigma}_h : \boldsymbol{\varepsilon}_h
= \sum_{i,j=1}^{d} (\sigma_h)_{ij} (\varepsilon_h)_{ij}.
\label{eq:post_sed}
\end{equation}
These derived fields are exported across the full computational domain for two-dimensional contour visualization and sampled along the line in front of the crack tip to extract the normal stress $\sigma_{22}$, normal strain $\varepsilon_{22}$, and strain energy density profiles.

\subsection{$hp$-convergence study}\label{subsec:ddm-mesh-convergence}

To verify the accuracy of the finite element implementation of the
nonlinear density-dependent elasticity model, we carry out a convergence study
using the method of manufactured solutions. The auxiliary problem is defined on
the unit square $\Omega=(0,1)^2$ and satisfies the equilibrium equation. 
The constitutive response implemented in the numerical scheme is written as
\begin{equation}\label{eq:ddm-constitutive-verification}
  \boldsymbol\sigma(\mathbf u)
  =\frac{1}{E_1}\,\boldsymbol\varepsilon(\mathbf u)
   +H(\nabla\!\cdot\mathbf u)\,\mathbf I,
  \qquad
  H(\xi)=- \frac{E_2}{E_1} \;  \frac{\xi(1+\beta|\xi|)}
  {E_1+2E_2+2E_2\beta|\xi|},
\end{equation}
where
$E_1=(1+\nu)/E$, and $E_2=-\nu/E$.
We prescribe the smooth displacement field
\begin{equation}\label{eq:ddm-mms-exact}
  \mathbf u_{\mathrm{ex}}(x,y)
  =\delta\begin{pmatrix}e^x-1\\ e^y-1\end{pmatrix},
  \qquad (x,y)\in\Omega,
\end{equation}
and construct the body force consistently from
$\mathbf f=-\nabla\!\cdot\boldsymbol\sigma(\mathbf u_{\mathrm{ex}})$.
The exact displacement is imposed as a Dirichlet condition on the entire
boundary $\partial\Omega$. The material and loading parameters are fixed at
$\delta=0.1$, $\beta=0.5$, $E=1$, and $\nu=0.35$ in both convergence
experiments. The nonlinear discrete equilibrium equations are solved by
Newton's method with a line search, and the resulting displacement errors
are evaluated using an independent higher-order quadrature rule.

Let $\mathbf u_{h,p}$ denote the computed displacement using continuous
quadrilateral $Q_p$ finite elements. The error is measured in the
$L^2(\Omega)$ norm,
\begin{equation}\label{eq:ddm-l2-error}
  e_{L^2}(h,p)
  =\|\mathbf u_{\mathrm{ex}}-\mathbf u_{h,p}\|_{L^2(\Omega)}
  =\left(\int_{\Omega}
  |\mathbf u_{\mathrm{ex}}-\mathbf u_{h,p}|^2\,\mathrm d\Omega\right)^{1/2}.
\end{equation}

\subsubsection{$h$-refinement study}
For the first experiment, the polynomial degree is fixed at $p=1$, and
the unit square is partitioned into uniform $N\times N$ meshes, with
$N=4,8,16,32,64,128,256,512$ and element width $h=1/N$.
The observed convergence order between two consecutive meshes is calculated as
\begin{equation}\label{eq:ddm-h-rate}
  r_h=\frac{\log\!\left(e_{L^2}(h_{j-1},1)/e_{L^2}(h_j,1)\right)}
  {\log(h_{j-1}/h_j)}.
\end{equation}
Table~\ref{tab:ddm-h-convergence} shows that the $L^2$ displacement error
falls by approximately a factor of four whenever the mesh width is halved.
The observed orders remain close to $2$, in agreement with the expected
second-order $L^2$ convergence of bilinear elements for a sufficiently
smooth displacement field. In particular, the error decreases from
$1.436074\times10^{-3}$ on the coarsest mesh to
$8.795478\times10^{-8}$ on the finest mesh, confirming the anticipated
spatial convergence behavior for the manufactured problem.

\begin{table}[H]
\centering
\caption{$h$-convergence of the displacement $L^2$ error for the
manufactured solution using bilinear $Q_1$ elements.}
\label{tab:ddm-h-convergence}
\small
\setlength{\tabcolsep}{6pt}
\begin{tabular}{rrrrrr}
\toprule
$N$ & Cells & DoFs & $h$ & $L^2$ error & Rate \\
\midrule
  4 & 16 & 50 & 0.250000 & $1.436074\times10^{-3}$ & -- \\
  8 & 64 & 162 & 0.125000 & $3.599515\times10^{-4}$ & 1.996 \\
  16 & 256 & 578 & 0.062500 & $9.004628\times10^{-5}$ & 1.999 \\
  32 & 1,024 & 2,178 & 0.031250 & $2.251522\times10^{-5}$ & 2.000 \\
  64 & 4,096 & 8,450 & 0.015625 & $5.629031\times10^{-6}$ & 2.000 \\
  128 & 16,384 & 33,282 & 0.007812 & $1.407272\times10^{-6}$ & 2.000 \\
  256 & 65,536 & 132,098 & 0.003906 & $3.518189\times10^{-7}$ & 2.000 \\
  512 & 262,144 & 526,338 & 0.001953 & $8.795478\times10^{-8}$ & 2.000 \\
\bottomrule
\end{tabular}
\end{table}

\subsubsection{$p$-refinement study}
\label{sec:p_refinement}

In the second convergence experiment, the spatial mesh is kept fixed at
$4\times4$ uniform quadrilateral elements ($h=1/4$), while the polynomial approximation
order is increased successively from $Q_1$ to $Q_6$. The manufactured
solution, material parameters, and boundary conditions remain unchanged
from the $h$-refinement experiment. The resulting comparison therefore
isolates the influence of polynomial enrichment on the discretization error.

Table~\ref{tab:p_convergence} reports the number of degrees of freedom and
the corresponding displacement errors in the $L^2(\Omega)$ norm
introduced above. Over the range of polynomial
degrees considered, the computed errors are closely described by an
exponential relation of the form
\[
    e_{L^2}(1/4,p) \approx A e^{-bp}, \qquad
    A=0.152793,\quad b=4.422509,\quad p=1,\ldots,6.
\]

\begin{table}[H]
\centering
\caption{$p$-convergence of the displacement $L^2$ error for the manufactured solution on a fixed $4\times4$ mesh (16 elements).}
\label{tab:p_convergence}
\begin{tabular}{@{}crr@{}}
\toprule
Degree & DoFs & $L^2$ error \\
\midrule
$Q_1$ & 50    & $1.436074\times10^{-3}$ \\
$Q_2$ & 162   & $2.261699\times10^{-5}$ \\
$Q_3$ & 338   & $3.267063\times10^{-7}$ \\
$Q_4$ & 578   & $3.931880\times10^{-9}$ \\
$Q_5$ & 882   & $4.006448\times10^{-11}$ \\
$Q_6$ & 1,250 & $3.526909\times10^{-13}$ \\
\bottomrule
\end{tabular}
\end{table}

\begin{figure}[H]
\centering
\includegraphics[width=0.72\textwidth]{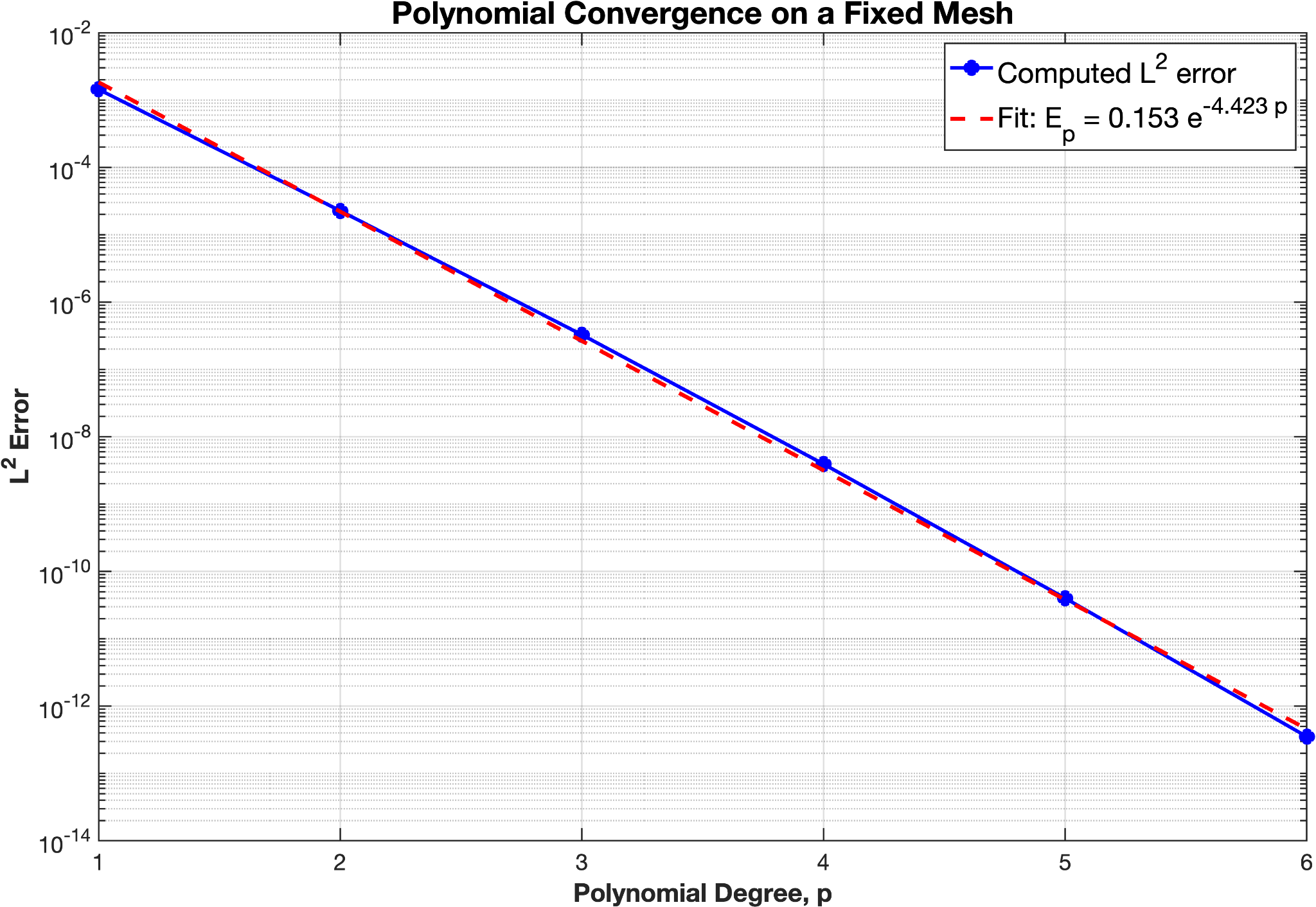}
\caption{Semilogarithmic plot of the $L^2$ displacement error and its exponential approximation as the polynomial degree increases on a fixed $4\times4$ mesh.}
\label{fig:p_convergence}
\end{figure}

Figure~\ref{fig:p_convergence} presents the errors together with the
exponential approximation on a semilogarithmic scale. The nearly linear
decrease in the plotted errors is consistent with exponential-type
convergence under $p$-refinement for the smooth manufactured solution. Together with the second-order
$h$-convergence reported above, these results support the consistency
and accuracy of the finite element approximation for the manufactured problem.

\subsection{Uniform tensile loading}
\begin{figure}[H]
    \centering
    \begin{subfigure}[t]{0.32\linewidth}
        \centering
        \includegraphics[width=\linewidth]
        {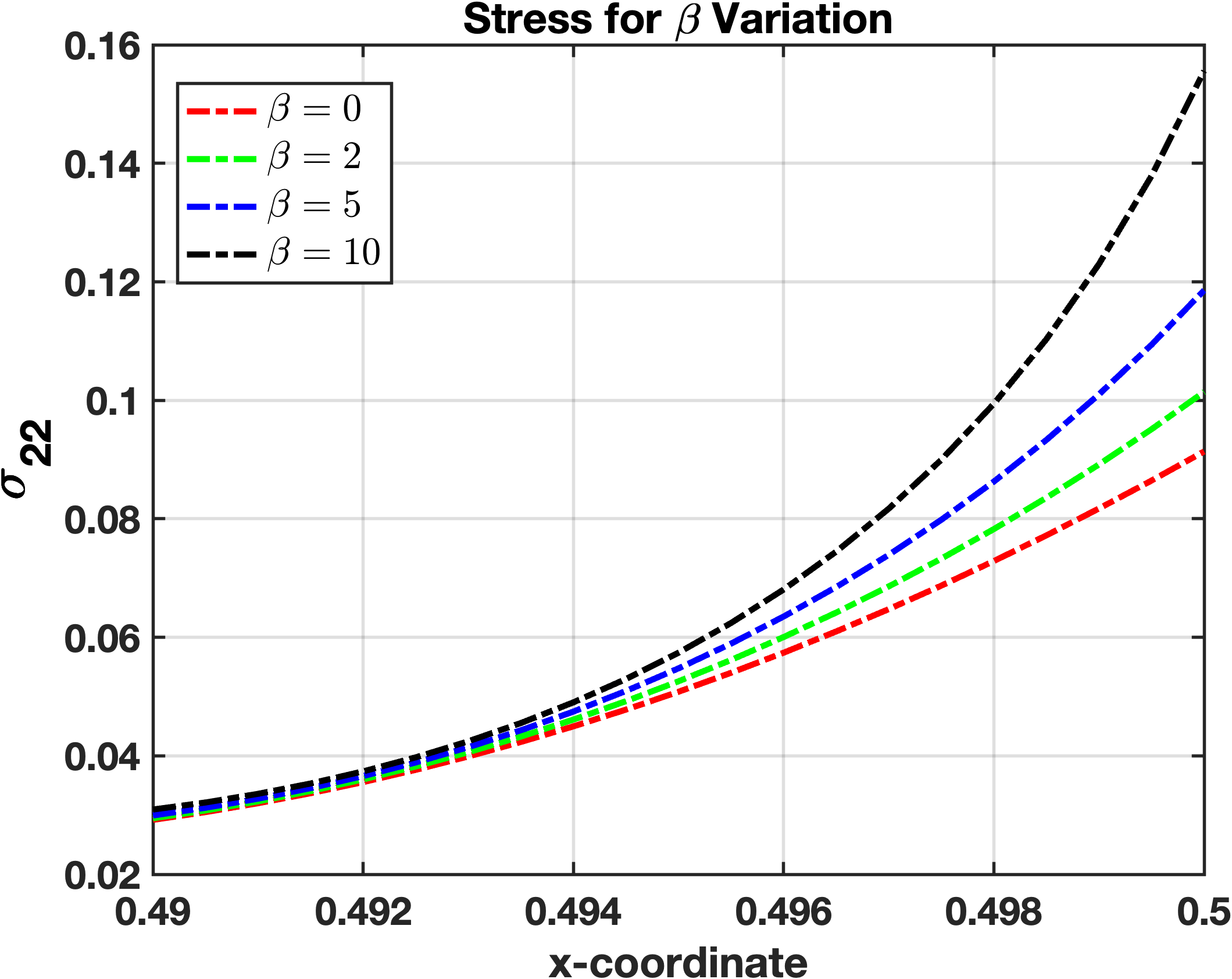}
        \caption{Normal stress $\sigma_{22}$.}
        \label{fig:UT_matlab_stress}
    \end{subfigure}%
    \hfill
    \begin{subfigure}[t]{0.32\linewidth}
        \centering
        \includegraphics[width=\linewidth]
        {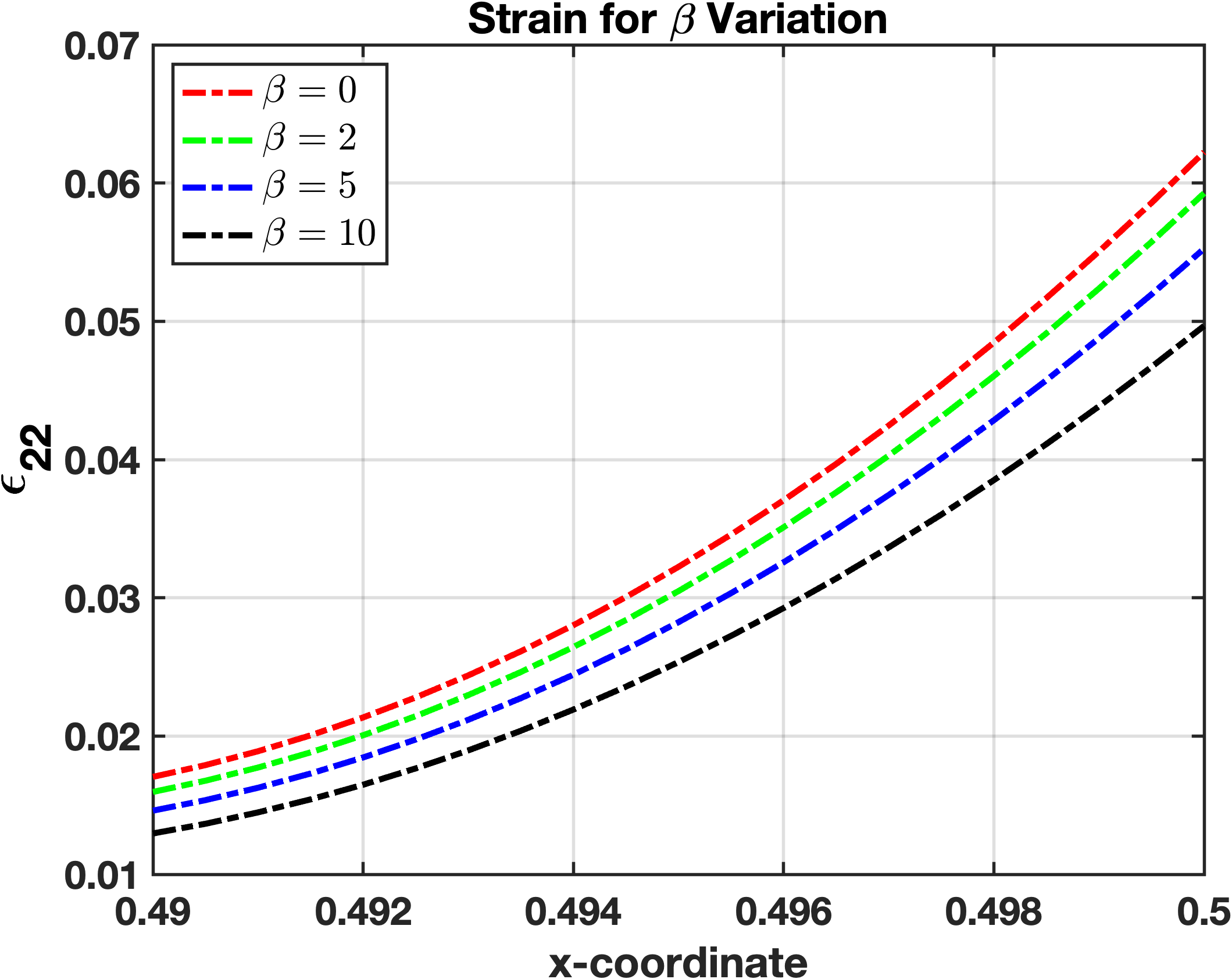}
        \caption{Normal strain $\varepsilon_{22}$.}
        \label{fig:UT_matlab_strain}
    \end{subfigure}
    \hfill
    \begin{subfigure}[t]{0.32\linewidth}
        \centering
        \includegraphics[width=\linewidth]
        {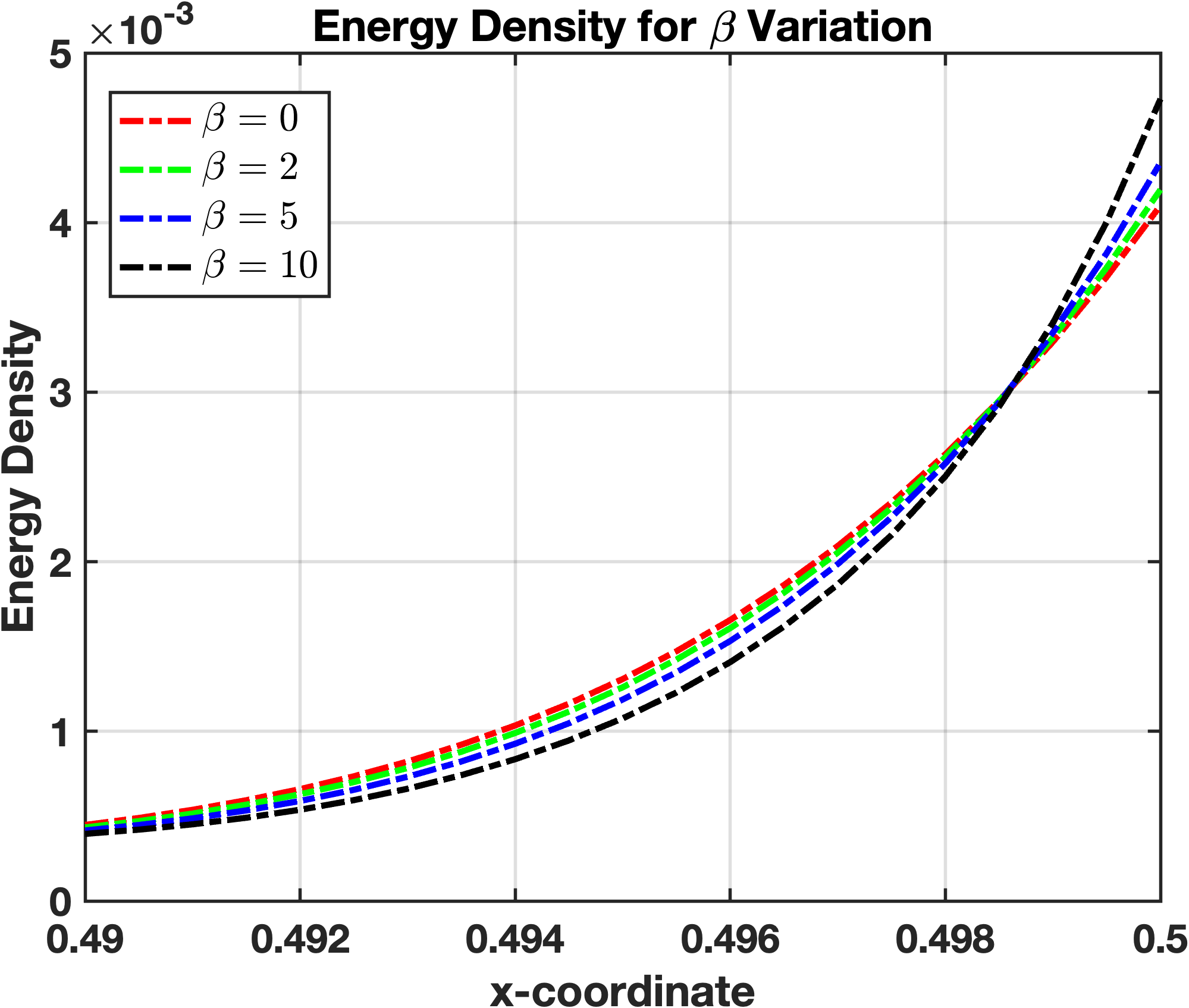}
        \caption{Strain Energy Density.}
        \label{fig:UT_matlab_strain}
    \end{subfigure}

    \caption{Evolution of the normal stress $\sigma_{22}$ and normal strain $\varepsilon_{22}$ along the uncracked ligament approaching the crack tip under uniform tensile loading.}
    \label{fig:UT_matlab}
\end{figure}

The influence of the nonlinearity parameter $\beta$ is most profoundly exhibited in the immediate vicinity of the crack tip. As depicted in Fig.~\ref{fig:UT_matlab}, elevating $\beta$ significantly amplifies the localized tensile stress concentration, while concurrently suppressing the normal strain $\varepsilon_{22}$. This inverse kinematic relationship is a fundamental hallmark of the strain-limiting constitutive framework. As the material approaches its prescribed deformational threshold, it inherently resists further kinematic stretching; consequently, the continuum must balance the uniform external loading through a severe, localized escalation in internal stress. The parameter $\beta$ thus governs a critical redistribution of the local mechanical state, transforming the classical singular strain field into a highly stress-dominated fracture process zone. Furthermore, this restricted local deformability under elevated stress mandates a higher concentration of strain energy density at the crack tip for larger $\beta$ configurations.
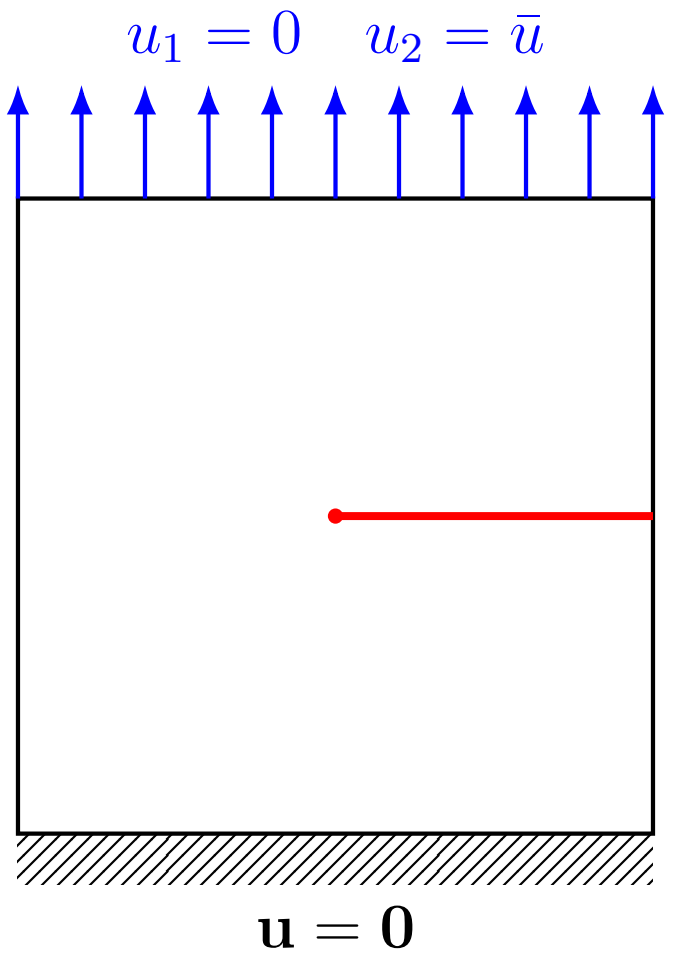
\begin{figure}[H]
\centering
\begin{tikzpicture}[x=4.2cm,y=4.2cm]
  \draw[thick] (0,0) rectangle (1,1);
  \draw[red,line width=1.5pt] (0.5,0.5)--(1,0.5);
  \fill[red] (0.5,0.5) circle[radius=0.012];
  \fill[pattern=north east lines] (0,-0.08) rectangle (1,0);
  
  \foreach \x in {0,0.1,...,1.1} \draw[-{Latex[length=2mm]},thick,blue] (\x,1)--++(0,0.18);
  
  \node[blue,above] at (0.5,1.18) {$u_1 = 0 \quad u_2=\bar u$};
  \node[below] at (0.5,-0.08) {$\mathbf u=\mathbf 0$};
\end{tikzpicture}
\caption{Uniform tensile boundary condition.}
\label{fig:ddm-schematic-i}
\end{figure}

\begin{figure}[H]
    \centering
    \begin{subfigure}[t]{0.48\textwidth}
        \centering
        \includegraphics[width=\textwidth]
        {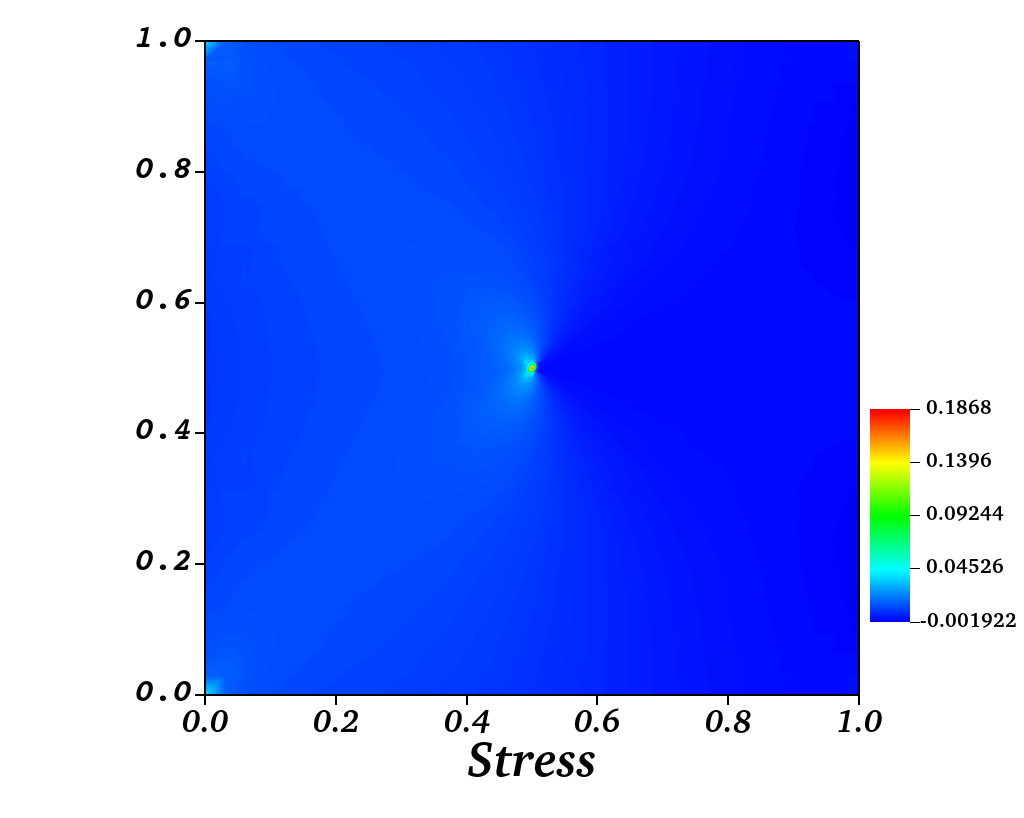}
        \caption{Stress field $\sigma_{22}$.}
        \label{fig:UT_visit_stress}
    \end{subfigure}
    \hfill
    \begin{subfigure}[t]{0.48\textwidth}
        \centering
        \includegraphics[width=\textwidth]
        {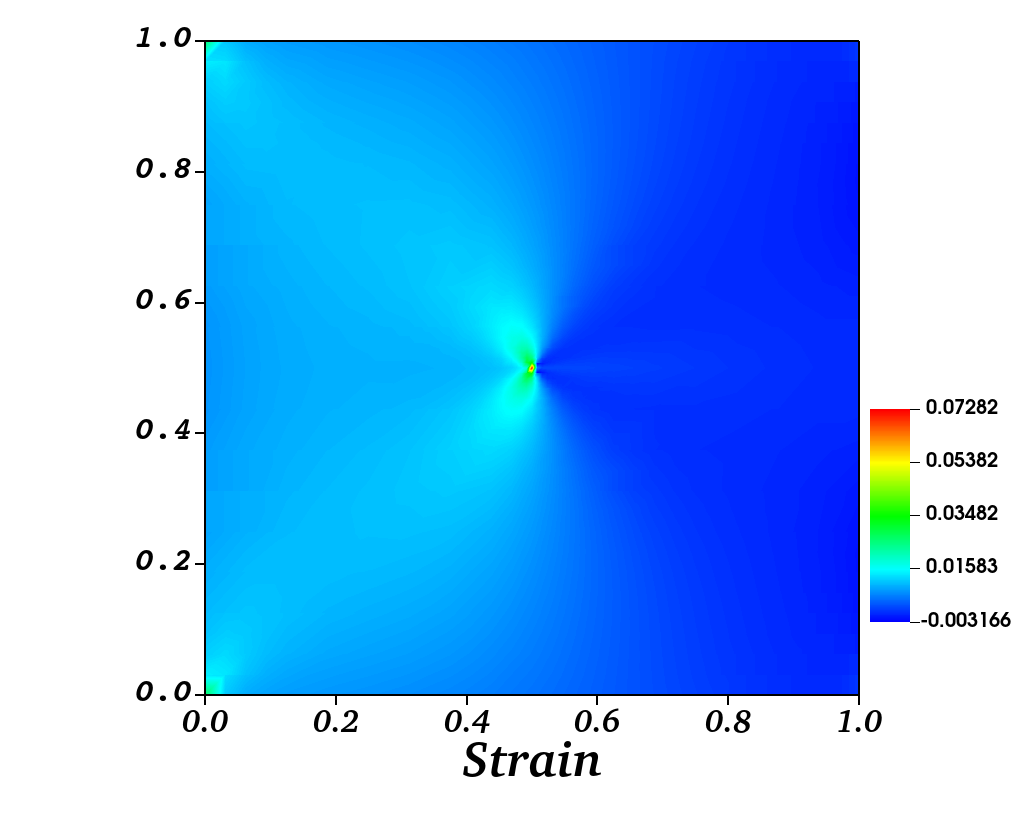}
        \caption{Strain field $\varepsilon_{22}$.}
        \label{fig:UT_visit_strain}
    \end{subfigure}

    \vspace{0.3cm}

    \begin{subfigure}[t]{0.48\textwidth}
        \centering
        \includegraphics[width=\textwidth]
        {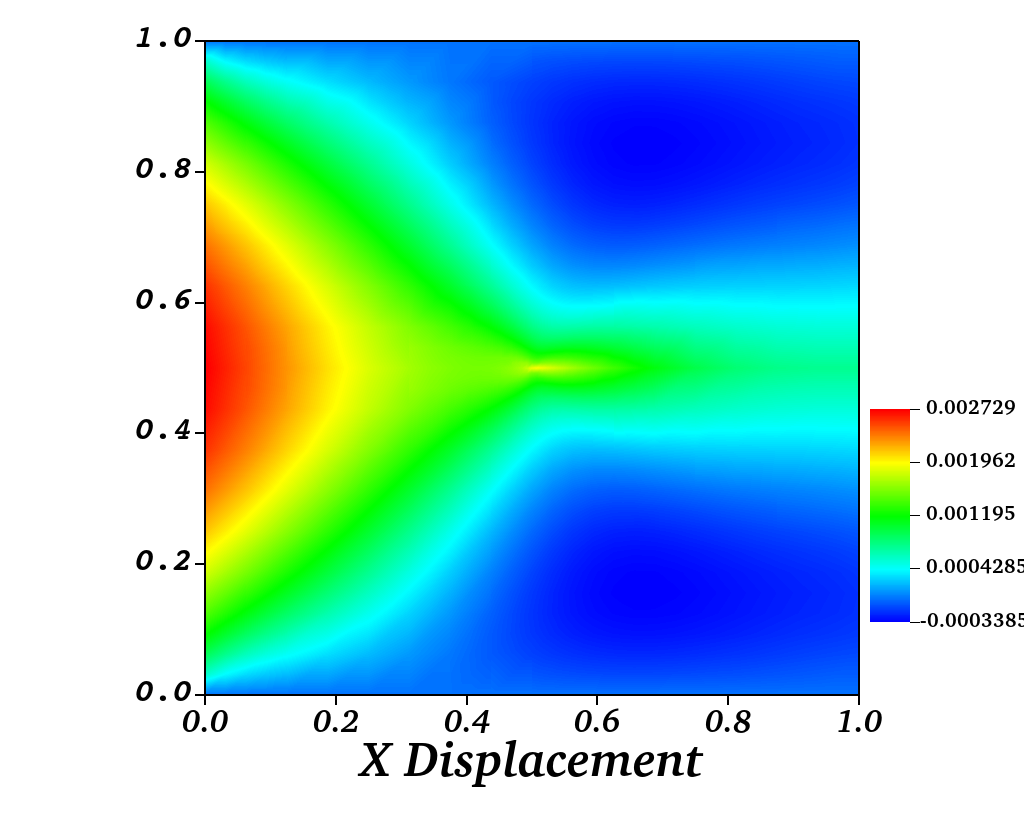}
        \caption{Horizontal displacement $u_1$.}
        \label{fig:UT_visit_ux}
    \end{subfigure}
    \hfill
    \begin{subfigure}[t]{0.48\textwidth}
        \centering
        \includegraphics[width=\textwidth]
        {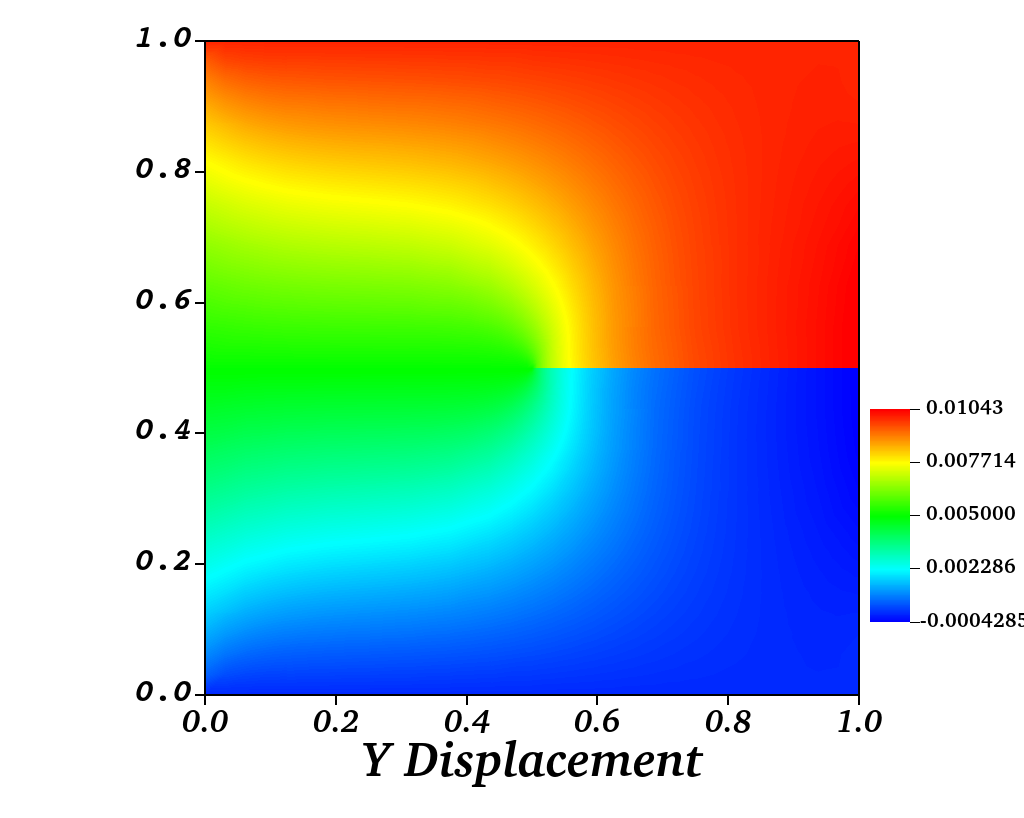}
        \caption{Vertical displacement $u_2$.}
        \label{fig:UT_visit_uy}
    \end{subfigure}

    \caption{Full-field distributions of stress, strain, and displacement for uniform tensile loading evaluated at a high nonlinearity threshold ($\beta=5$).}
    \label{fig:UT_visit}
\end{figure}

The full-field spatial distributions presented in Fig.~\ref{fig:UT_visit} corroborate this pronounced tip localization for a nonlinear state ($\beta=5$). Both the stress and strain fields exhibit intense, symmetric concentrations tightly confined around the crack tip, distinctly capturing the geometry of the fracture process zone. The corresponding displacement fields retain the global macroscopic character of the imposed tensile loading while accurately reflecting the kinematic discontinuities introduced by the crack. Specifically, the vertical displacement field ($u_2$) reveals a sharp jump across the opposing crack faces, indicative of a pure Mode I crack opening displacement under the applied load. Simultaneously, the horizontal displacement field ($u_1$) illustrates the macroscopic Poisson's effect—manifesting as inward lateral contraction—alongside a distinct spatial redistribution pattern as the material deforms and flows around the stress-relieved crack flanks.

\subsection{Uniform compressive loading}
\begin{figure}[H]
    \centering
    \begin{subfigure}[t]{0.32\linewidth}
        \centering
        \includegraphics[width=\linewidth]{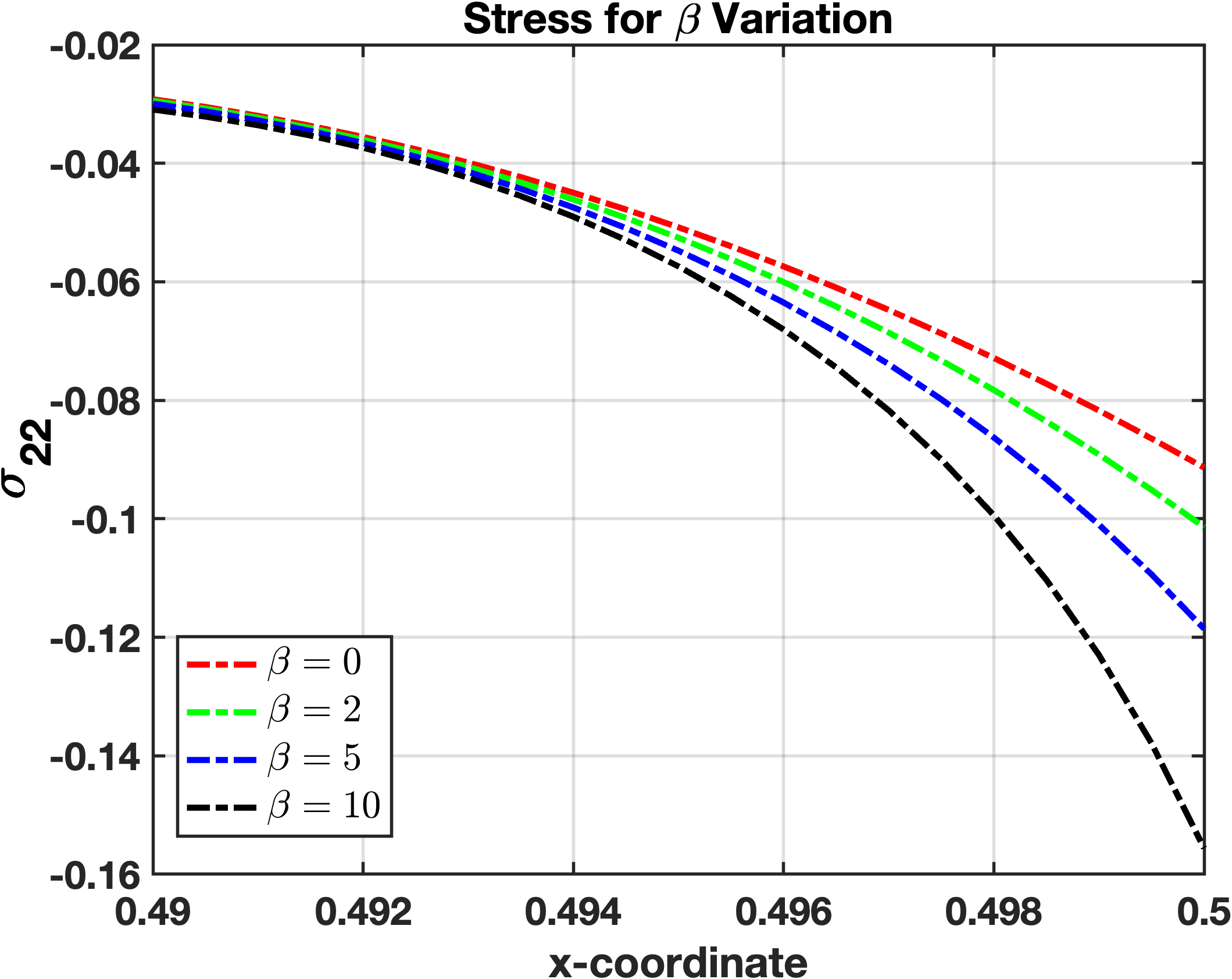}
        \caption{Normal stress $\sigma_{22}$.}
        \label{fig:UC_matlab_stress}
    \end{subfigure}%
    \hfill
    \begin{subfigure}[t]{0.32\linewidth}
        \centering
        \includegraphics[width=\linewidth]{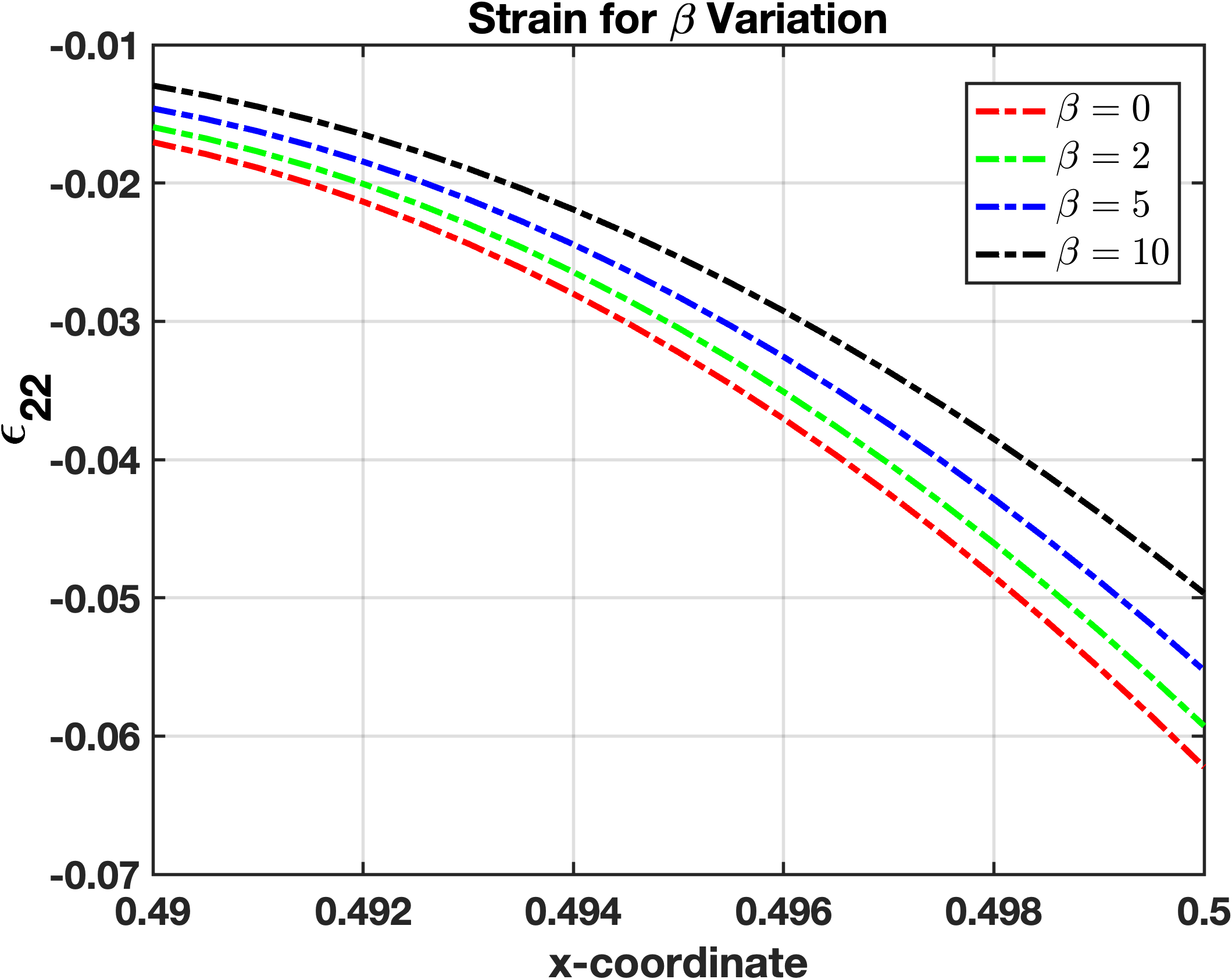}
        \caption{Normal strain $\varepsilon_{22}$.}
        \label{fig:UC_matlab_strain}
    \end{subfigure}%
    \hfill
    \begin{subfigure}[t]{0.32\linewidth}
        \centering
        \includegraphics[width=\linewidth]{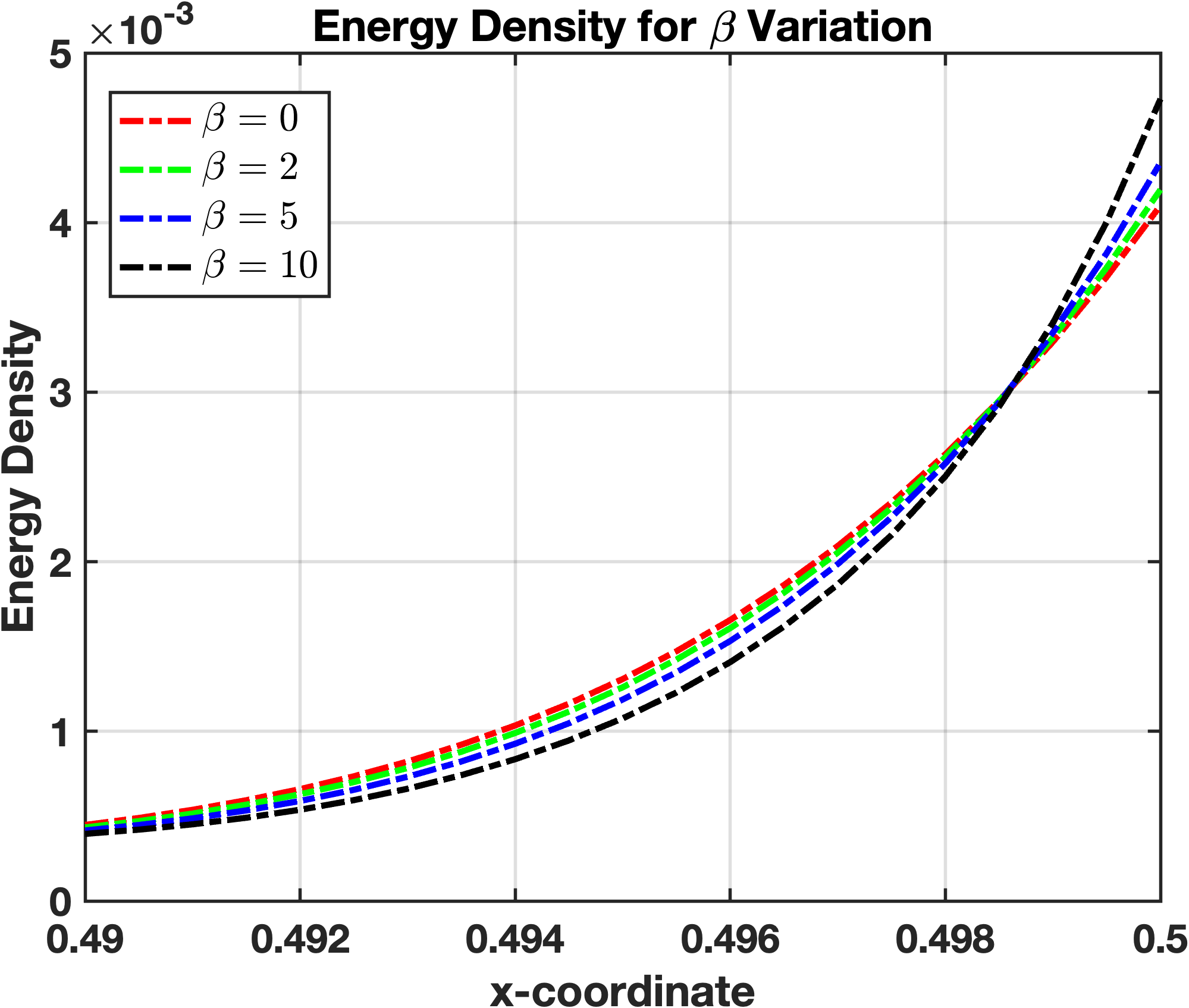}
        \caption{Strain Energy Density (SED).}
        \label{fig:UC_matlab_sed}
    \end{subfigure}
    
    \caption{Evolution of the normal stress $\sigma_{22}$, normal strain $\varepsilon_{22}$, and Strain Energy Density (SED) along the uncracked ligament approaching the crack tip. The plots demonstrate the influence of the nonlinearity parameter $\beta$ under uniform compressive loading.}
    \label{fig:UC_matlab}
\end{figure}

Under uniform compressive loading, the elastostatic fields ahead of the crack tip exhibit a profound sensitivity to the constitutive nonlinearity. As depicted in Figure~\ref{fig:UC_matlab}, both the normal stress $\sigma_{22}$ and normal strain $\varepsilon_{22}$ trace strictly negative trajectories along the sampled interval, with their magnitudes increasing sharply as $x \to 0.5$. The influence of the material parameter $\beta$ dictates a distinct kinematic divergence in this singular regime. Consistent with the strain-limiting framework, an increase in $\beta$ enforces a strict boundary on the material's local deformability; the magnitude of the compressive strain is heavily restricted near the tip. Consequently, to maintain static equilibrium under the applied load while accommodating this kinematic truncation, the continuum must sustain severe internal force concentrations. This is evidenced by the amplified compressive stress magnitudes observed for higher $\beta$ values. Furthermore, this constrained mechanical response alters the local thermodynamic landscape, driving a steep increase in the SED precisely at the singularity for highly nonlinear ($\beta=10$) configurations. 

\begin{figure}[H]
\centering
\begin{tikzpicture}[x=4.2cm,y=4.2cm]
  \draw[thick] (0,0) rectangle (1,1);
  \draw[red,line width=1.5pt] (0.5,0.5)--(1,0.5);
  \fill[red] (0.5,0.5) circle[radius=0.012];
  \fill[pattern=north east lines] (0,-0.08) rectangle (1,0);

  \foreach \x in {0,0.1,...,1.1}
    \draw[-{Latex[length=2mm]},thick,blue] (\x,1)--++(0,-0.18);

  \node[blue,above] at (0.5,1.18) {$ u_1=0, \quad u_2= - \, \bar u$};
  \node[below] at (0.5,-0.08) {$\mathbf u=\mathbf 0$};
\end{tikzpicture}
\caption{Uniform compressive boundary condition.}
\label{fig:ddm-schematic-ii}
\end{figure}
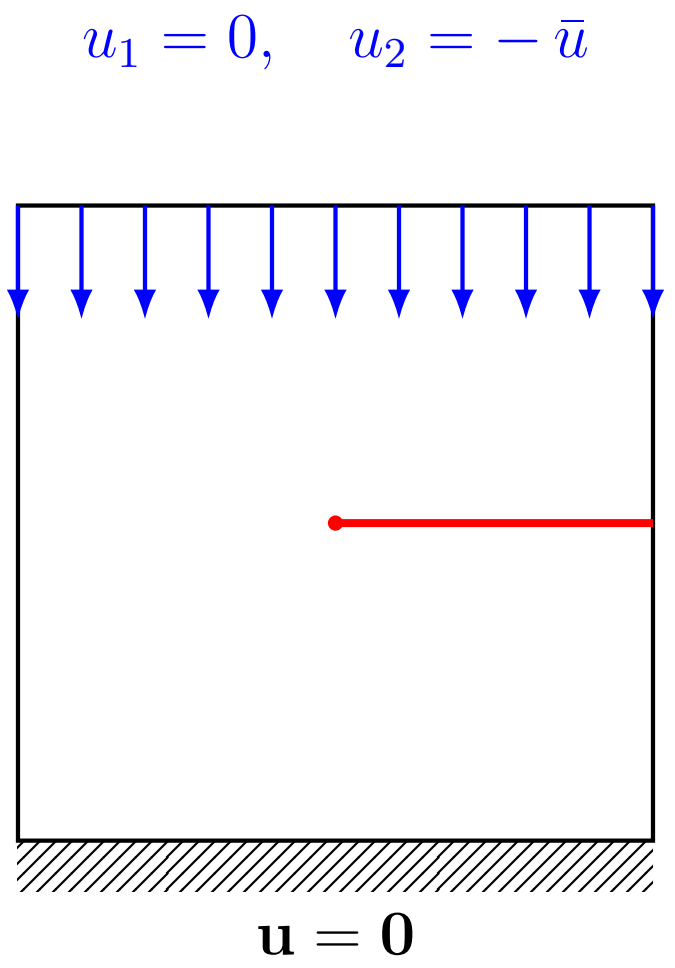 

\begin{figure}[H]
    \centering
    \begin{subfigure}[t]{0.48\textwidth}
        \centering
        \includegraphics[width=\textwidth]{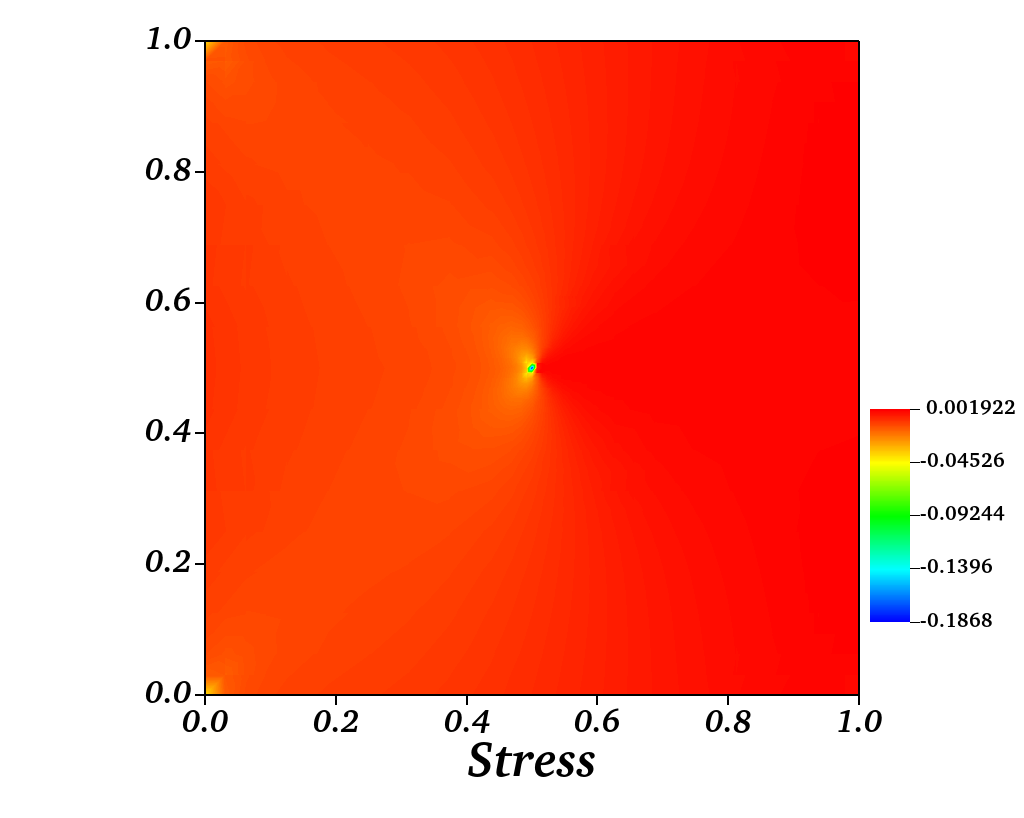}
        \caption{Stress field.}
        \label{fig:UC_visit_stress}
    \end{subfigure}
    \hfill
    \begin{subfigure}[t]{0.48\textwidth}
        \centering
        \includegraphics[width=\textwidth]{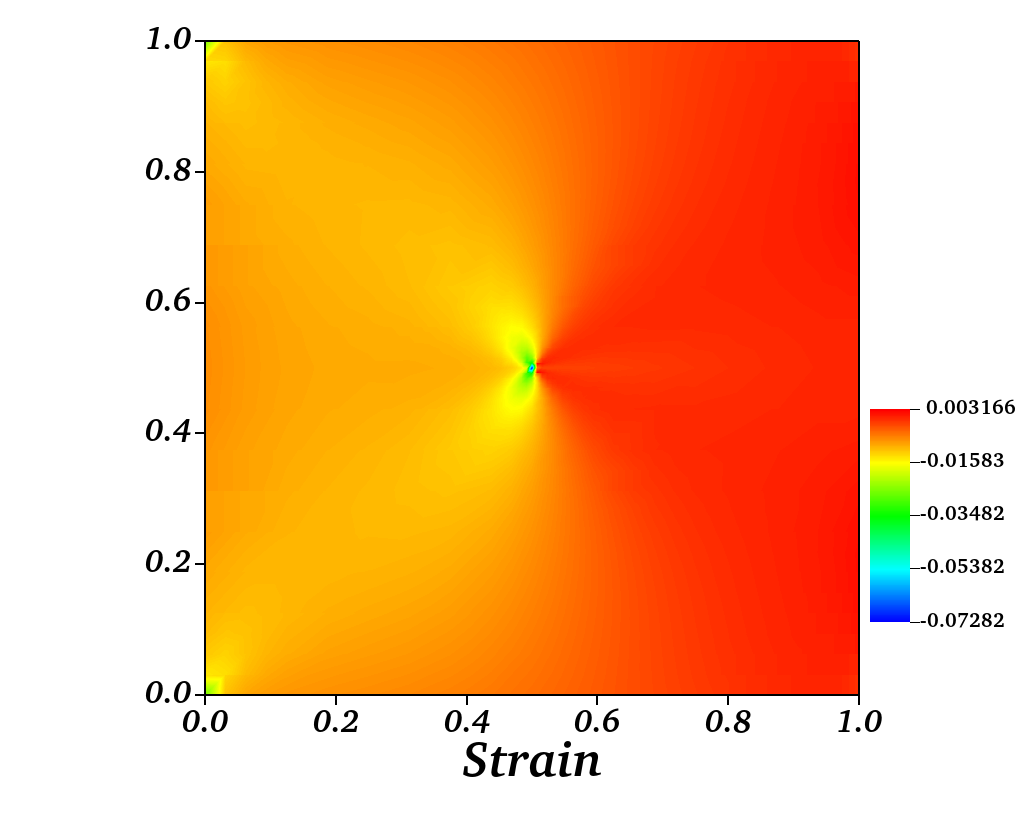}
        \caption{Strain field.}
        \label{fig:UC_visit_strain}
    \end{subfigure}

    \vspace{0.3cm}

    \begin{subfigure}[t]{0.48\textwidth}
        \centering
        \includegraphics[width=\textwidth]{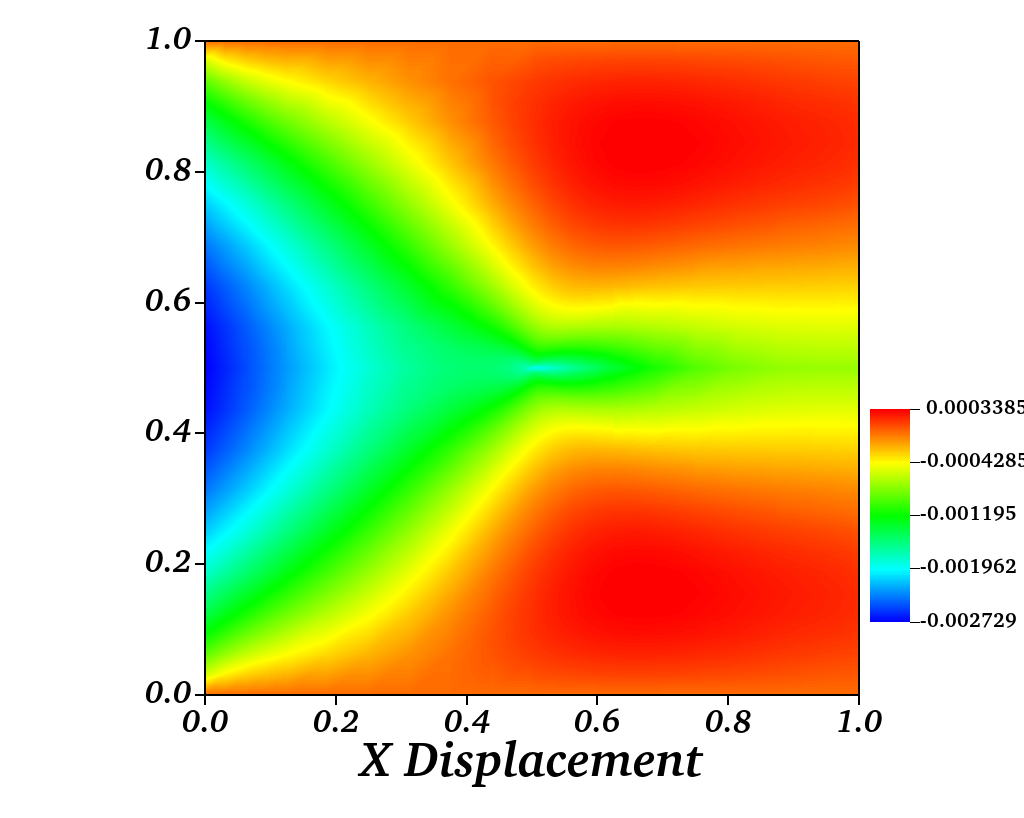}
        \caption{Horizontal displacement $u_1$.}
        \label{fig:UC_visit_ux}
    \end{subfigure}
    \hfill
    \begin{subfigure}[t]{0.48\textwidth}
        \centering
        \includegraphics[width=\textwidth]{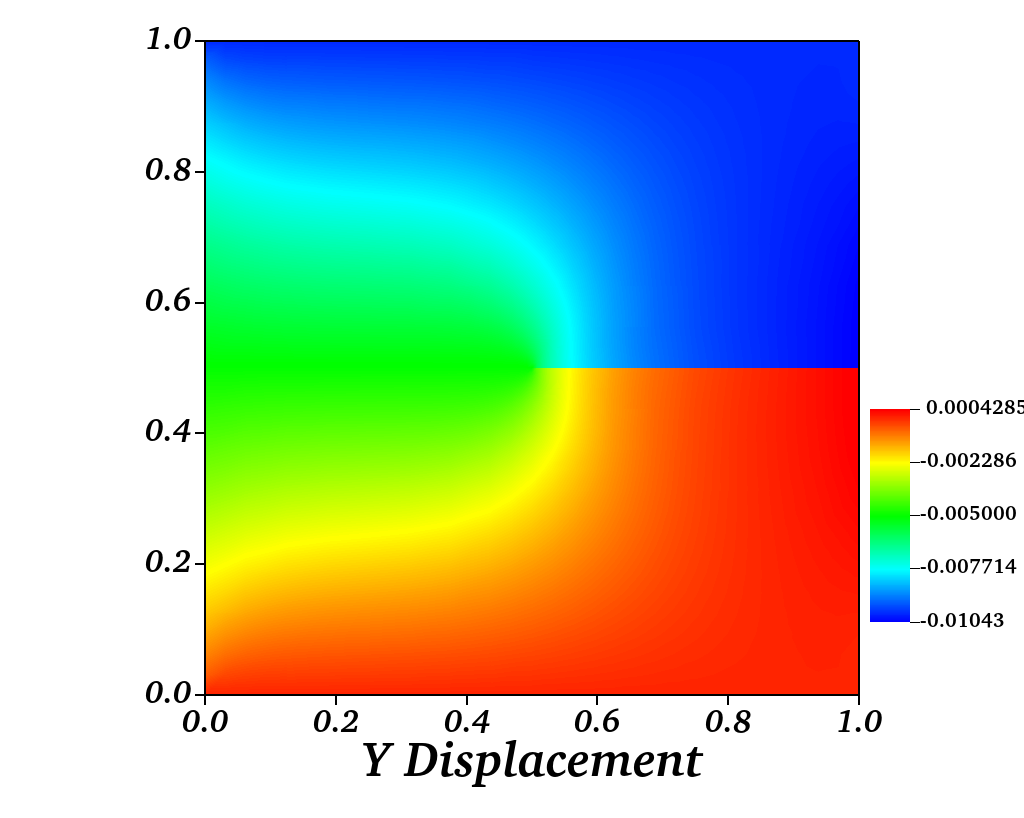}
        \caption{Vertical displacement $u_2$.}
        \label{fig:UC_visit_uy}
    \end{subfigure}

    \caption{Two-dimensional spatial distributions of the stress, strain, and displacement fields for uniform compressive loading evaluated at $\beta=5$.}
    \label{fig:UC_visit}
\end{figure}

Figure~\ref{fig:UC_visit} delineates the full two-dimensional spatial distribution of the deformation fields, confirming that the high-gradient fracture process zone remains strongly governed by the tip geometry. The stress and strain contours exhibit localized, lobular concentrations characteristic of crack-tip singular fields, perfectly mirroring the spatial topology of the tensile case but with an inverted polarity. 

From a kinematic perspective, the displacement fields offer critical insight into the continuum response. The horizontal displacement field $u_1$ captures the expected lateral material expansion (Poisson effect) induced by the primary vertical compression. Simultaneously, the vertical displacement field $u_2$ illustrates the physical convergence of the upper and lower crack flanks. It is mathematically rigorous to note that because the current boundary value problem is formulated as a variational equality subject to free-surface Neumann conditions, it does not explicitly enforce non-penetration constraints (i.e., Signorini unilateral contact conditions). As a result, the deep compressive loads induce a theoretical geometric interpenetration of the crack faces. Despite this idealized crack closure mechanics, the numerical scheme flawlessly resolves the localized strain-limiting phenomena without spurious oscillations, preserving the structural symmetry and stability of the density-dependent constitutive model.

\subsection{Piecewise-slope tensile loading}
\begin{figure}[H]
    \centering
    \begin{subfigure}[t]{0.32\linewidth}
        \centering
        \includegraphics[width=\linewidth]{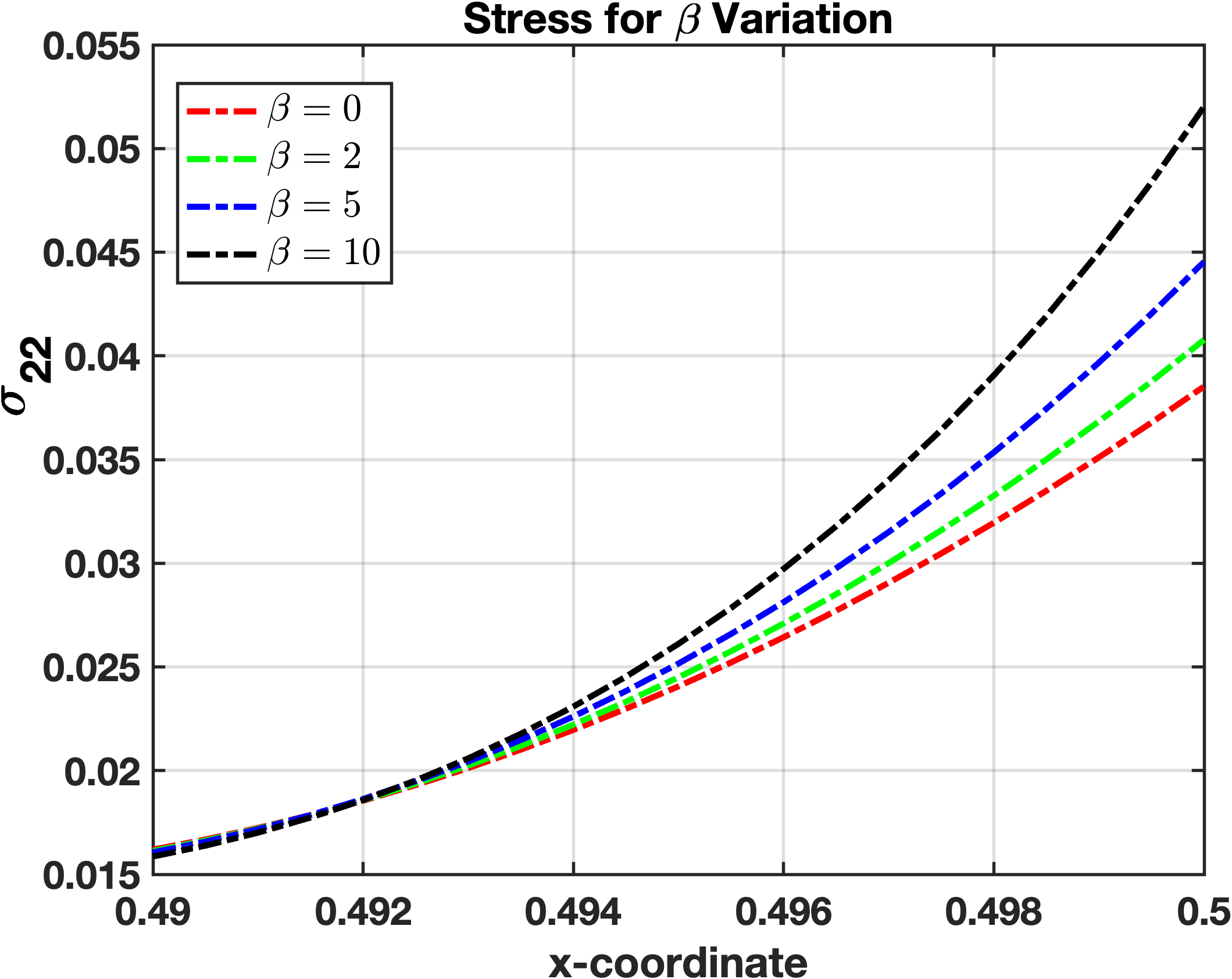}
        \caption{Normal stress $\sigma_{22}$.}
        \label{fig:PT_matlab_stress}
    \end{subfigure}%
    \hfill
    \begin{subfigure}[t]{0.32\linewidth}
        \centering
        \includegraphics[width=\linewidth]{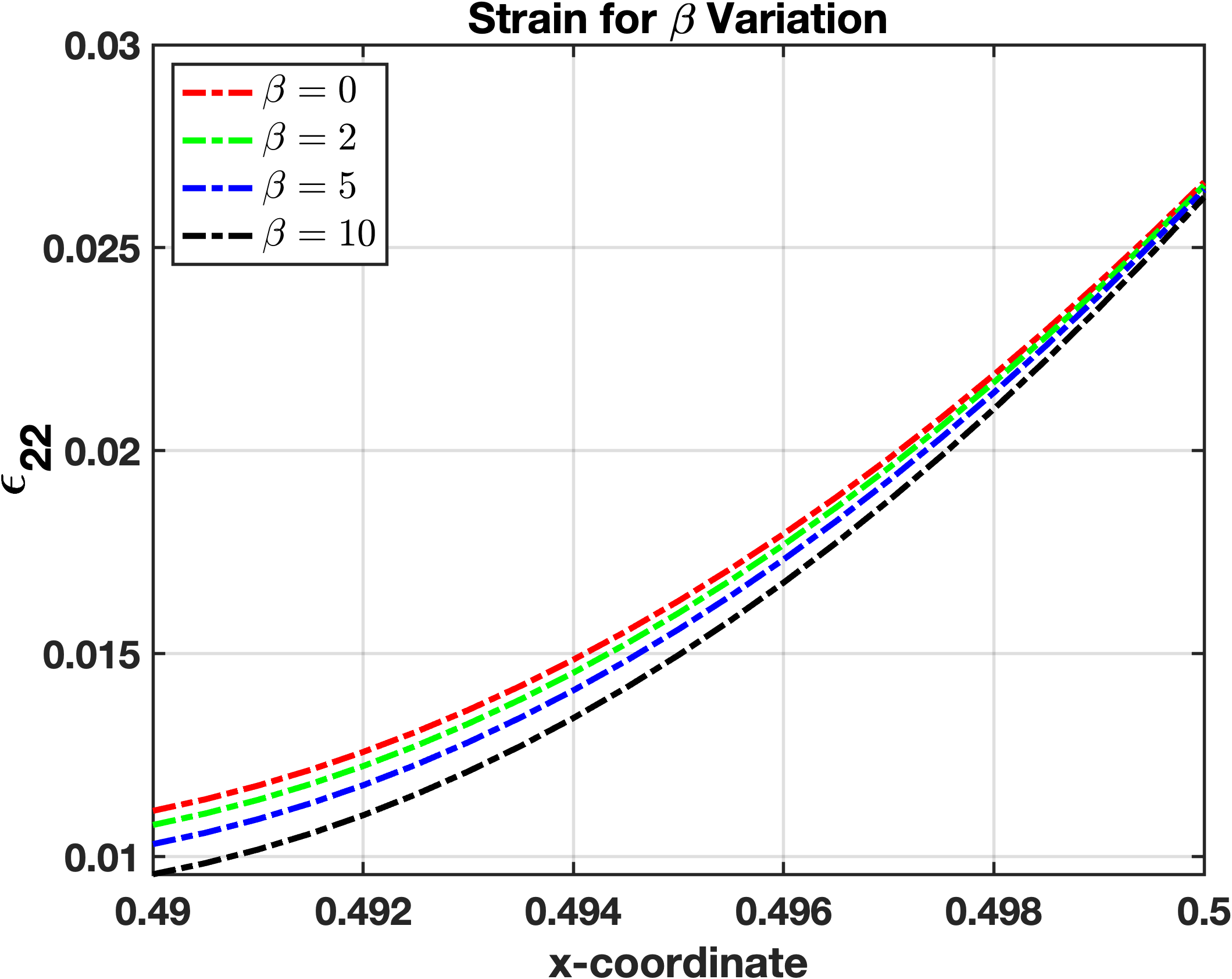}
        \caption{Normal strain $\varepsilon_{22}$.}
        \label{fig:PT_matlab_strain}
    \end{subfigure}%
    \hfill
    \begin{subfigure}[t]{0.32\linewidth}
        \centering
        \includegraphics[width=\linewidth]{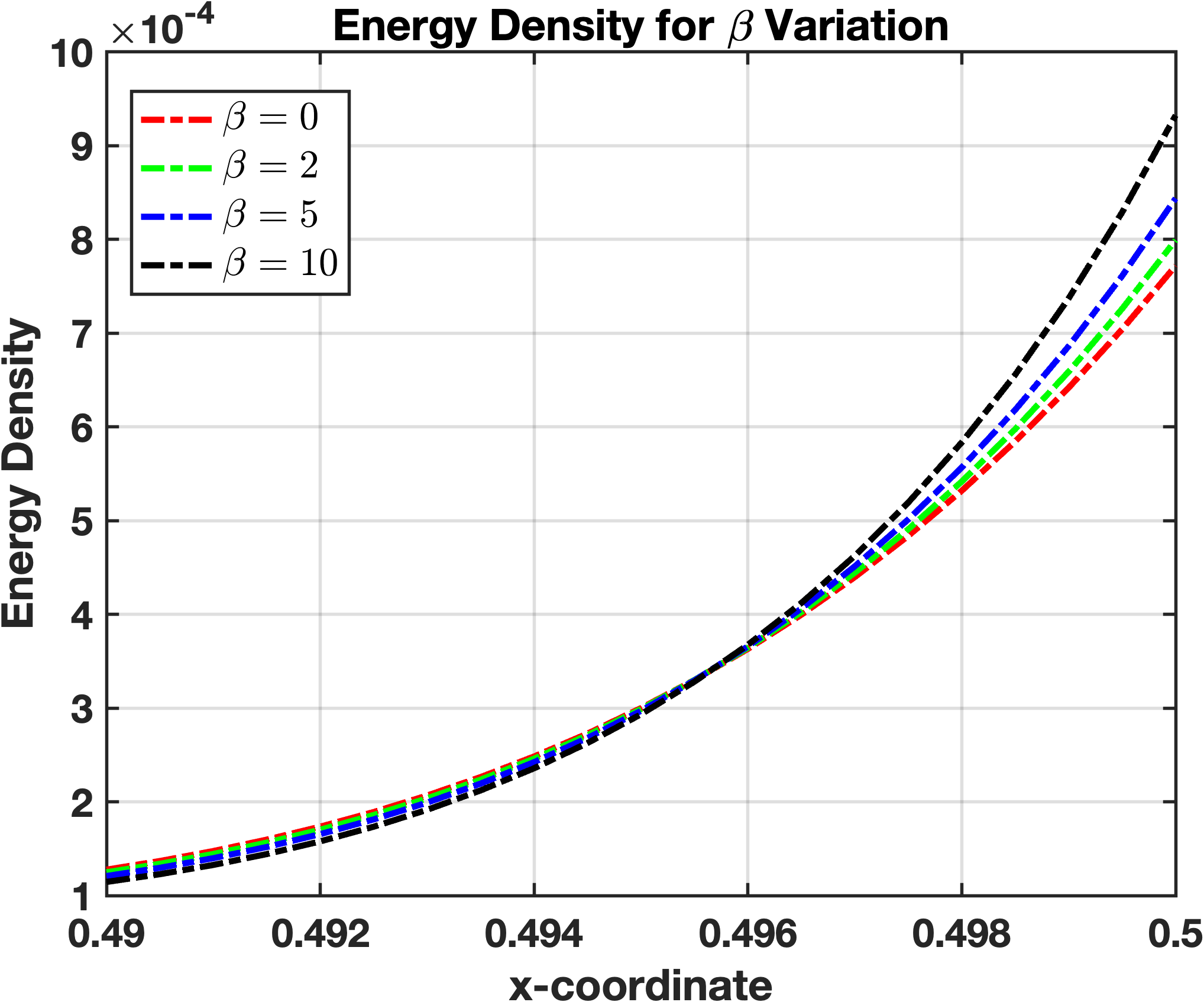}
        \caption{Strain Energy Density (SED).}
        \label{fig:PT_matlab_sed}
    \end{subfigure}
    
    \caption{Evolution of the normal stress $\sigma_{22}$, normal strain $\varepsilon_{22}$, and Strain Energy Density (SED) along the uncracked ligament approaching the crack tip. The plots demonstrate the influence of the nonlinearity parameter $\beta$ under piecewise-slope tensile loading.}
    \label{fig:PT_matlab}
\end{figure}

The imposition of a piecewise-slope tensile loading introduces a non-uniform macroscopic traction profile, fundamentally altering how far-field loads are transmitted to the fracture process zone. As evaluated along the uncracked ligament ($x \to 0.5$), the near-tip elastostatic fields reveal a pronounced decoupling in their sensitivity to the constitutive nonlinearity. As illustrated in Figure~\ref{fig:PT_matlab}, increasing the strain-limiting parameter $\beta$ triggers a severe amplification in the normal stress $\sigma_{22}$, with the stress trajectories diverging significantly near the singularity. In stark contrast, the normal strain profiles ($\varepsilon_{22}$) remain tightly clustered and nearly coalesce at the crack tip, displaying minimal separation regardless of the $\beta$ magnitude. This kinematic bounding confirms that the material's local deformability is strictly governed by the implicit constitutive relation. Because the material cannot freely stretch to accommodate the loading gradients, the continuum is forced to balance the boundary tractions through highly localized force concentrations. Consequently, this stress-dominant response dictates the local thermodynamics, culminating in a steep exponential surge in the SED precisely at the crack tip for highly nonlinear ($\beta=10$) configurations.

\begin{figure}[H]
\centering
\begin{tikzpicture}[x=4.2cm,y=4.2cm]
  \draw[thick] (0,0) rectangle (1,1);
  \draw[red,line width=1.5pt] (0.5,0.5)--(1,0.5);
  \fill[red] (0.5,0.5) circle[radius=0.012];
  \fill[pattern=north east lines] (0,-0.08) rectangle (1,0);

\foreach \x in {0.1,0.2,...,0.9} {
  \pgfmathsetmacro{\h}{0.36*(1-abs(2*\x-1))}
  \draw[-{Latex[length=2mm]},thick,blue] (\x,1)--++(0,\h);
}

\fill[blue] (-0.018,1) -- (0.018,1) -- (0,1.035) -- cycle;
\fill[blue] (0.982,1) -- (1.018,1) -- (1,1.035) -- cycle;

  \node[blue,above] at (0.5,1.36) {$u_1=0, \quad u_2=2\bar u\min(x,1-x)$};
  \node[below] at (0.5,-0.08) {$\mathbf u=\mathbf 0$};
\end{tikzpicture}
\caption{Piecewise-slope-up boundary condition.}
\label{fig:ddm-schematic-iii}
\end{figure}
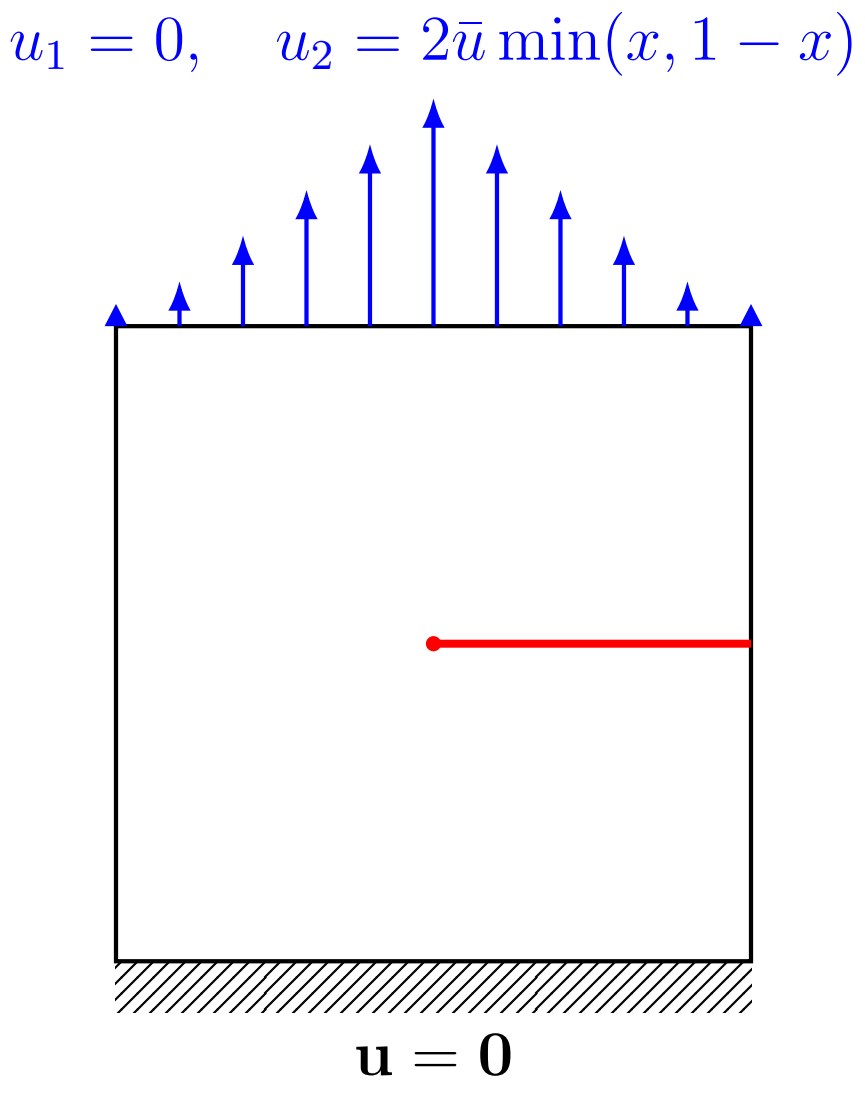 

\begin{figure}[H]
    \centering
    \begin{subfigure}[t]{0.48\textwidth}
        \centering
        \includegraphics[width=\textwidth]{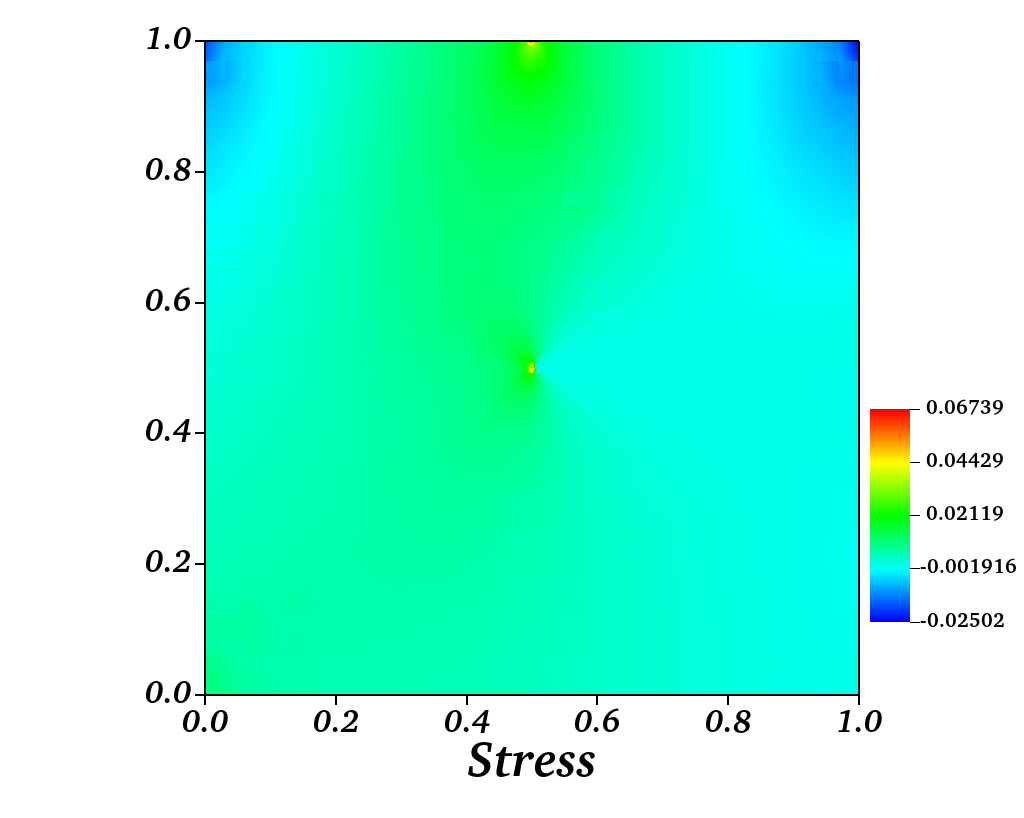}
        \caption{Stress field.}
        \label{fig:PT_visit_stress}
    \end{subfigure}
    \hfill
    \begin{subfigure}[t]{0.48\textwidth}
        \centering
        \includegraphics[width=\textwidth]{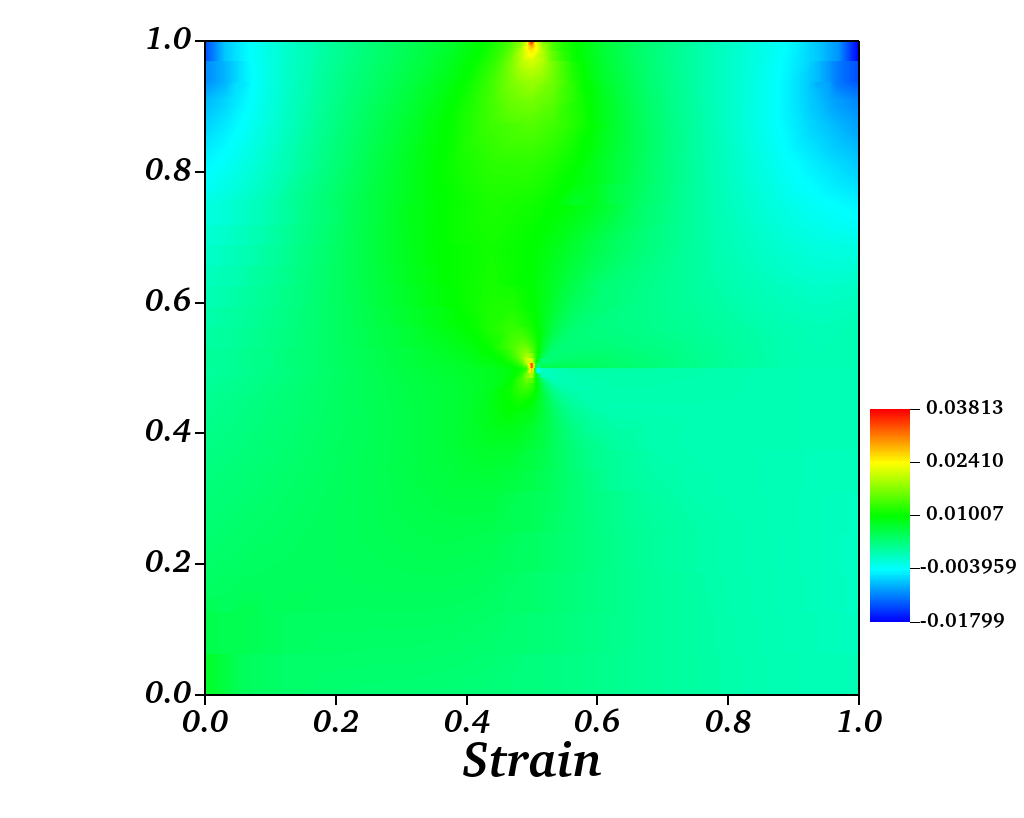}
        \caption{Strain field.}
        \label{fig:PT_visit_strain}
    \end{subfigure}

    \vspace{0.3cm}

    \begin{subfigure}[t]{0.48\textwidth}
        \centering
        \includegraphics[width=\textwidth]{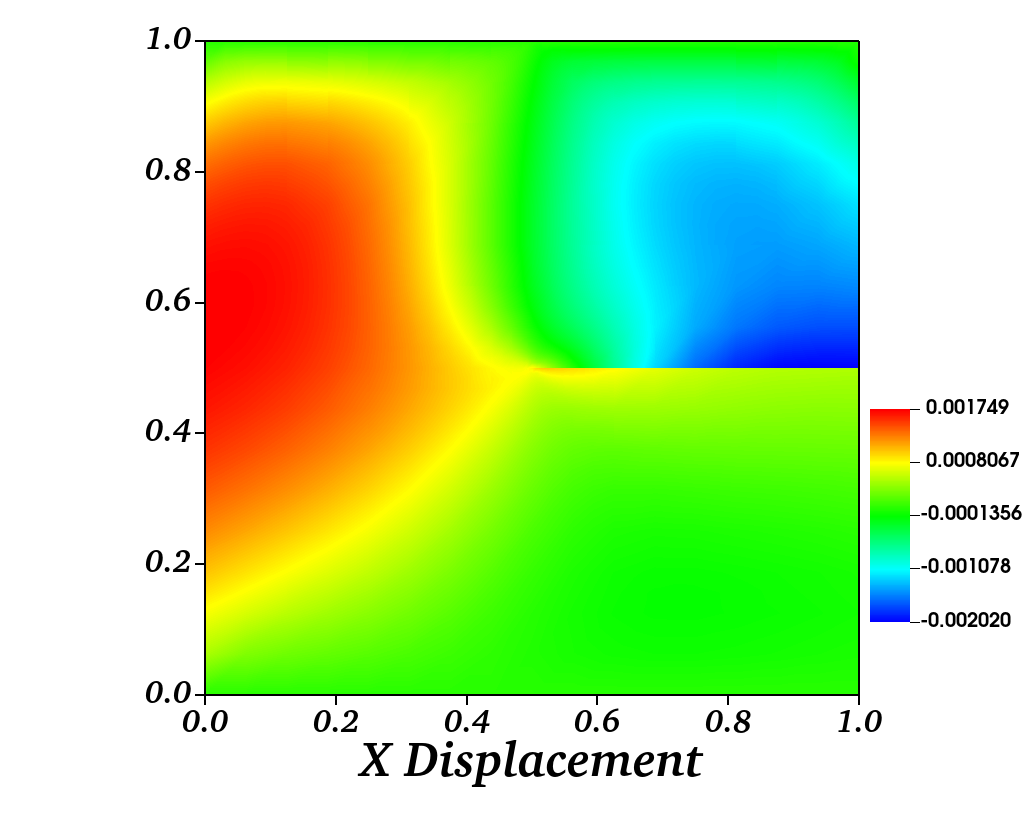}
        \caption{Horizontal displacement $u_1$.}
        \label{fig:PT_visit_ux}
    \end{subfigure}
    \hfill
    \begin{subfigure}[t]{0.48\textwidth}
        \centering
        \includegraphics[width=\textwidth]{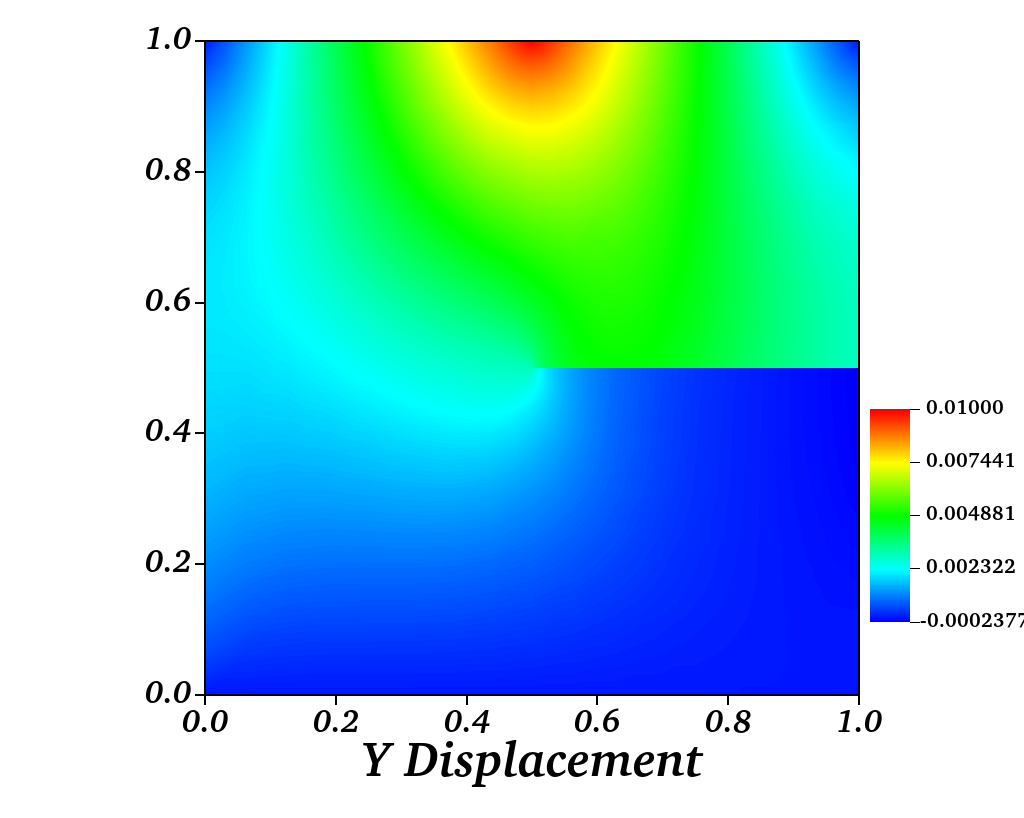}
        \caption{Vertical displacement $u_2$.}
        \label{fig:PT_visit_uy}
    \end{subfigure}

    \caption{Two-dimensional spatial distributions of the stress, strain, and displacement fields for piecewise-slope tensile loading evaluated at $\beta=5$.}
    \label{fig:PT_visit}
\end{figure}

The two-dimensional spatial distributions presented in Figure~\ref{fig:PT_visit} provide critical insights into the competition between global boundary conditions and local geometric singularities. The macroscopic influence of the piecewise-slope profile is immediately evident in the displacement fields: the vertical displacement $u_2$ exhibits its maximum structural deflection centrally along the top boundary, directly mirroring the applied non-uniform traction. Similarly, the lateral Poisson contraction mapped in the horizontal displacement $u_1$ displays an asymmetric topological shift compared to uniform loading regimes.

Despite these global kinematic distortions, the intense gradients of the fracture process zone remain decisively governed by the crack-tip geometry. The spatial contours of both the stress and strain fields show a high localization surrounding the singularity at $(0.5, 0.5)$. Furthermore, the vertical displacement field captures the discontinuous mathematical jump across the crack ($y=0.5$), confirming that the adaptive discretization scheme robustly resolves the intricate mixed influence of the non-uniform boundary profile and the density-dependent constitutive singularity without sacrificing structural symmetry.

\subsection{Piecewise-slope compressive loading}

\begin{figure}[H]
    \centering
    \begin{subfigure}[t]{0.32\linewidth}
        \centering
        \includegraphics[width=\linewidth]{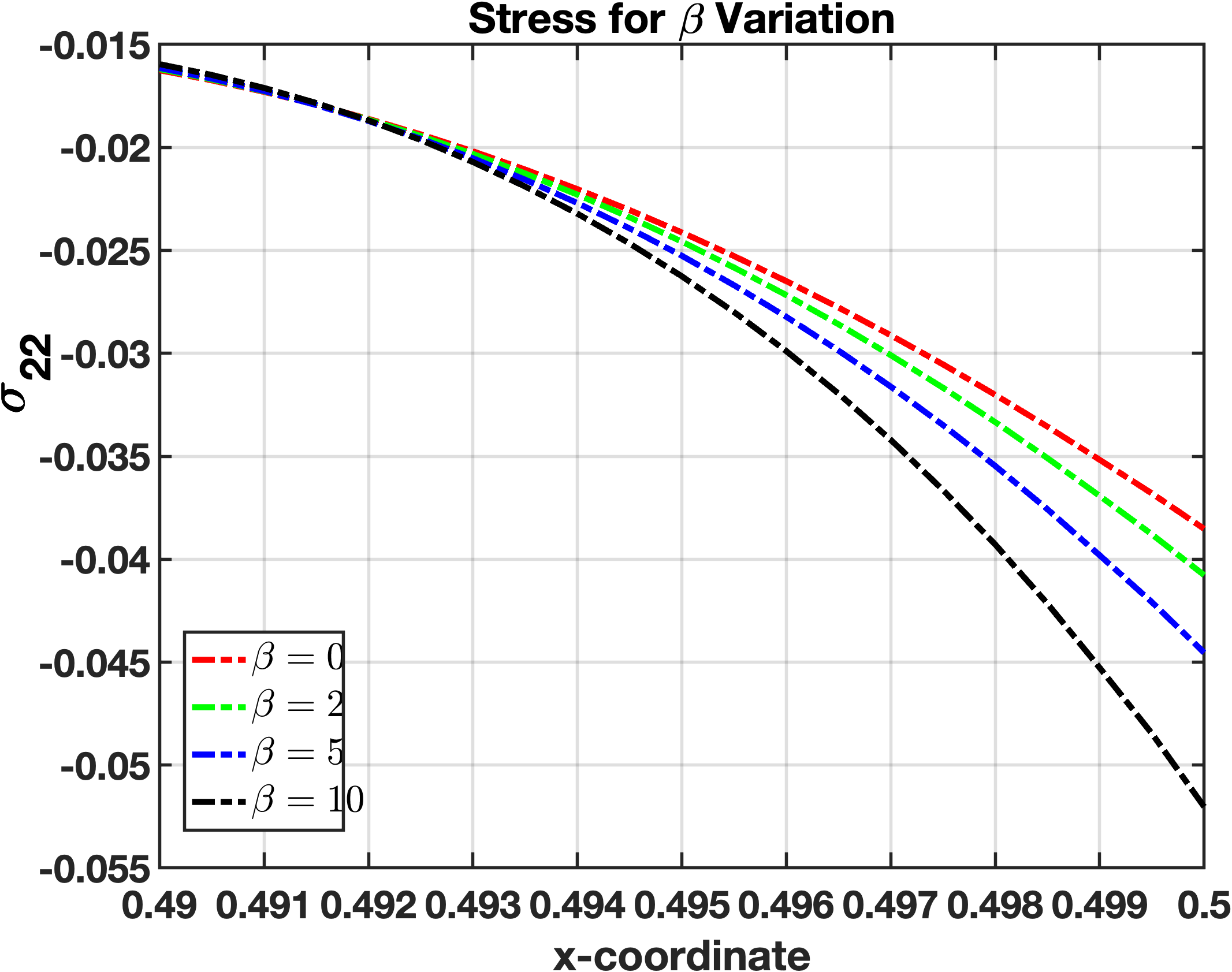}
        \caption{Normal stress $\sigma_{22}$.}
        \label{fig:PC_matlab_stress}
    \end{subfigure}%
    \hfill
    \begin{subfigure}[t]{0.32\linewidth}
        \centering
        \includegraphics[width=\linewidth]{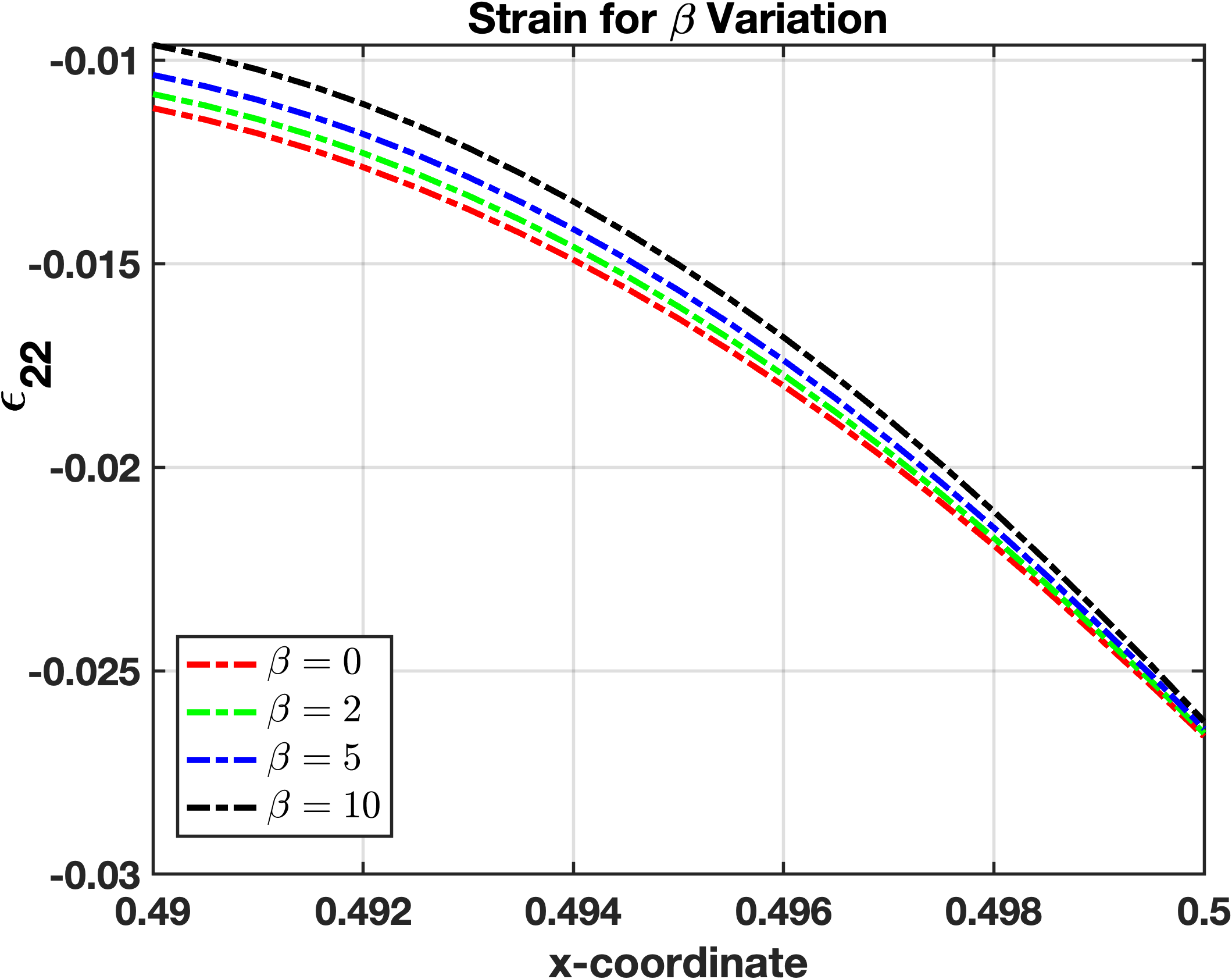}
        \caption{Normal strain $\varepsilon_{22}$.}
        \label{fig:PC_matlab_strain}
    \end{subfigure}%
    \hfill
    \begin{subfigure}[t]{0.32\linewidth}
        \centering
        \includegraphics[width=\linewidth]{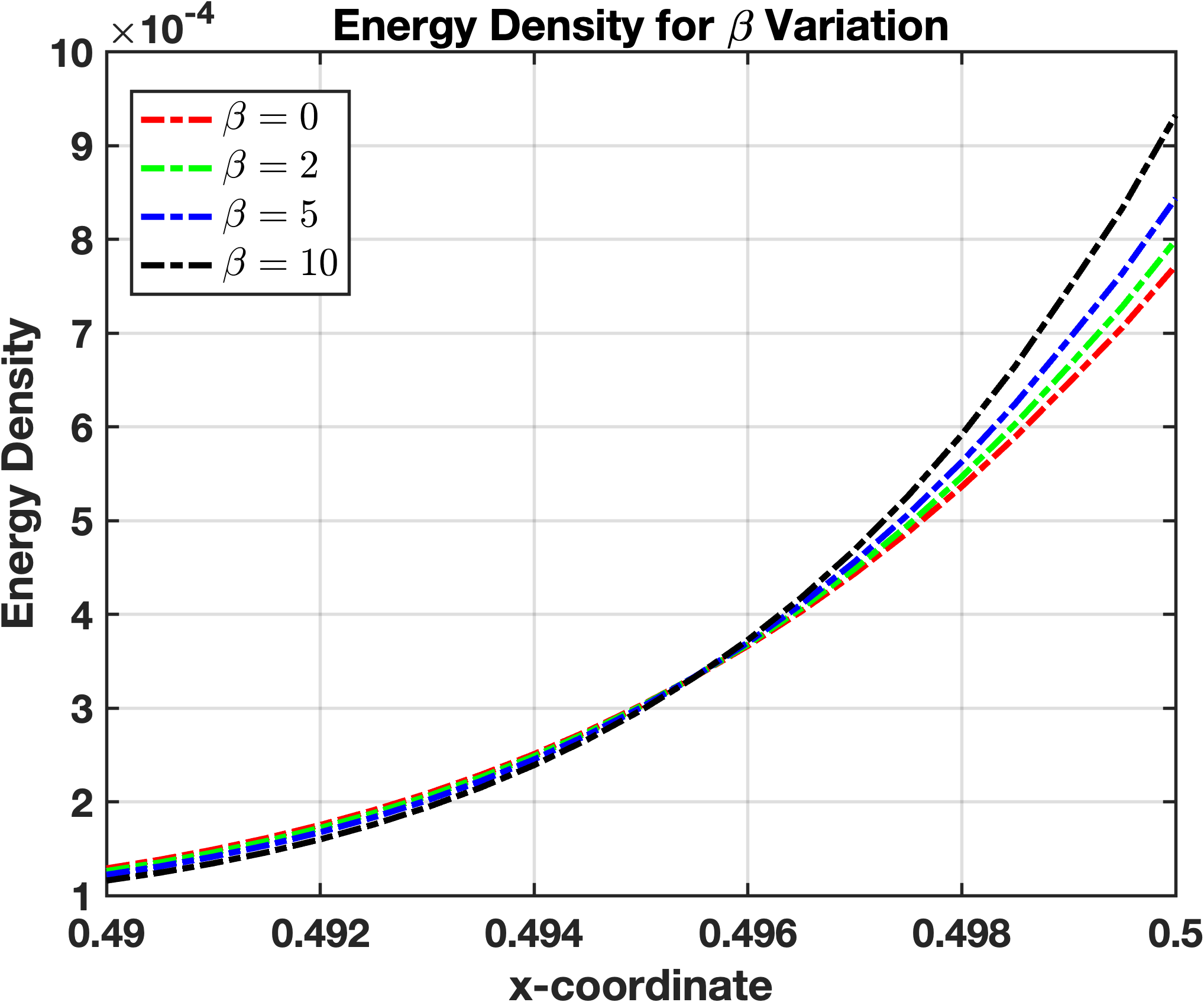}
        \caption{Strain Energy Density (SED).}
        \label{fig:PC_matlab_sed}
    \end{subfigure}
    
    \caption{Evolution of the normal stress $\sigma_{22}$, normal strain $\varepsilon_{22}$, and Strain Energy Density (SED) along the uncracked ligament approaching the crack tip. The plots demonstrate the influence of the nonlinearity parameter $\beta$ under a non-uniform, piecewise-slope compressive load.}
    \label{fig:PC_matlab}
\end{figure}

The application of a piecewise-slope compressive load yields a profound localization of near-tip fields that mirror the macroscopic non-uniformity of the boundary tractions. As observed along the uncracked ligament in Figure~\ref{fig:PC_matlab}, both the normal stress $\sigma_{22}$ and normal strain $\varepsilon_{22}$ exhibit strictly negative trajectories. The constitutive nonlinearity, governed by $\beta$, enforces a strict kinematic restriction on the material's compressive deformability. As $\beta$ increases, the magnitude of the compressive strain is heavily bounded, forcing the strain curves to cluster tightly as they approach the singularity. To sustain the localized equilibrium under this restricted deformation mapping, the material must develop severe internal compressive stress concentrations. This is clearly evidenced by the stark divergence and amplification of the $\sigma_{22}$ magnitudes for higher $\beta$ values near $x=0.5$ \cite{YoonMallikarjunaiahBhatta2024}. Consequently, this highly constrained, stress-dominant mechanical response drives a rapid exponential increase in the local SED precisely at the crack tip, particularly for the highly nonlinear ($\beta=10$) state.

\begin{figure}[H]
\centering
\begin{tikzpicture}[x=4.2cm,y=4.2cm]
  \draw[thick] (0,0) rectangle (1,1);
  \draw[red,line width=1.5pt] (0.5,0.5)--(1,0.5);
  \fill[red] (0.5,0.5) circle[radius=0.012];
  \fill[pattern=north east lines] (0,-0.08) rectangle (1,0);

  \foreach \x in {0.1,0.2,...,0.9} {
    \pgfmathsetmacro{\h}{0.36*(1-abs(2*\x-1))}
    \draw[-{Latex[length=2mm]},thick,blue] (\x,1)--++(0,-\h);
  }

  \fill[blue] (-0.018,1) -- (0.018,1) -- (0,0.965) -- cycle;
  \fill[blue] (0.982,1) -- (1.018,1) -- (1,0.965) -- cycle;

  \node[blue,above] at (0.5,1.08) {$u_1=0, \quad u_2=-2\bar u\min(x,1-x)$};
  \node[below] at (0.5,-0.08) {$\mathbf u=\mathbf 0$};
\end{tikzpicture}
\caption{Piecewise-slope-down boundary condition.}
\label{fig:ddm-schematic-iv}
\end{figure}
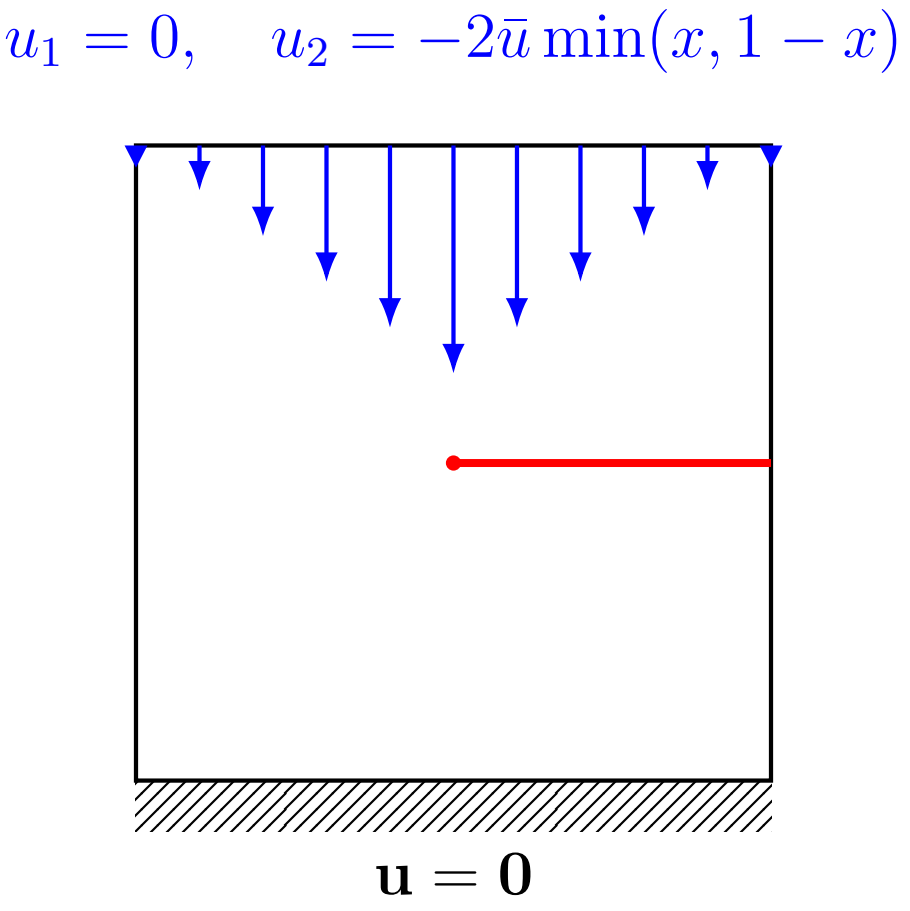 

\begin{figure}[H]
    \centering
    \begin{subfigure}[t]{0.48\textwidth}
        \centering
        \includegraphics[width=\textwidth]{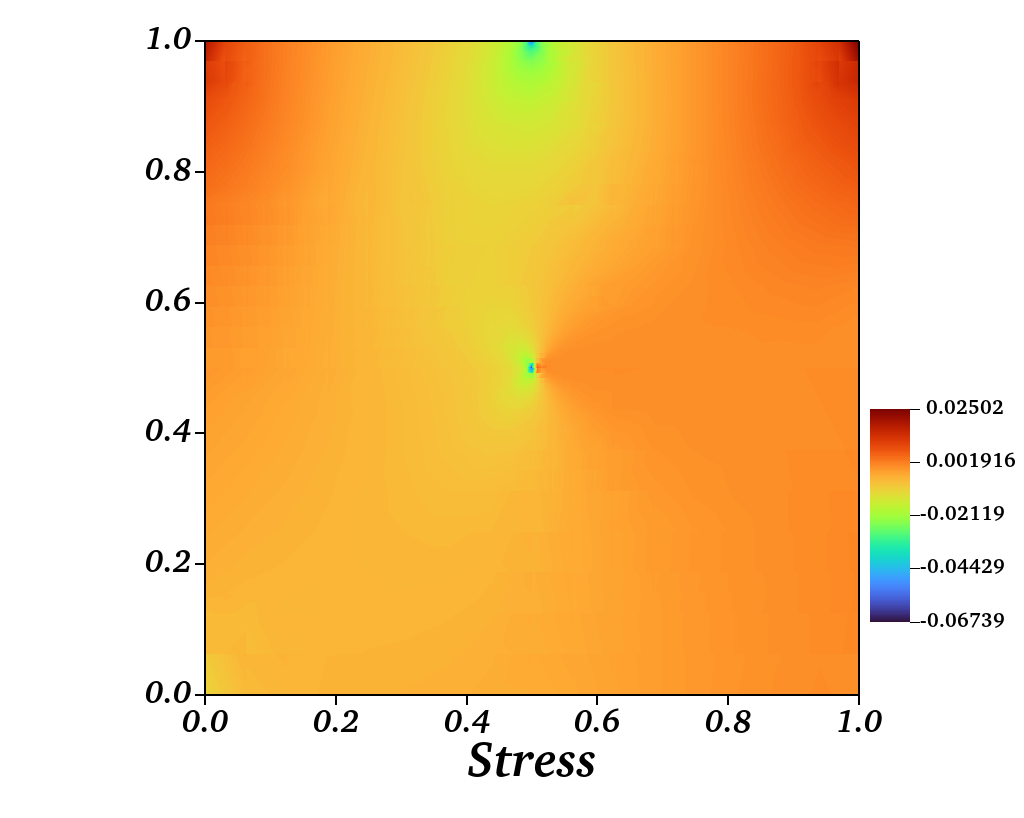}
        \caption{Stress field.}
        \label{fig:PC_visit_stress}
    \end{subfigure}
    \hfill
    \begin{subfigure}[t]{0.48\textwidth}
        \centering
        \includegraphics[width=\textwidth]{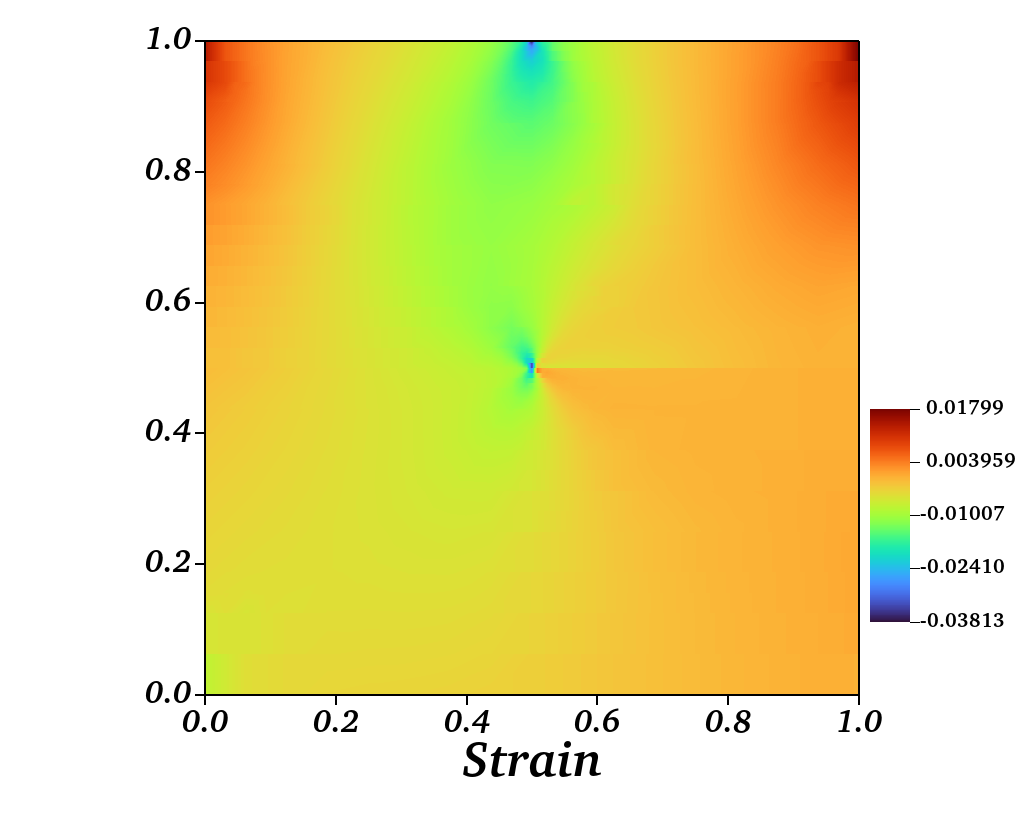}
        \caption{Strain field.}
        \label{fig:PC_visit_strain}
    \end{subfigure}

    \vspace{0.3cm}

    \begin{subfigure}[t]{0.48\textwidth}
        \centering
        \includegraphics[width=\textwidth]{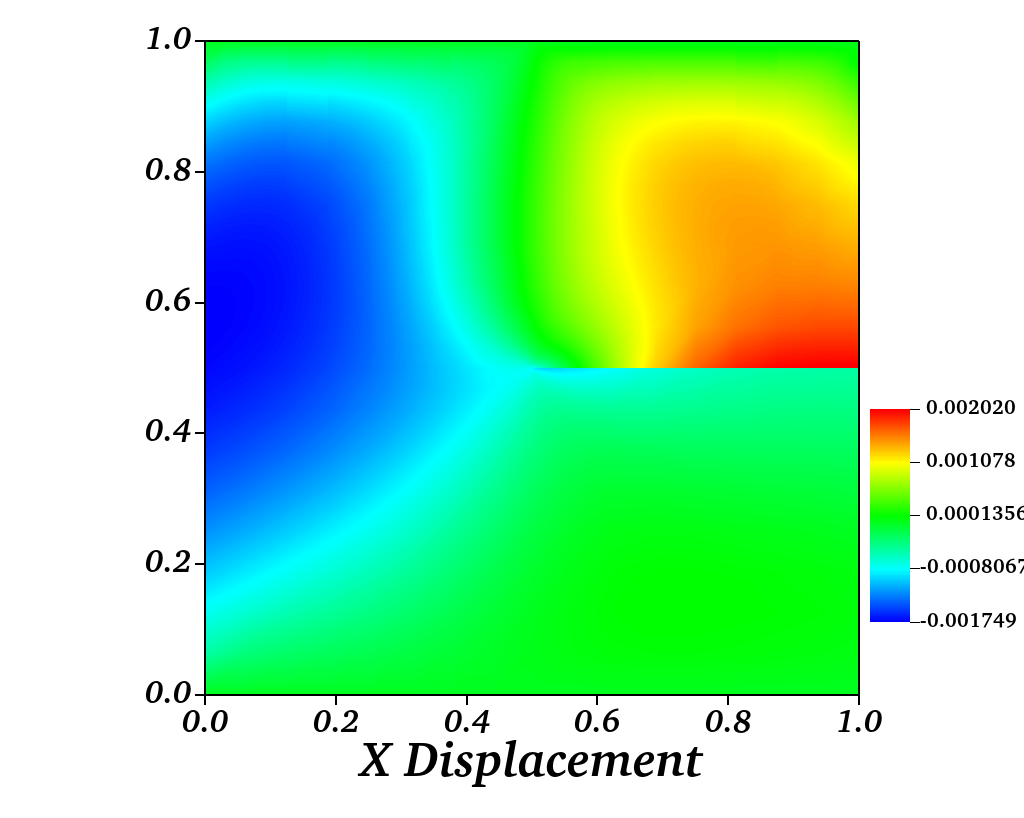}
        \caption{Horizontal displacement $u_1$.}
        \label{fig:PC_visit_ux}
    \end{subfigure}
    \hfill
    \begin{subfigure}[t]{0.48\textwidth}
        \centering
        \includegraphics[width=\textwidth]{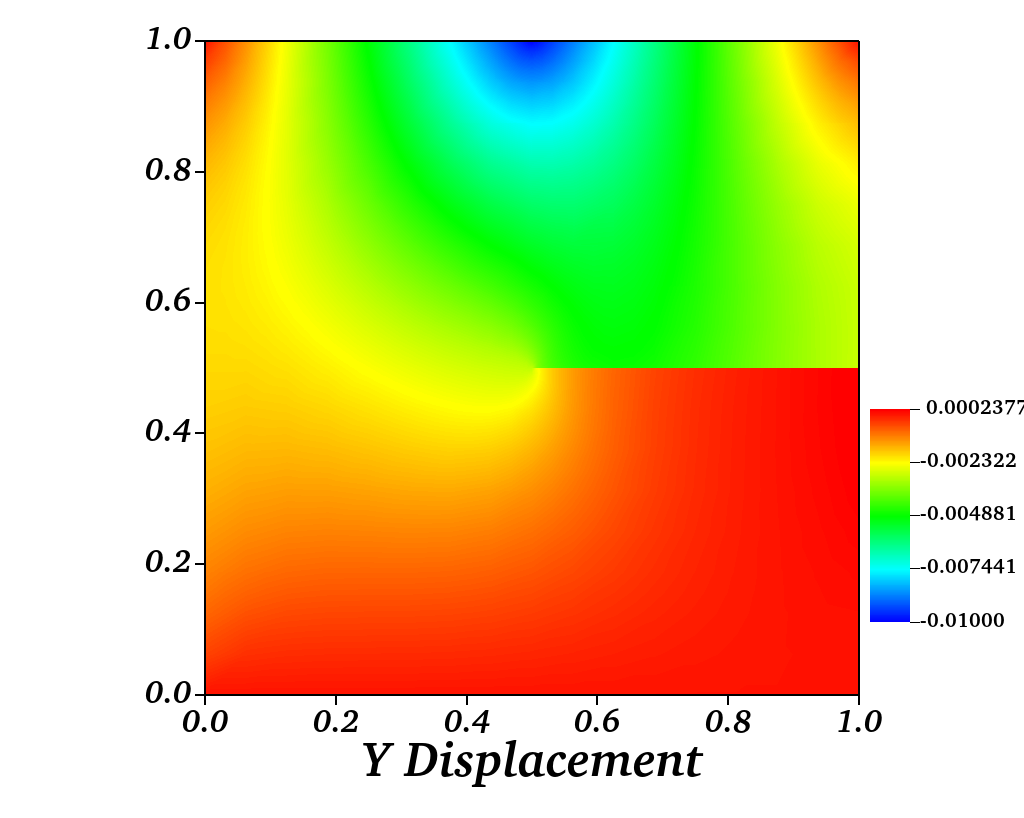}
        \caption{Vertical displacement $u_2$.}
        \label{fig:PC_visit_uy}
    \end{subfigure}

    \caption{Two-dimensional spatial distributions of the stress, strain, and displacement fields for piecewise-slope compressive loading evaluated at $\beta=5$.}
    \label{fig:PC_visit}
\end{figure}

The two-dimensional spatial contours detailed in Figure~\ref{fig:PC_visit} elucidate the complex interplay between the macroscopic piecewise-slope boundary profile and the dominant crack-tip singularity. The displacement fields provide a clear visualization of the structural response: the vertical displacement $u_2$ exhibits a symmetric deflection pattern along the top boundary, directly mapping the non-uniform compressive loading \cite{GhoshBhattaMallikarjunaiah2026}. Furthermore, the $u_2$ field captures the theoretical interpenetration of the crack faces ($y=0.5$), a mathematical consequence of formulating the boundary value problem without explicit Signorini non-penetration constraints. The lateral Poisson expansion, visualized in the horizontal displacement $u_1$, also reflects this global asymmetry, differing significantly from uniform loading profiles.

Despite the asymmetric transmission of these boundary loads through the bulk domain, the intense gradients of the fracture process zone remain decisively localized by the crack-tip geometry. Both the stress and strain fields preserve their characteristic concentrated, lobular patterns immediately surrounding the singularity. This confirms that while the global displacement topology is heavily influenced by the macroscopic loading profile, the extreme local fields governing potential failure are strictly dictated by the interaction between the geometric discontinuity and the density-dependent, strain-limiting constitutive law.

\subsection{Point load at the center}

\begin{figure}[H]
    \centering
    \begin{subfigure}[t]{0.32\linewidth}
        \centering
        \includegraphics[width=\linewidth]{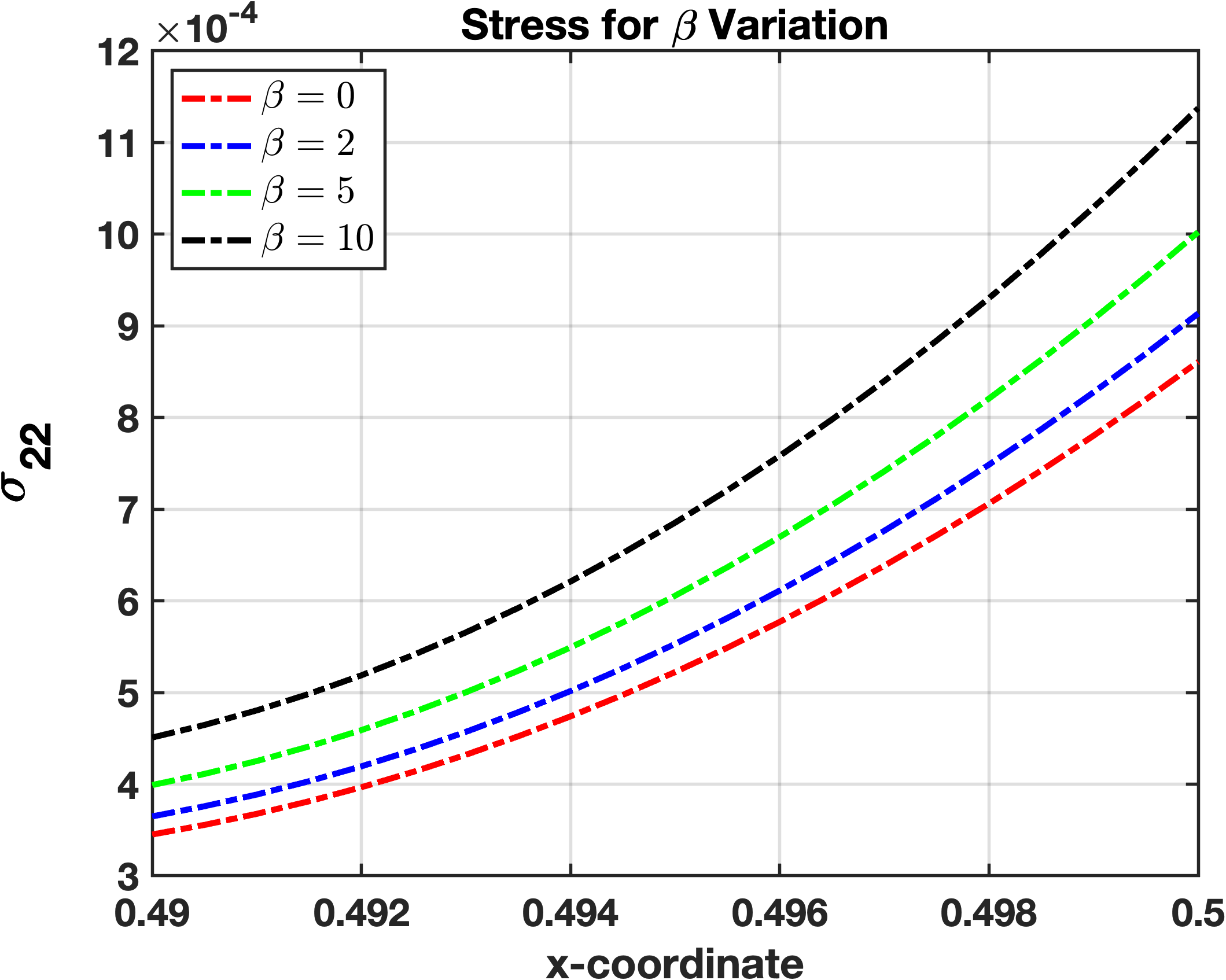}
        \caption{Normal stress $\sigma_{22}$.}
        \label{fig:PL_matlab_stress}
    \end{subfigure}%
    \hfill
    \begin{subfigure}[t]{0.32\linewidth}
        \centering
        \includegraphics[width=\linewidth]{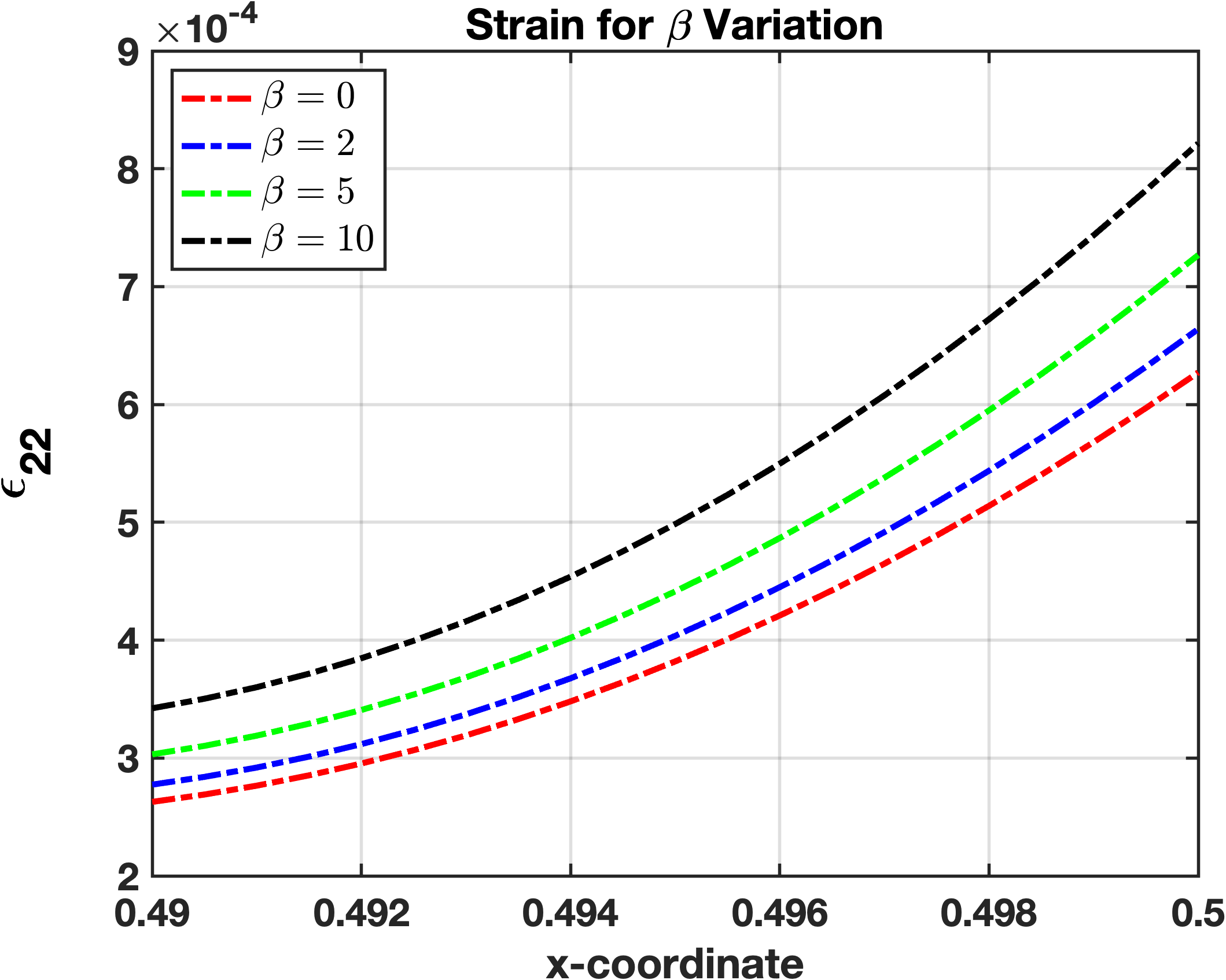}
        \caption{Normal strain $\varepsilon_{22}$.}
        \label{fig:PL_matlab_strain}
    \end{subfigure}%
    \hfill
    \begin{subfigure}[t]{0.32\linewidth}
        \centering
        \includegraphics[width=\linewidth]{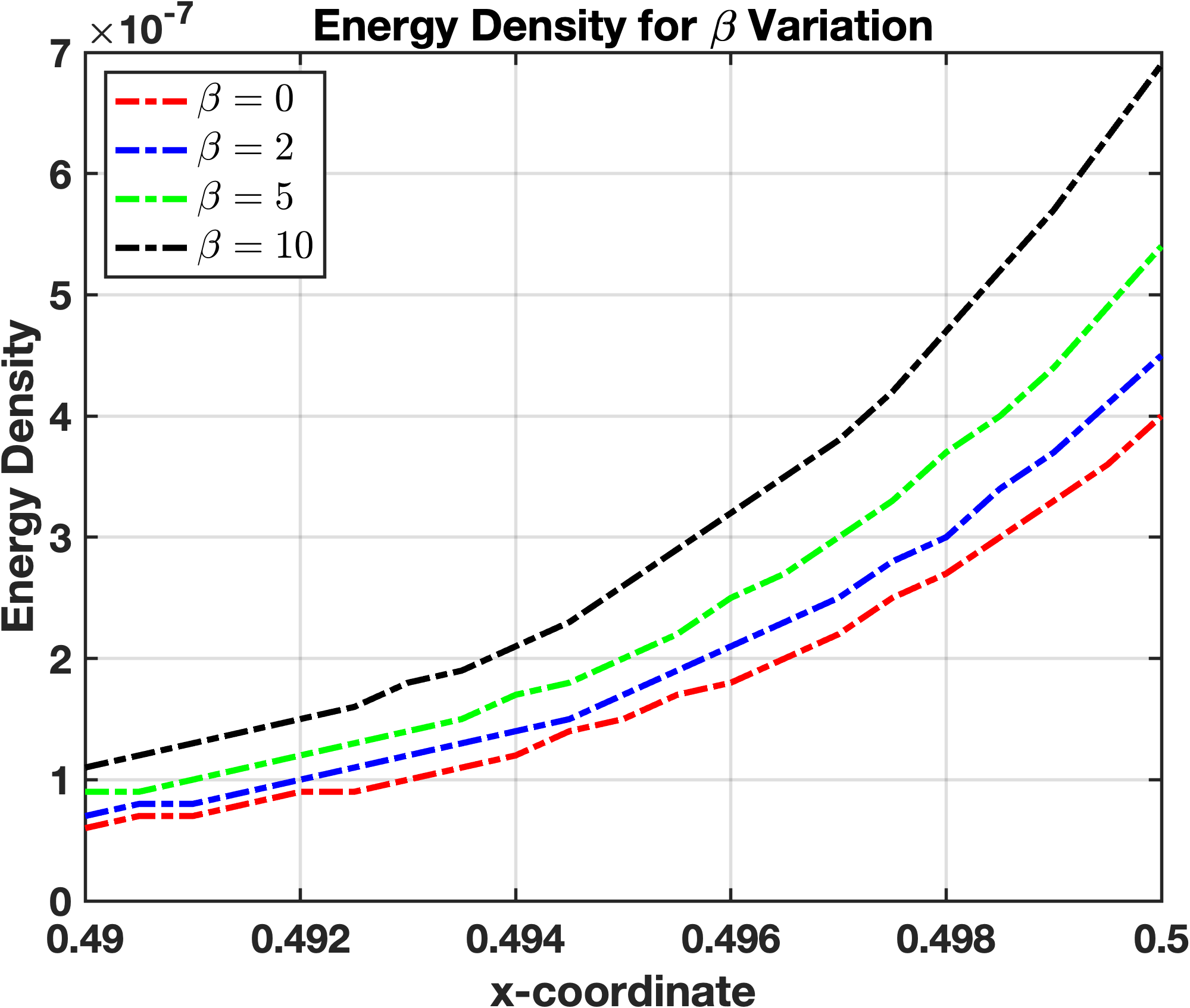}
        \caption{Strain Energy Density.}
        \label{fig:PL_matlab_sed}
    \end{subfigure}
    
    \caption{Evolution of the normal stress $\sigma_{22}$, normal strain $\varepsilon_{22}$, and SED along the sampling line approaching the crack tip. The plots demonstrate the influence of the nonlinearity parameter $\beta$ under a localized central point load.}
    \label{fig:PL_matlab}
\end{figure}

The application of a concentrated point load at the center of the upper boundary introduces a profound mechanical complexity: a dual-singularity domain. The continuum must simultaneously negotiate the severe boundary singularity at the load application point and the internal geometric singularity at the crack tip. As depicted in Figure~\ref{fig:PL_matlab}, the normal stress $\sigma_{22}$, normal strain $\varepsilon_{22}$, and SED strictly increase as the sampling line approaches $x=0.5$. However, unlike the previously analyzed uniform boundary conditions, the influence of the strain-limiting parameter $\beta$ under a point load manifests a highly complex, non-monotone structural response precisely at the crack tip. 

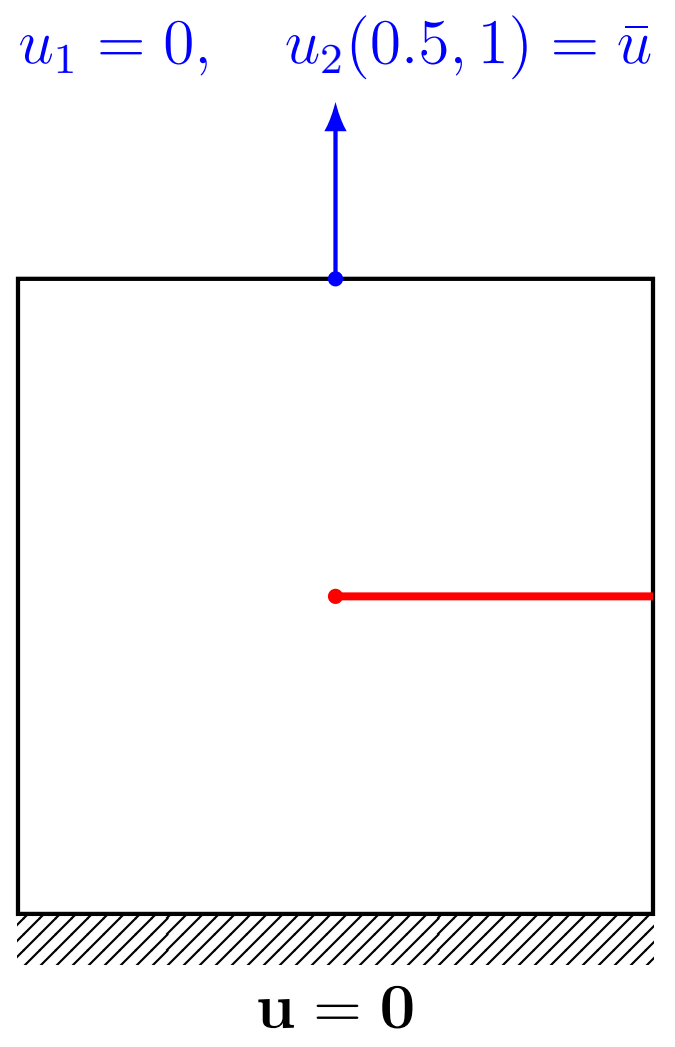
\begin{figure}[H]
\centering
\begin{tikzpicture}[x=4.2cm,y=4.2cm]
  \draw[thick] (0,0) rectangle (1,1);
  \draw[red,line width=1.5pt] (0.5,0.5)--(1,0.5);
  \fill[red] (0.5,0.5) circle[radius=0.012];
  \fill[pattern=north east lines] (0,-0.08) rectangle (1,0);
  \draw[-{Latex[length=2mm]},thick,blue] (0.5,1)--++(0,0.28);
  \fill[blue] (0.5,1) circle[radius=0.012];
  \node[blue,above] at (0.5,1.28) {$u_1=0, \quad u_2(0.5,1)=\bar u$};
  \node[below] at (0.5,-0.08) {$\mathbf u=\mathbf 0$};
\end{tikzpicture}
\caption{Center-point-displacement boundary condition.}
\label{fig:ddm-schematic-v}
\end{figure}

\begin{table}[H]
    \centering
    \begin{tabular}{|c|c|c|}
        \hline
        $\beta$ & $\sigma_{22}$ & $\varepsilon_{22}$ \\
        \hline
        $0$  & 0.02585762 & 0.01833390 \\
        $2$  & 0.02291417 & 0.01540787 \\
        $5$  & 0.02682407 & 0.01719116 \\
        $10$ & 0.03137257 & 0.01627325 \\
        \hline
    \end{tabular}
    \caption{Finite element estimates of the normal stress $\sigma_{22}$ and normal strain $\varepsilon_{22}$ evaluated precisely at the nominal crack-tip location $(0.5,0.5)$ for the concentrated-load case.}
    \label{tab:point_source}
\end{table}

The tabulated finite element estimates in Table~\ref{tab:point_source} explicitly capture this non-monotone interplay. From the linear elastic baseline ($\beta=0$) to a mildly nonlinear state ($\beta=2$), both the local stress and strain decrease. Mechanically, the initial activation of the strain-limiting constitutive law alters the global compliance of the domain; this slight stiffening diffuses the load path originating from the boundary point load, momentarily shielding the crack tip. However, as the nonlinearity intensifies ($\beta=5$ to $\beta=10$), strict kinematic truncation dominates the constitutive response. Because the material is fundamentally restricted from stretching to accommodate the localized energy, the continuum is forced to channel massive internal forces directly between the two singularities to preserve static equilibrium. Consequently, at the highly nonlinear limit ($\beta=10$), the normal stress $\sigma_{22}$ severely eclipses its linear elastic counterpart, while the normal strain $\varepsilon_{22}$ remains firmly bounded below the $\beta=0$ baseline. This stress-dominant state directly correlates with the steep exponential surge in the near-tip Strain Energy Density observed for $\beta=10$.

\section{Conclusion}
\label{sec:conclusion}

In this work, a comprehensive mathematical and computational framework was established to investigate the fracture mechanics of isotropic elastic solids governed by a density-dependent, nonlinear constitutive relation. By extending classical linear elasticity to incorporate volumetric strain-dependent moduli via the relation~\eqref{eq:magnitude_constitutive}, a well-posed boundary value problem was formulated and rigorously analyzed. To resolve the complex, multiscale field gradients induced by the interaction between the material nonlinearity and the crack-tip singularity, a robust $hp$-adaptive finite element strategy was deployed. This high-fidelity framework enabled a systematic evaluation of how diverse macroscopic loading topologies---specifically, uniform tension and compression, piecewise-slope distributions, and localized concentrated loads---dictate the mechanical response near a stationary crack tip.

The numerical investigations reveal that the manifestation of material nonlinearity is profoundly coupled to the spatial distribution and direction of the applied external load. For distributed loading configurations, such as the piecewise-slope cases, the near-tip stress fields exhibit a pronounced sensitivity to the nonlinear constitutive parameter, whereas the corresponding localized strain profiles demonstrate a tightly bounded saturation effect across the scales analyzed. Conversely, under concentrated loading, the continuum exhibits a complex, non-monotone response at the crack tip. This phenomenon arises from the simultaneous, overlapping influence of the load-application singularity and the crack-front localization. These comparative results mathematically demonstrate that the influence of the density-dependent constitutive parameter cannot be treated as an intrinsic, isolated scalar effect; rather, its impact is highly dynamic and strictly modulated by the global loading topology.

While the present study provides critical insights into stationary crack-tip fields within an infinitesimal-strain regime, it also delineates natural boundaries for future theoretical extensions. The current variational formulation assumes a traction-free, fixed crack geometry and does not mathematically preclude crack-face interpenetration. Consequently, future efforts must reformulate the boundary value problem as a variational inequality to rigorously enforce Signorini contact and friction conditions, particularly under compressive or mixed-mode shear loading. Furthermore, because the local fields reported herein are discrete numerical approximations of an ideal geometric crack, extending this nonlinear constitutive framework to incorporate phase-field damage regularization or quasi-static crack propagation criteria will be vital for predicting macroscopic structural failure.

Ultimately, these findings underscore a fundamental paradigm shift in fracture modeling utilizing algebraically nonlinear constitutive relationships: accurately evaluating localized mechanical fields in density-dependent elastic media requires a holistic analytical approach that simultaneously resolves the constitutive nonlinearity, the localized crack geometry, and the exact spatial distribution of the applied boundary loads.

\section*{Acknowledgments}
SG thanks the University of Texas Rio Grande Valley for providing a Presidential Research Fellowship during his PhD studies. SMM gratefully acknowledges the financial support of the National Science Foundation under Grant No.~2316905.

\bibliographystyle{unsrt}
\bibliography{references}

@book{brezis2011functional,
  title={Functional analysis, Sobolev spaces and partial differential equations},
  author={Brezis, Haim},
  volume={2},
  number={3},
  year={2011},
  publisher={Springer}
}

@book{showalter2013monotone,
  title={Monotone operators in Banach space and nonlinear partial differential equations},
  author={Showalter, Ralph Edwin},
  volume={49},
  year={2013},
  publisher={American Mathematical Soc.}
}

@article{itou2023generalization,
  title={A generalization of the Kelvin--Voigt model with pressure-dependent moduli in which both stress and strain appear linearly},
  author={Itou, Hiromichi and Kovtunenko, Victor A and Rajagopal, K. R.},
  journal={Mathematical Methods in the Applied Sciences},
  volume={46},
  number={14},
  pages={15641--15654},
  year={2023},
  publisher={Wiley Online Library}
}

@article{itou2022investigation,
  title={Investigation of implicit constitutive relations in which both the stress and strain appear linearly, adjacent to non-penetrating cracks},
  author={Itou, Hiromichi and Kovtunenko, Victor A and Rajagopal, K. R.},
  journal={Mathematical Models and Methods in Applied Sciences},
  volume={32},
  number={07},
  pages={1475--1492},
  year={2022},
  publisher={World Scientific}
}

@article{itou2019well,
  title={Well-posedness of the problem of non-penetrating cracks in elastic bodies whose material moduli depend on the mean normal stress},
  author={Itou, Hiromichi and Kovtunenko, Victor A and Rajagopal, K. R.},
  journal={International Journal of Engineering Science},
  volume={136},
  pages={17--25},
  year={2019},
  publisher={Elsevier}
}

@article{Rajagopal2003,
  author  = {Rajagopal, K. R.},
  title   = {On implicit constitutive theories},
  journal = {Applications of Mathematics},
  year    = {2003},
  volume  = {48},
  number  = {4},
  pages   = {279--319},
  doi     = {10.1023/A:1026062615145}
}

@article{Rajagopal2007,
  author  = {Rajagopal, K. R.},
  title   = {The elasticity of elasticity},
  journal = {Zeitschrift f{\"u}r Angewandte Mathematik und Physik},
  year    = {2007},
  volume  = {58},
  number  = {2},
  pages   = {309--317},
  doi     = {10.1007/s00033-006-6084-5}
}

@article{RajagopalSrinivasa2007,
  author  = {Rajagopal, K. R. and Srinivasa, A. R.},
  title   = {On the response of non-dissipative solids},
  journal = {Proceedings of the Royal Society A: Mathematical, Physical and Engineering Sciences},
  year    = {2007},
  volume  = {463},
  number  = {2078},
  pages   = {357--367},
  doi     = {10.1098/rspa.2006.1760}
}

@article{RajagopalSrinivasa2009,
  author  = {Rajagopal, K. R. and Srinivasa, A. R.},
  title   = {On a class of non-dissipative materials that are not hyperelastic},
  journal = {Proceedings of the Royal Society A: Mathematical, Physical and Engineering Sciences},
  year    = {2009},
  volume  = {465},
  number  = {2102},
  pages   = {493--500},
  doi     = {10.1098/rspa.2008.0319}
}

@article{Rajagopal2010,
  author  = {Rajagopal, K. R.},
  title   = {On a new class of models in elasticity},
  journal = {Mathematical and Computational Applications},
  year    = {2010},
  volume  = {15},
  number  = {4},
  pages   = {506--528},
  doi     = {10.3390/mca15040506}
}

@article{Rajagopal2011Limiting,
  author  = {Rajagopal, K. R.},
  title   = {Non-linear elastic bodies exhibiting limiting small strain},
  journal = {Mathematics and Mechanics of Solids},
  year    = {2011},
  volume  = {16},
  number  = {1},
  pages   = {122--139},
  doi     = {10.1177/1081286509357272}
}

@article{Rajagopal2014,
  author  = {Rajagopal, K. R.},
  title   = {On the nonlinear elastic response of bodies in the small strain range},
  journal = {Acta Mechanica},
  year    = {2014},
  volume  = {225},
  number  = {6},
  pages   = {1545--1553},
  doi     = {10.1007/s00707-013-1015-y}
}

@article{RajagopalSaccomandi2016,
  author  = {Rajagopal, K. R. and Saccomandi, Giuseppe},
  title   = {A novel approach to the description of constitutive relations},
  journal = {Frontiers in Materials},
  year    = {2016},
  volume  = {3},
  pages   = {36},
  doi     = {10.3389/fmats.2016.00036}
}

@article{Rajagopal2018Linearization,
  author  = {Rajagopal, K. R.},
  title   = {A note on the linearization of the constitutive relations of non-linear elastic bodies},
  journal = {Mechanics Research Communications},
  year    = {2018},
  volume  = {93},
  pages   = {132--137},
  doi     = {10.1016/j.mechrescom.2017.08.002}
}

@article{Rajagopal2021,
  author  = {Rajagopal, K. R.},
  title   = {An implicit constitutive relation for describing the small strain response of porous elastic solids whose material moduli are dependent on the density},
  journal = {Mathematics and Mechanics of Solids},
  year    = {2021},
  volume  = {26},
  number  = {8},
  pages   = {1138--1146},
  doi     = {10.1177/10812865211021465}
}

@article{RajagopalSaccomandi2022,
  author  = {Rajagopal, K. R. and Saccomandi, Giuseppe},
  title   = {Implicit nonlinear elastic bodies with density dependent material moduli and its linearization},
  journal = {International Journal of Solids and Structures},
  year    = {2022},
  volume  = {234--235},
  pages   = {111255},
  doi     = {10.1016/j.ijsolstr.2021.111255}
}

@article{AlagappanArumugamRajagopal2023,
  author  = {Alagappan, P. and Arumugam, J. and Rajagopal, K. R.},
  title   = {A note on the response of elastic bodies whose material moduli depend on the density and the mechanical pressure},
  journal = {Applications in Engineering Science},
  year    = {2023},
  volume  = {16},
  pages   = {100162},
  doi     = {10.1016/j.apples.2023.100162}
}

@article{ItouKovtunenkoRajagopal2021,
  author  = {Itou, Hiromichi and Kovtunenko, Victor A. and Rajagopal, K. R.},
  title   = {On an implicit model linear in both stress and strain to describe the response of porous solids},
  journal = {Journal of Elasticity},
  year    = {2021},
  volume  = {144},
  number  = {1},
  pages   = {107--118},
  doi     = {10.1007/s10659-021-09831-x}
}

@article{RajagopalWalton2011,
  author  = {Rajagopal, K. R. and Walton, Jay R.},
  title   = {Modeling fracture in the context of a strain-limiting theory of elasticity: A single anti-plane shear crack},
  journal = {International Journal of Fracture},
  year    = {2011},
  volume  = {169},
  number  = {1},
  pages   = {39--48},
  doi     = {10.1007/s10704-010-9581-7}
}

@article{KulvaitMalekRajagopal2013,
  author  = {Kulvait, Vojt{\v e}ch and M{\'a}lek, Josef and Rajagopal, K. R.},
  title   = {Anti-plane stress state of a plate with a {V}-notch for a new class of elastic solids},
  journal = {International Journal of Fracture},
  year    = {2013},
  volume  = {179},
  number  = {1--2},
  pages   = {59--73},
  doi     = {10.1007/s10704-012-9772-5}
}

@article{BulicekMalekRajagopalWalton2015,
  author  = {Bul{\'i}{\v c}ek, Miroslav and M{\'a}lek, Josef and Rajagopal, K. R. and Walton, Jay R.},
  title   = {Existence of solutions for the anti-plane stress for a new class of ``strain-limiting'' elastic bodies},
  journal = {Calculus of Variations and Partial Differential Equations},
  year    = {2015},
  volume  = {54},
  number  = {2},
  pages   = {2115--2147},
  doi     = {10.1007/s00526-015-0859-5}
}

@article{GouMallikarjunaiahRajagopalWalton2015,
  author  = {Gou, K. and Mallikarjunaiah, S. M. and Rajagopal, K. R. and Walton, Jay R.},
  title   = {Modeling fracture in the context of a strain-limiting theory of elasticity: A single plane-strain crack},
  journal = {International Journal of Engineering Science},
  year    = {2015},
  volume  = {88},
  pages   = {73--82},
  doi     = {10.1016/j.ijengsci.2014.04.018}
}

@article{MallikarjunaiahWalton2015,
  author  = {Mallikarjunaiah, S. M. and Walton, Jay R.},
  title   = {On the direct numerical simulation of plane-strain fracture in a class of strain-limiting anisotropic elastic bodies},
  journal = {International Journal of Fracture},
  year    = {2015},
  volume  = {192},
  number  = {2},
  pages   = {217--232},
  doi     = {10.1007/s10704-015-0006-5}
}

@article{MaiWalton2015,
  author  = {Mai, Tina and Walton, Jay R.},
  title   = {On strong ellipticity for implicit and strain-limiting theories of elasticity},
  journal = {Mathematics and Mechanics of Solids},
  year    = {2015},
  volume  = {20},
  number  = {2},
  pages   = {121--139},
  doi     = {10.1177/1081286514544254}
}

@article{BonitoGiraultSuli2020,
  author  = {Bonito, Andrea and Girault, Vivette and S{\"u}li, Endre},
  title   = {Finite element approximation of a strain-limiting elastic model},
  journal = {IMA Journal of Numerical Analysis},
  year    = {2020},
  volume  = {40},
  number  = {1},
  pages   = {29--86},
  doi     = {10.1093/imanum/dry065}
}

@article{YoonLeeMallikarjunaiah2021,
  author  = {Yoon, Hyun C. and Lee, Sanghyun and Mallikarjunaiah, S. M.},
  title   = {Quasi-static anti-plane shear crack propagation in nonlinear strain-limiting elastic solids using phase-field approach},
  journal = {International Journal of Fracture},
  year    = {2021},
  volume  = {227},
  number  = {2},
  pages   = {153--172},
  doi     = {10.1007/s10704-020-00501-y}
}

@article{YoonMallikarjunaiah2022,
  author  = {Yoon, Hyun C. and Mallikarjunaiah, S. M.},
  title   = {A finite-element discretization of some boundary value problems for nonlinear strain-limiting elastic bodies},
  journal = {Mathematics and Mechanics of Solids},
  year    = {2022},
  volume  = {27},
  number  = {2},
  pages   = {281--307},
  doi     = {10.1177/10812865211020789}
}

@article{LeeYoonMallikarjunaiah2022,
  author  = {Lee, Sanghyun and Yoon, Hyun Chul and Mallikarjunaiah, S. M.},
  title   = {Finite element simulation of quasi-static tensile fracture in nonlinear strain-limiting solids with the phase-field approach},
  journal = {Journal of Computational and Applied Mathematics},
  year    = {2022},
  volume  = {399},
  pages   = {113715},
  doi     = {10.1016/j.cam.2021.113715}
}

@article{YoonVasudevaMallikarjunaiah2022,
  author  = {Yoon, Hyun Chul and Vasudeva, Karthik K. and Mallikarjunaiah, S. M.},
  title   = {Finite element model for a coupled thermo-mechanical system in nonlinear strain-limiting thermoelastic body},
  journal = {Communications in Nonlinear Science and Numerical Simulation},
  year    = {2022},
  volume  = {108},
  pages   = {106262},
  doi     = {10.1016/j.cnsns.2022.106262}
}

@article{YoonMallikarjunaiahBhatta2024,
  author  = {Yoon, Hyun C. and Mallikarjunaiah, S. M. and Bhatta, Dambaru},
  title   = {Finite element solution of crack-tip fields for an elastic porous solid with density-dependent material moduli and preferential stiffness},
  journal = {Advances in Mechanical Engineering},
  year    = {2024},
  volume  = {16},
  number  = {2},
  pages   = {16878132241231792},
  doi     = {10.1177/16878132241231792}
}

@article{GhoshBhattaMallikarjunaiah2026,
  author  = {Ghosh, Saugata and Bhatta, Dambaru and Mallikarjunaiah, S. M.},
  title   = {Computational insights into orthotropic fracture: crack-tip fields in strain-limiting materials under non-uniform loads},
  journal = {Zeitschrift f{\"u}r Angewandte Mathematik und Physik},
  year    = {2026},
  volume  = {77},
  pages   = {34},
  doi     = {10.1007/s00033-025-02672-x}
}

@article{BhattaGhoshMallikarjunaiah2026,
  author  = {Bhatta, Dambaru and Ghosh, Saugata and Mallikarjunaiah, S. M.},
  title   = {Computational investigation of crack-tip fields in a compressed nonlinear strain-limiting material},
  journal = {Mathematics and Mechanics of Solids},
  year    = {2026},
  note    = {OnlineFirst},
  doi     = {10.1177/10812865261417686}
}

@article{Mallikarjunaiah2026,
  author  = {Mallikarjunaiah, S. M.},
  title   = {An {$hp$}-adaptive finite element framework for static cracks: the impact of pointwise density variations on mode {I}, mode {II}, and mixed-mode fracture},
  journal = {Zeitschrift f{\"u}r Angewandte Mathematik und Physik},
  year    = {2026},
  volume  = {77},
  pages   = {203},
  doi     = {10.1007/s00033-026-02850-5}
}

@article{KellyGagoZienkiewiczBabuska1983,
  author  = {Kelly, D. W. and Gago, J. P. de S. R. and Zienkiewicz, O. C. and Babu{\v s}ka, I.},
  title   = {A posteriori error analysis and adaptive processes in the finite element method: Part I---error analysis},
  journal = {International Journal for Numerical Methods in Engineering},
  year    = {1983},
  volume  = {19},
  number  = {11},
  pages   = {1593--1619},
  doi     = {10.1002/nme.1620191103}
}

@article{ArndtBangerthEtAl2021,
  author  = {Arndt, Daniel and Bangerth, Wolfgang and Davydov, Denis and Heister, Timo and Heltai, Luca and Kronbichler, Martin and Maier, Matthias and Pelteret, Jean-Paul and Turcksin, Bruno and Wells, David},
  title   = {The {deal.II} finite element library: Design, features, and insights},
  journal = {Computers \& Mathematics with Applications},
  year    = {2021},
  volume  = {81},
  pages   = {407--422},
  doi     = {10.1016/j.camwa.2020.02.022}
}

\end{document}